\documentclass[reqno,11pt,a4paper]{article}
\usepackage{verbatim}   % useful for program listings
\usepackage{color}      % use if color is used in text
\usepackage{subfigure}  % use for side-by-side figures
\usepackage{hyperref}   % use for hypertext links, including those to external documents and URLs
\usepackage{amsmath}
\usepackage{algorithm}
\usepackage{setspace}
\usepackage[noend]{algpseudocode}
\usepackage{mathtools}
\usepackage{amsfonts}
\usepackage{amssymb}
\usepackage{natbib}
\usepackage{booktabs}
\usepackage{multirow}
\usepackage[pdftex,dvipsnames]{xcolor}  % Coloured text etc.
\usepackage{xargs}                      % Use more than one optional parameter in a new commands
\usepackage{cleveref}
\usepackage[pdftex]{graphicx}
\usepackage[colorinlistoftodos,prependcaption,textsize=tiny]{todonotes}
\usepackage[body={16cm,22cm},centering]{geometry}
\usepackage{amsthm,amsfonts,mathrsfs,graphicx}
\usepackage[font=small,labelfont=bf]{caption}

\newcommand{\eqnum}{\refstepcounter{equation}\textup{\tagform@{\theequation}}}

\newcommandx{\info}[2][1=]{\todo[linecolor=OliveGreen,backgroundcolor=OliveGreen!25,bordercolor=OliveGreen,#1]{#2}}

\algnewcommand\algorithmicinput{\textbf{Input:}}
\algnewcommand\INPUT{\item[\algorithmicinput]}

\DeclareMathOperator{\spann}{span}

\def\R{\mathbb{R}}

\DeclareMathOperator{\prox}{Prox}
\DeclareMathOperator{\sgn}{sgn}

\def\be{\begin{equation}}
\def\ee{\end{equation}}
\def\bea{\begin{eqnarray}}
\def\eea{\end{eqnarray}}
\def\bsec{\begin{subequations}\begin{eqnarray} }
\def\esec{\end{eqnarray} \end{subequations} }
\def\veci{{\mathbf i}}

\def\vu{{\mathbf u}}
\def\vecy{{\mathbf y}}

\def\vv{{\mathbf v}}

\def\vp{{\mathbf p}}
\def\vq{{\mathbf q}}

\newcommand{\bw}{\mathbf{w}}
\newcommand{\bx}{\mathbf{x}}
\newcommand{\by}{\mathbf{y}}
\newcommand{\bu}{\mathbf{u}}
\newcommand{\bff}{\mathbf{f}}

\newcommand{\bg}{\mathbf{g}}
\newcommand{\basis}{\mathfrak{I}}

\title{Optimal Feedback Law Recovery by Gradient-Augmented Sparse Polynomial Regression}
\author{Behzad Azmi\thanks{RICAM, Austrian Academy of Sciences, Altenbergerstra\ss e 69, A-4040 Linz, Austria. E-mail: behzad.azmi@ricam.oeaw.ac.at}\and
Dante Kalise\thanks{School of Mathematical Sciences, University of Nottingham, NG7 2QL, United Kingdom. E-mail: dante.kalise@nottingham.ac.uk}\and
Karl Kunisch\thanks{RICAM, Austrian Academy of Sciences and Department of Mathematics, University of Graz, Heinrichstra\ss e 36, A-8010 Graz, Austria. E-mail:karl.kunisch@uni-graz.at}}
\date{}                     %% if you don't need date to appear

\begin{document}
\date{\today}
\maketitle

\begin{abstract}
	A {\color{black}sparse regression} approach for the computation of high-dimensional optimal feedback laws arising in deterministic nonlinear control is proposed. The approach exploits the control-theoretical link between Hamilton-Jacobi-Bellman PDEs characterizing the value function of the optimal control problems, and first-order optimality conditions via Pontryagin's Maximum Principle. The latter is used as a representation formula to recover the value function and its gradient at arbitrary points in the space-time domain through the solution of a two-point boundary value problem. After generating a dataset consisting of different state-value pairs, a hyperbolic cross polynomial model for the value function is fitted using a LASSO regression. An extended set of low and high-dimensional numerical tests in nonlinear optimal control reveal that enriching the dataset with gradient information reduces the number of training samples, and that the sparse polynomial regression consistently yields a feedback law of lower complexity. 
\end{abstract}
%\begin{keywords}
%	Optimal Feedback Control, Optimality Conditions, Hamilton-Jacobi-Bellman PDE, Polynomial Approximation, Sparse Optimization
%\end{keywords}

\section{Introduction}\label{intro}

A large class of design problems including the synthesis of autopilot and guidance systems, {\color{black}the stabilization of chemical reactions}, and the control of fluid flow phenomena, among others, can be cast as {\color{black}deterministic nonlinear} optimal control problems. {\color{black} In this framework, we synthesize a time-dependent control signal $\bu(t)$ in $\mathbb{R}^m$ by solving a dynamic optimization problem which we formulate as
	\begin{align}\label{eq:optcost}
		\underset{u(\cdot)\in L^2(t_0,T;\R^m)}{\min}J(\bu;t_0,\bx) := \int^{\top}_{t_0} \ell(\by(t))+\beta\|\bu(t))\|_2^2\,dt\,,\qquad \beta>0\,,
	\end{align}
	subject to $\by(t)$ in $\mathbb{R}^n $  being the solution to the control-affine nonlinear dynamics
	\begin{align}\label{eq:state}
		\frac{d}{dt}\vecy(t)&=\bff(\vecy(t))+\bg(\by(t))\vu(t)\,,\qquad \vecy(t_0)=\bx\,.
	\end{align}
	We assume that the running cost $\ell : \mathbb{R}^n \to \mathbb{R}$, the dynamics $\bff: \mathbb{R}^n \to \mathbb{R}^n$, and  $\bg: \mathbb{R}^n \to \mathbb{R}^{n \times m}$, are continuously differentiable.}
The numerical realization of control laws by solving the dynamic optimization problem above is a topic at the interface between control theory, computational optimization, and numerical analysis. While this problem dates back to the birth of Calculus of Variations, it was during the second half of the 20th century when two major methodological breakthroughs shaped our understanding of optimal control theory, namely, the development of Pontryagin's Maximum Principle and the theory of Dynamic Programming (see ~\citep{Peschhist} for a historical survey on this topic). On the one hand, Pontryagin's Maximum Principle (PMP) ~\citep{PMP} yields first-order optimality conditions for \eqref{eq:optcost}-\eqref{eq:state} in the form of a two-point boundary value problem for a forward-backward coupling between optimal state,{\color{black} adjoint $\vp^*=(p_1^*,\ldots,p_n^*)$, and control variables, denoted by $(\vecy^*(t),\vp^*(t),\vu^*(t))$ respectively, which in short reads
	\begin{align}\tag{TPBVP}\label{eq:tpbvp}
		\begin{cases}
			\frac{d}{dt}\vecy^*(t)=\bff(\vecy^*(t))+\bg(\by^*(t))\vu^*(t)\,,\qquad \\
			\vecy^*(t_0)=\bx\,,\\
			-\frac{d}{dt}p_i^*(t)=\sum\limits_{j=1}^np^*_j(t)\left(\partial_{y_i}(f_j(\by^*(t))+g_j(\by^*(t))\bu^*(t)\right)+\partial_{y_i}\ell(\by^*(t))\,,\qquad i=1,\ldots,n\,,\\
			p_i^*(T)=0\,,
		\end{cases}
	\end{align}
	closed with the optimality condition
	\begin{align}\label{eq:optcon}
		\bu^*(t)=-\frac{1}{2\beta}\bg^{\top}(\by^*(t))\vp^*(t)\,,\qquad\forall t\in(t_0,T)\,.
	\end{align}
	This procedure yields an optimal state-adjoint-control triple originating from the initial condition $\bx$, in what is known as an open-loop control. On the other hand, the Dynamic Programming approach synthesizes the optimal control as 
	\begin{equation}\label{eq:feedhjb}
		\bu^*(t,\bx) = \underset{u\in R^m}{\operatorname{argmin}}\left\{\beta \|\bu\|_2^2  + \nabla V(t,\bx)^{\top}\left(\bg(\bx)\bu\right)\right\} = - \frac{1}{2\beta}\bg^{\top}(\bx)\nabla V(t,\bx)\,,
	\end{equation}
	where $V(t,\bx): [0,T]\times \R^n\longrightarrow \R$ is the value function of the problem
	\begin{equation}
		V(t,\bx):= \inf_{\bu(\cdot)} \{J(\bu;t,\bx) \text{ subject to \eqref{eq:state}} \}\,,
	\end{equation}
	which in turn satisfies a first-order, nonlinear Hamilton-Jacobi-Bellman (HJB) partial differential equation of the form
	\begin{align}\label{HJBe}\tag{HJB}
		\begin{cases}
			\partial_t V(t,\bx)-\frac{1}{2\beta}\nabla V(t,\bx)^{\top}\bg(\bx)\bg^{\top}(\bx)\nabla V(t,\bx)+\nabla V(t,\bx)^{\top}\bff(\bx)+\ell(\bx)=0\,,\\
			V(T,\bx)=0\,,
		\end{cases}
\end{align}}
to be solved over the state space of the dynamics $\mathcal{X}\subset\R^n$ ~\citep[Chapter 1]{BCD}. Here and throughout we assume that $V$ is $C^1$, see e.g.\cite{fran2002}. This approach expresses the optimal control as a feedback map or closed-loop form, i.e. $\vu^*=\vu^*(t,\by(t))$, yielding a control law that is optimal in the whole state space. From a practical viewpoint, optimal trajectories obtained from the solution of \eqref{eq:tpbvp} are not robust with respect to disturbances, {\color{black}justifying the necessity} of a HJB-based feedback synthesis. However, the solution of HJB PDEs, especially for high-dimensional dynamical systems, comes at a formidable computational cost, often referred in the literature as the \textsl{curse of dimensionality} ~\citep{curse}. Over the last years, a number of works have reported remarkable progress in the solution of high-dimensional HJB PDEs, including the use of sparse grids~\citep{axelsparse,bokasparse}, tree structure algorithms~\citep{allatree}, max-plus methods ~\citep{maxplus1,maxplus2}, polynomial approximation ~\citep{poli1,poli2}, tensor decomposition techniques~\citep{tensor1,tensor2,tensor3,tensor4,tensor5}, and an evergrowing literature on artificial neural networks~\citep{ml1,ml2,ml3,ml4,ml5}.

\paragraph{On the link between PMP and HJB PDEs.} In this work, we follow an alternative approach that circumvents the curse of dimensionality in the computation of optimal feedback laws by exploiting the relation between the HJB PDE and the PMP. We interpret system (TPBVP) as a representation formula for the solution of \eqref{HJBe} . There exists an extensive literature dating back to ~\citep[Chapter 1]{PMP} discussing the relation between dynamic programming and first-order optimality conditions. In the simplest version of the statement, assuming the solution of the HJB PDE is $\mathcal{C}^2$, it can be shown that the forward-backward dynamics originating from the PMP correspond to the characteristic curves of the HJB equation.{\color{black} This result was further improved in ~\citep{BJ86}, where the PMP was derived from the viscosity solution of a first-order HJB PDE. The precise result linking the solution of the PMP with the characteristic curves of the HJB PDE, and the identification of the adjoint variable as the gradient of the value function can be found in ~\citep[Theorems II.9 and II.10]{Subbotina}. We refer the reader to ~\citep[Section 3.4]{BCD} for an exhaustive revision of the main results in this topic. More concretely, the computation of a given $V(t^i,\bx^i)$ can be realized by solving (TPBVP) setting $t_0=t_i$ and the initial condition $\by(t_0)=\bx^i$, and evaluating the optimal cost \eqref{eq:optcost} using the optimal triple $(\by^*(t),\vp^*(t),\bu^*(t))$. Moreover, the optimal adjoint verifies $\vp^*(t)=\nabla V(t,\by^*(t))$ ~\citep{cv87,cf91,CristianiPMP}. The PMP can be thus understood as a representation formula for the value function and its gradient along an optimal trajectory.}

\paragraph{A data-driven method for computing optimal feedback laws.} We restrict our attention to a class of smooth and unconstrained nonlinear optimal control problems where the aforedescribed link between PMP and the HJB PDE is direct, and we use it to generate a characteristic-based, causality-free method to approximate $V(t,\bx)$ and as a by-product $\vu(t,\bx)$, without solving \eqref{HJBe}. To do this, we sample a set of initial conditions $\{(t^i,\bx^i)\}_{i=1}^{N}$, for which we compute both $V(t^i,\bx^i)$ and $\nabla V(t^i,\bx^i)$ by realizing the optimal trajectory through PMP. This is done by following a reduced gradient approach ~\citep{rgkk10}, in which forward-backward iterative solves of \eqref{eq:tpbvp} are combined with a gradient descent method to find the minimizer of $J(\bu;t_0,\bx)$. Having collected a dataset  $\{t^i,\bx^i,V(t^i,\bx^i),\nabla V(t^i,\bx^i)\}_{i=1}^{N}$ enriched with gradient information, we fit a polynomial model for the value function
\begin{align}\label{eq:polans}
	V_{\theta}(t,\bx)=\sum\limits_{i=1}^q\theta_i\Phi_i(t,\bx)=\langle\theta,\Phi\rangle\,,
\end{align} 
with $\Phi(t,\bx)=(\Phi_1(t,\bx),\ldots,\Phi_q(t,\bx))$  are elements of a suitable polynomial basis, and the parameters $\theta=(\theta_1,\ldots,\theta_q)$ obtained from a LASSO regression
\begin{align}\label{LASSO}
	\min_{\theta \in \mathbb{R}^q }\|  [\Phi;\nabla\Phi]\,\theta-[V;\nabla V] \|^2_2 + \lambda \| \theta\|_{1, \mathbf{w}}\,,
\end{align}
where the matrix $[\Phi;\nabla\Phi]\in\R^{(n+1)N\times q}$ and the vector $[V;\nabla V]\in\R^{(n+1)N}$ include value function and gradient data. Finally, the optimal feedback map is recovered as
\begin{align}\label{sparsefeed}
	\vu^*(t,\bx)=\underset{u}{\operatorname{argmin}}\;\;\{\beta \|\bu\|_2^2  + \nabla V_{\theta}(t,\bx)^{\top}\left(\bg(\bx)\bu\right)\}\,.
\end{align}
{\color{black} Recall from the discussion above that the gradient $\nabla V(\bx)$ can be obtained at no extra cost from open-loop evaluations of $V(\bx)$.}

\paragraph{Related literature and contributions.} The idea of using open-loop solves to build a numerical representation of an optimal feedback law dates back at least to ~\citep{beeler,itointerp}, where a low-dimensional feedback law was constructed directly by interpolating the optimal control from a set of collocation points. More recently, this idea has been revisited and improved in ~\citep{kang1,kang2} by exploiting the extremely parallelizable structure of the open-loop solves together with a sparse grid interpolant, scaling up to 6-dimensional examples. In ~\citep{nakazim,nakazim2} this idea has been further developed by replacing the use of a grid-based interpolant with artificial neural networks which are trained using a dataset consisting of both the value function and gradient evaluations, presenting computational results up to dimension 30. In a similar vein, the works ~\citep{chow1,chow2,chow3,chow4} have proposed the use of representation formulas for HJB PDEs, ranging from the celebrated Lax-Hopf formula to variations of the PMP, in conjunction with efficient convex optimization techniques for the solution of point values of the HJB PDE on the fly. These results have been further studied in ~\citep{Yegorov_2018}, and used in ~\citep{Yegorov_2019,lars2020} for the construction of control Lyapunov functions both with sparse grids and neural networks.
Our work, while in line with the aforedescribed, proposes a different methodology which is summarized in the following ingredients:
\begin{itemize}
	\item {\color{black} An enriched dataset containing both value function and gradient information. Approximating a function based on measurements of both the function values and its derivatives dates back to Hermite and spline interpolants. Here, we follow the approach presented in ~\citep{MR3975885}, where gradient information is used in a sparse regression framework, showing that gradient-augmented sparse regression can reduce the amount of samples required to reach a certain training error. In our case, samples are obtained from the solution of open-loop solves that are realized through a reduced gradient approach. The reduction of the number of samples due to the inclusion of gradient information is particularly relevant for high-dimensional nonlinear optimal control problems, as sampling generation can be particularly costly.} 
	\item The value function, and as a consequence the feedback law, is approximated with a polynomial ansatz \eqref{eq:polans}. This choice is backed by the extensive literature concerning power series approximations of the value function ~\citep{Alb61,KreAH13,Luk69,BreKP19}. In fact, it is well-known that for linear dynamics and quadratic cost functions, the value function corresponds to a quadratic form, which is contained in the span of our approximation space. In ~\citep{poli1,poli2}, we have studied a Galerkin approach for HJB PDEs arising in nonlinear control and games with polynomial approximation functions. In these works we have constructed a polynomial basis limited by the total degree of the monomials, solving up to 14-dimensional tests. Here instead, borrowing a leaf from the vast literature on polynomial approximation theory ~\citep{poliap1,poliap2,poliap3,poliap5}, we consider a hierarchical basis defined through hyperbolic cross approximation, for which we report tests up to dimension 80 at moderate computational cost. {\color{black} Following results such as \citep[Theorem 4.1]{MR3975885}, the sparse regression of a feedback law with gradient-augmented information and a hyperbolic cross expansion has an error in the $H^1$-norm that is governed by the projection of the original value function over the hyperbolic cross set\footnote{assuming a sufficiently large training set, which we generate offline.} In general, given an optimal control problem, it is difficult to establish a priori bounds for this projection error. This is mainly due to the fact that the regularity of $V$ as a function of $f, g$ and $\ell$ is difficult to quantify. Of course, for the linear-quadratic control problem the value function is a quadratic function and therefore it can be accurately represented by expansion in the hyperbolic cross with low polynomial degree. Sufficient conditions for local $C{2,1}$-regularity  of $V$ are given in \cite{canfran2013}. Higher order regularity will depend on the problem structure, see e.g. \cite{BreKP19}. As we restrict our problem formulation to control-affine dynamics without control or state constraint and smooth value functions, circumventing the fully nonlinear case, we expect that the value function can be correctly captured by quadratic and high-order polynomial terms.}
	\item The use of a polynomial expansion for the value function, which is linear in the coefficients, allows us to fit the value function model through a LASSO regression framework \eqref{LASSO}. This least squares problem can easily account for the use of gradient data and has a well-understood numerical realization ~\citep{MR3975885, Laurent2019}. Moreover, we include an $\ell_1$ penalty on the expansion coefficients, leading to the synthesis of low complexity feedback laws \eqref{sparsefeed}, which is crucial for a fast online computation of the feedback action. 
	\item Our extensive numerical assessment of the methodology includes a class of genuinely non-linear, fully coupled, high-dimensional dynamics arising in agent-based modelling ~\citep{mfc1}, which ultimately connects with the control of non-local transport equations arising in mean field control ~\citep{mfc2,FS14,SGPavl18}.
\end{itemize}

\paragraph{Stabilization with static feedback laws.} For large optimization horizon T,  and $\ell(\by)=\|\by\|^2$, the cost \eqref{eq:optcost} can be considered as an approximation to the asymptotic stabilization problem were $T=\infty$. This scenario will be the focus of our numerical tests. This also motivates that in the rest
of the paper we will restrict our presentation to the approximation of $V(0,\bx)=V(\bx)$ and $\nabla V(0,\bx)=\nabla V(\bx)$, and to the associated static feedback law $\bu^*(0,\bx)=\bu^*(\bx)$.  This approximation also relates to the one  done in nonlinear model predictive control, where after an open-loop solve is computed, the initial optimal control $u^*(0)$ is used to evolve the state equation for a short period , after which the open-loop optimization is re-computed with an updated initial state.  It can be shown that as the prediction horizon $T$ increases, the optimal control approaches the stationary feedback laws, see for instance    ~\citep{AK2016,GR2008,AR12,WRGA14}. For the reader interested in obtaining the complete time-dependent optimal feedback law, we discuss at the end of Section \ref{sec:approx} how to extend the proposed methodology.

The rest of the paper is structured as follows. In Section 2, we describe our numerical methodology, including the numerical generation of the dataset, the polynomial ansatz for the value function and the model fit through LASSO regression. Then, in Section 3 we present an exhaustive numerical assessment of the proposed methodology, including the synthesis of high-dimensional optimal feedback laws for nonlinear PDEs and multiagent systems.

\section{Data-driven Recovery of Feedback Laws}\label{sec:approx} 
In this section, we develop the different building blocks of the proposed approach. We first discuss how to generate the training dataset by solving a set of open-loop optimal control problems with a reduced gradient approach. Next, we build a polynomial model for the value function based on a hyperbolic cross approximation. Having specified both the data and the model, we fit our model with a LASSO regression. At the end of the section, we explain how to modify the proposed framework to recover time-dependent feedback laws.

\subsection{Generating a dataset with a reduced gradient approach} \label{sec:rdg} 
We begin by generating a dataset $\{\bx^j,V(0,\bx^j),\nabla V(0,\bx^j)\}_{j=1}^{N}$ which is obtained from solving open-loop optimal control problems of the form \eqref{eq:optcost} (with $t_0=0$ and $T$ sufficiently large) through the use of first-order optimality conditions \eqref{eq:state}-\eqref{eq:optcon}. For this purpose we follow a reduced gradient approach with a Barzilai-Borwein update ~\citep{BB,BB2}, which is summarized as follows.

Assuming that the solution operator  $ \vecy(\vu)=\vecy( \vu,\bx )$ corresponding to the state equation \eqref{eq:state}  is well-defined and continuously differentiable, we can rewrite \eqref{eq:optcost}-\eqref{eq:state} as the following unconstrained dynamic optimization problem depending solely on the control variable $\vu$ 
\begin{equation}
	\underset{\vu(\cdot)}{\min}\;\mathcal{J}(\vu) = \underset{\vu(\cdot)}{\min}\;J(\vecy(\vu),\vu) =\underset{(\vecy(\cdot),
		\vu(\cdot))}{\min} \{ J(\vecy,\vu) : \text{ subject  to } e(\vecy, \vu)=0 \},
\end{equation}
where
\begin{equation}
	e(\vecy,\vu):=
	\begin{pmatrix}
		\frac{d}{dt}\vecy(t) -\left(\bff(\vecy(t))+\bg(\by(t))\vu(t)\right)\\\
		\vecy (0)-\bx
	\end{pmatrix}.
\end{equation}
Formally, we obtain the directional derivative of   $\mathcal{J}$ at $\bar{\vu} \in L^2(0,T;\mathbb{R}^m)$  in a direction  $\delta \vu \in L^2(0,T;\mathbb{R}^m)$ by computing
\begin{equation}
	\label{e19}
	\mathcal{J}'(\bar{\vu})\delta\vu  = (\mathcal{G}(\bar{\vu}),\delta \vu) = ((\vecy'(\bar{\vu}))^* \partial_{\vecy} J(\bar{\vv})+\partial_{\vu}J(\bar{\vv}),\delta \vu),
\end{equation}
where  $\bar{\vv}:=(\bar{\vecy} ,\bar{\vu})$ with $\bar{\vecy} := \vecy(\bar{\vu})$,  $\mathcal{G}$ denotes the gradient of  $\mathcal{J}$,  $(\cdot,\cdot)$ stands for the scalar product in the space of controls $L^2(0,T;\mathbb{R}^m)$, and the superscript $*$ corresponds to the adjoint operator. Moreover,  the term $\vecy'(\bar{\vu})$ is given  by 
\begin{equation}
	\label{e20}
	\vecy'(\bar{\vu})\delta \vu = -(\partial_{\vecy}e( \bar{\vv}))^{-1}\partial_{\vu}e(\bar{\vv})\delta \vu.
\end{equation}
It can be shown that   $(\partial_{\vecy}e(\vecy(\vu),\vu))^{-1}$, defined  by $ (\phi,\vq_0)  \mapsto \vq $,   is the solution operator of the following linearised equation 
\begin{equation}
	\begin{cases}
		\frac{d}{dt}\vq(t)-\partial_{\vecy}\left(\bff(\vecy(t))+\bg(\by(t))\vu(t)\right)\vq(t) =\phi ,\\
		\vq (0)=\vq_0,&
	\end{cases}
\end{equation}
and that the adjoint operator  $(\partial_{\vecy}e(\vecy(\vu),\vu))^{-*}$ defined by $ (\psi,\vp_T)  \mapsto \vp $, is the solution operator to the following backward-in-time equation
\begin{equation}\label{eq:adjoint}
	\begin{cases}
		-\frac{d}{dt}\vp(t)-(\partial_{\vecy}\left(\bff(\vecy(t))+\bg(\by(t))\vu(t)\right)^{\top}\vp(t) =\psi ,\\
		\vp (T)=\vp_T\,.&
	\end{cases}
\end{equation}
Putting these elements together, we are now in a position to compute the gradient  $\mathcal{G}$  of  $\mathcal{J}$ at $\bar{\vu}$. Using \eqref{e19} and \eqref{e20},  we obtain   
\begin{equation}
	\mathcal{G}(\bar{\vu}) = \partial_{\vu}J(\bar{\vv})-(\partial_{\vu}e(\bar{\vv}))^{*}(\partial_{\vecy}e(\bar{\vv}))^{-*}\partial_{\vecy} J(\bar{\vv}) \,, 
\end{equation}
and therefore, for almost every $t \in (0,T)$, we have
\begin{equation}
	\label{e21}
	\mathcal{G}(\bar{\vu})(t) =  \bg(\bar{\by}(t))^{\top} \bar{\vp}(t)+2\beta \bar{\vu}(t),
\end{equation}
where $ \bar{\vp}$ is the solution to 
\begin{equation}
	\label{e22}
	\begin{cases}
		-\frac{d}{dt}\vp(t)-(\partial_{\vecy}\left(\bff(\bar{\vecy}(t))+\bg(\bar{\vecy}(t))\bar{\vu}(t)\right)^{\top}\vp(t)  =\partial_{\vecy}\ell(\bar{\vecy}(t)),\\
		\vp (T)=0.&
	\end{cases}
\end{equation}
Having a realization of the reduced gradient, we follow the Barzilai-Borwein gradient method 
for finding the stationary point ${\vu}^*$ of  $\mathcal{J}$  ( i.e. $\mathcal{G}({\vu}^*) = 0 $).  In this method, the stepsizes are chosen according to be either
\begin{equation}
	\label{e1a}
	\alpha^{BB1}_k := \frac{(\mathcal{S}_{k-1},\mathcal{Y}_{k-1})}{(\mathcal{S}_{k-1},\mathcal{S}_{k-1})}, \quad \text{ or } \quad
	\alpha^{BB2}_k := \frac{ (\mathcal{Y}_{k-1},\mathcal{Y}_{k-1})}{(\mathcal{S}_{k-1},\mathcal{Y}_{k-1})},
\end{equation}
where  $\mathcal{G}_k:=\mathcal{G}({\vu}_k)$,  $\mathcal{S}_{k-1}:= {\vu}_k-{\vu}_{k-1}$ and $\mathcal{Y}_{k-1}:=\mathcal{G}_k-\mathcal{G}_{k-1}$. With these specifications, we introduce Algorithm \ref{BBa}, which is used for solving the open-loop problems. 
\begin{algorithm}[htbp]
	\caption{Barzilai-Borwein two-point step-size gradient method}\label{BBa}
	\begin{algorithmic}[1]
		\INPUT   Choose ${\vu}_{-1}: = 0$  and  ${\vu}_0 := -\mathcal{G}(0)$, tolerance $tol>0$.
		\State Set $k=0$.
		\While{$\|\mathcal{G}_k\|\geq tol\,,$} 
		\State Compute $\by_k(\vu_k)$ via \eqref{eq:state}.
		\State Compute $\vp_k(\by_k,\vu_k)$ via \eqref{e22}.
		\State Compute $\mathcal{G}_k =  \mathcal{G}({\vu}_k)$ using \eqref{e21} with $(\by_k,\vp_k,\vu_k)$.
		\State  Choose 
		\begin{equation*}
			\alpha_k = 
			\begin{cases}
				\alpha^{BB1}_k  &   \text{ for  odd } k, \\
				\alpha^{BB2}_k  &   \text{ for  even } k. 
			\end{cases}
		\end{equation*}
		\State Set $d_k =  \frac{1}{\alpha_k}\mathcal{G}_k$.
		
		\State Compute the step-size $\eta_k>0$ based on the non-monotone linesearch given in ~\citep{DH2001}. 
		\State  Set ${\vu}_{k+1} = {\vu}_k- \eta_k d_k$, $k=k+1$, and go to Step 2.
		\EndWhile
	\end{algorithmic}
\end{algorithm}
Note that the formulas above are written for continuous-time dynamical systems. In practice, Algorithm \ref{BBa} is expected to be used in conjunction with a suitable numerical integrator for an accurate approximation of both the state and its adjoint. Finally, fixing an initial condition $\bx^j$, after Algorithm \ref{BBa} has converged to $\bu^*$ the dataset is completed with $V(\bx^j)=\mathcal{J}(\bu^*)$ and $\nabla V(\bx^j)=\vp^*(0)$.
\subsection{Building a polynomial model for the value function}
Having generated a dataset for recovering the value function associated to the optimal control problem, we now turn our attention to deriving a suitable model for regression. Our approximation of the static value function $V(\bx):\mathbb{R}^n \to \mathbb{R}$ follows the ideas presented in ~\citep{MR3975885,MR3962896}. 

\begin{figure}[t!]
	\centering
	\includegraphics[width=0.33\textwidth]{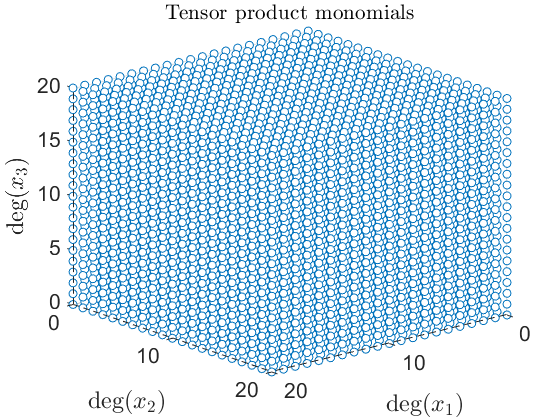}\hfill
	\includegraphics[width=0.33\textwidth]{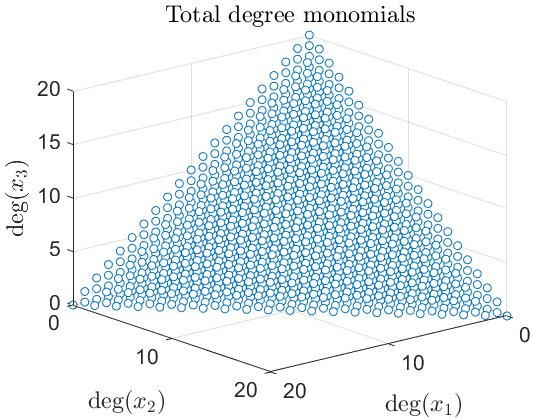}\hfill
	\includegraphics[width=0.33\textwidth]{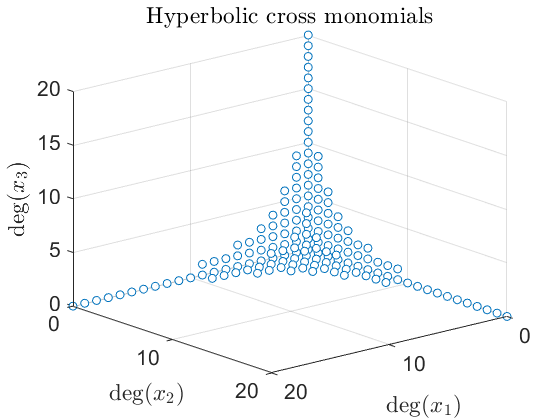}
	\caption{Alternatives for generating a high-dimensional polynomial basis. From left to right: direct tensorization of a 1-dimensional basis, truncation by total degree, hyperbolic cross approximation.}\label{hc}
\end{figure}

Let $\mathcal{D} \subset \mathbb{R}^n$ be a bounded domain and $\{\Phi_\mathbf{i}\}_{\mathbf{i} \in \mathbb{N}^n_0}$  be a tensor-product orthonormal  basis of  $L^2(\mathcal{D})$. We consider bases which are polynomial, using for instance Legendre or Chebyshev polynomials.  Concretely, assume that $\mathcal{D} : = (-1,1)^n$ and that $\{\phi_i\}^{\infty}_{i = 0}$ is one-dimensional orthonormal basis of $L^2(-1,1)$. Then, the corresponding tensor-product basis of $L^2(\mathcal{D})$ is defined by 
\begin{equation}
	\Phi_{\mathbf{i}}(\mathbf{x}) := \prod^n_{j =1} \phi_{i_j}(x_i)\,,\quad  \text{ with } \quad  \mathbf{i} =(i_1,i_2,\dots,i_n) \in \mathbb{N}^n,  \quad \mathbf{x}=(x_1,x_2,\dots,x_n),  
\end{equation} 
where $\mathbb{N}_0: = \mathbb{N}\cup \{ 0\}$.  Assuming that  $V(\bx)\in L^2(\mathcal{D}) \cap L^{\infty}(\mathcal{D})$, we can write 
\begin{equation}
	V(\bx) = \sum_{\mathbf{i} \in \mathbb{N}^n_0} \theta_{\mathbf{i}} \Phi_{\mathbf{i}}\,,
\end{equation}
with $\theta_{\mathbf{i}} = (V(\bx),\Phi_{\mathbf{i}} )_{L^2(\mathcal{D})}$ for every $\mathbf{i} \in \mathbb{N}^n_0$. We approximate $V(\bx)$ by considering a truncated basis $\{\Phi_\mathbf{i}\}_{\mathbf{i} \in \mathfrak{I}}$ with a finite  multi-index set  $\mathfrak{I} \subset \mathbb{N}^n_0$  with cardinality $| \mathfrak{I}| = q < \infty$. Hence, we write
\begin{equation}
	V(\bx) = V^{\mathfrak{I}}_{\theta}+e_{^{\mathfrak{I}}} = \sum_{\mathbf{i} \in \mathfrak{I}} \theta_{\mathbf{i}} \Phi_{\mathbf{i}} +\sum_{\mathbf{i} \notin  \mathfrak{I}} \theta_{\mathbf{i}} \Phi_{\mathbf{i}} 
\end{equation}
with $\{ \theta_{\mathbf{i}}\}_{\mathbf{i} \in \mathbb{N}_0^n} \in \ell^	2( \mathbb{N}_0^n)$. 
In this work we are particularly interested in the case where $\mathcal{D}$ is a high-dimensional space and computing $V(\bx)$ by solving a HJB PDE is not a feasible alternative. Therefore, the selection of a basis $\basis$ whose cardinality scales reasonably well in high-dimensions, while maintaining an acceptable level of accuracy, is a fundamental criterion in our model selection. Figure \ref{hc} illustrates some of the typical options for generating a multi-dimensional polynomial basis. Generating a basis by directly taking the tensor product of polynomials up to a certain degree $s$ leads to
\begin{equation}
	\label{tp}
	\mathfrak{I}_{TP}(s) = \left\lbrace \mathbf{i}=(i_1,i_2,\dots, i_n) \in \mathbb{N}_0^n : \|\veci\|_{\infty}\leq s  \right\rbrace\,,
\end{equation}
with $|\mathfrak{I}_{TP}(s)|=(s+1)^n$, scaling exponentially in the dimension, limiting its applicability to $n\leq 5$ unless additional low-rank structures are assumed ~\citep{tensor4}. This exponential increase in the dimension can be mitigated by considering a total degree truncation 
\begin{equation}
	\label{td}
	\mathfrak{I}_{TD}(s) = \left\lbrace \mathbf{i}=(i_1,i_2,\dots, i_n) \in \mathbb{N}_0^n : \|\veci\|_1\leq s  \right\rbrace\,,
\end{equation}
with cardinality
\begin{equation}
	|\mathfrak{I}_{TD}(s)| = \sum\limits_{j=1}^s {n+j-1\choose j} \,.
\end{equation}
This combinatorial dependence on the dimension allows to solve moderately high-dimensional problems, in our experience for $n\leq 15$ ~\citep{poli1,poli2}. In this work, we opt for a basis constructed with the \textsl{hyperbolic cross} index set, defined as
\begin{equation}
	\label{hcset}
	\mathfrak{I} = \mathfrak{I}(s) = \left\lbrace \mathbf{i}=(i_1,i_2,\dots, i_n) \in \mathbb{N}_0^n : \prod^n_{j = 1}(i_j+1) \leq s+1  \right\rbrace.
\end{equation}
While there are no explicit formulas for the cardinality of $\basis$, different upper bounds ~\citep{aw17} such as
\begin{equation}
	|\mathfrak{I}(s)|\leq\min\left\{2s^34^n,e^2s^{2+\log_2(n)}\right\} \,,
\end{equation}
indicate that it scales reasonably well for high-dimensional problems. For reference,  in this paper we report results up to $n=80$ at moderate computational cost. With 80 dimensions and degree 4, the tensor product basis would contain, $8.27\times 10^{55}$ elements, the total degree basis $1.93\times 10^6$ elements, while the upper bound  above for the hyperbolic cross above is $7.56\times 10^{5}$. Besides the dimensionality argument, the hyperbolic cross set is also an adequate basis regarding best approximation properties ~\citep{MR3975885} in conjunction with the $\ell_1$ regression framework we will discuss in the following section. For the rest of the paper, we will adopt the notation $\mathbf{i}_1 , \mathbf{i}_2, \dots, \mathbf{i}_q(s)$ for an order of multi-indices in  $\mathfrak{I}(s)$, and we write $\theta_{\mathfrak{I}} = \{ \theta_{\mathbf{i}}\}_{\mathbf{i} \in \mathfrak{I}} =\{ \theta_{\mathbf{i}_k}\}^q_{k = 1} $.

\subsection{Gradient-augmented regression}
The last building block of our approach consists in fitting the polynomial expansion presented above with the dataset containing information of both the value function and its gradient. We first present an un-augmented linear least squares approach which will be used for numerical comparison. Based on the data generation procedure presented in section \ref{sec:rdg}, we assume the existence of a dataset consisting of $N$ samples
\begin{equation}
	D =  \left\lbrace \mathbf{x}^j ,  V^j  \right\rbrace^{N}_{j = 1},
\end{equation} 
where  $V^j: = V(\mathbf{x}^j)$.  By defining  $\mathbf{V} \in \mathbb{R}^{N}$, $\mathbf{A} \in \mathbb{R}^{N\times q}$, and  $\mathbf{e}\in \mathbb{R}^{N}$   by
\begin{equation}
	\mathbf{V} := \frac{1}{\sqrt{N}}\left( V(\mathbf{x}^j)\right)^{N}_{j = 1}, \quad  \mathbf{A} := \frac{1}{\sqrt{N}}\left( \Phi_{\mathbf{i}_k}(\mathbf{x}^j)\right)^{N,q}_{ j,k = 1}, \quad \text{ and } \mathbf{e} := \frac{1}{\sqrt{N}}\left( e_{\mathfrak{I}}(\mathbf{x}^j)\right)^{N}_{j = 1},
\end{equation}
we write the following linear system to be satisfied by our model
\begin{align}
	\label{e3}
	\mathbf{V}  = \left[\frac{1}{\sqrt{N}} V_{\theta}(\bx^i)\right]_{i=1}^{N}=  \mathbf{A}\theta_{\mathfrak{I}} + \mathbf{e}\,,   
\end{align}
from which $V(\bx)$ can be approximated on the subspace $\{\spann(\{  \Phi_{\mathbf{i}} \}_{\mathbf{i} \in \mathfrak{I}})\}$,  by solving the linear least squares problem
\begin{equation}
	\label{e4}
	\tag{$P_{\ell_2}$}
	\min_{\theta \in \mathbb{R}^q }\|  \mathbf{A}\theta- \mathbf{V} \|^2_2.
\end{equation}
In order to  enhance sparsity in the vector of coefficients, we consider  the \textit{weighted LASSO} regression 
\begin{equation}
	\label{e5}
	\tag{$P_{\ell_1}$}
	\min_{\theta \in \mathbb{R}^q }\|  \mathbf{A}\theta- \mathbf{V} \|^2_2 + \lambda \| \theta\|_{1, \mathbf{w}},
\end{equation}
where $\lambda>0\,,$ and the weights $\mathbf{w}:=\{ w_{\mathbf{i}}\}_{\mathbf{i} \in \mathbb{N}^n_0}$ with $w_{\mathbf{i}}\geq 1$ define  the $\| \cdot\|_{1, \mathbf{w}}$ norm as 
\begin{equation}
	\|\theta\|_{1, \mathbf{w}}  = \sum_{ \mathbf{i} \in \mathbb{N}^n_0}  w_{\mathbf{i}}  |\theta_{\mathbf{i}}|.
\end{equation}
Denoting by  $\theta_{\ell_2} \in \mathbb{R}^q$ and  $\theta_{\ell_1} \in \mathbb{R}^q$ the solutions to problems \eqref{e4} and \eqref{e5}, respectively, we recover the following approximations of  $V(\bx)$
\begin{equation}
	\label{e9}
	V_{\ell_2}(\bx) = \sum^q_{ k=1} (\theta_{\ell_2})_{ \mathbf{i}_k}\Phi_{\mathbf{i}}(\bx) =  \sum_{ \mathbf{i} \in  \mathfrak{I}} (\theta_{\ell_2})_{ \mathbf{i}}\Phi_{\mathbf{i}}(\bx) \quad  \text{ and } \quad  V_{\ell_1} = \sum_{ \mathbf{i} \in  \mathfrak{I}} (\theta_{\ell_1})_{ \mathbf{i}}\Phi_{\mathbf{i}}(\bx). 
\end{equation}
Building on the fact that our data generation procedure retrieves the value function and its gradient, we recast the regression problems by using the augmented data ~\citep{MR3975885}. We assume $V(\bx) \in H^1(\mathcal{D})$.
For the augmented dataset
\begin{equation}
	D_{aug} =  \left\lbrace \mathbf{x}^j , V^j , V_x^j \right\rbrace^{N}_{j = 1},    \quad  \text{ with }   V^j_x = \left(  \frac{\partial V_T}{\partial x_1}(\mathbf{x}^j),  \frac{\partial V_T}{\partial x_1}(\mathbf{x}^j), \dots, \frac{\partial V_T}{\partial x_n}(\mathbf{x}^j) \right)^{\top} \text{ for } j = 1,\dots, N,
\end{equation} 
and, for every $m = 0,\dots, n$,  we define    
\begin{equation}
	\mathbf{A_m} := \frac{1}{\sqrt{N}}\left(  \frac{\partial \Phi_{\mathbf{i}_k} }{\partial x_m}(\mathbf{x}^j) \right)^{N,q}_{ j,k = 1}, \quad    \mathbf{V}_m := \frac{1}{\sqrt{N}}\left(\frac{\partial  V_T}{\partial x_m}(\mathbf{x}^j) \right)^{N}_{j = 1},         \text{ and }    \mathbf{e}_m :=\frac{1}{\sqrt{N}} \left(\frac{\partial e_{\mathfrak{I}}}{\partial x_m}(\mathbf{x}^j) \right)^{N}_{j = 1} ,
\end{equation}
where 
\begin{equation}
	\mathbf{A_0} := \frac{1}{\sqrt{N}}\left(  \Phi_{\mathbf{i}_k}(\mathbf{x}^j) \right)^{N,q}_{ j,k = 1},  \quad  \mathbf{V}_0 := \frac{1}{\sqrt{N}}\left( V_T(\mathbf{x}^j) \right)^{N}_{j = 1},   \text{ and  }  \mathbf{e}_0 := \frac{1}{\sqrt{N}}\left( e_{\mathfrak{I}}(\mathbf{x}^j)\right)^{N}_{j = 1}.
\end{equation}
Then by assembling
\begin{equation}
	\bar{\mathbf{A}} := \begin{pmatrix}
		\mathbf{A_0} \\
		\mathbf{A_1}  \\
		\vdots\\
		\mathbf{A_n} 
	\end{pmatrix}, 
	\quad 
	\bar{\mathbf{V}} := \begin{pmatrix}
		\mathbf{V_0} \\
		\mathbf{V_1}  \\
		\vdots\\
		\mathbf{V_n} 
	\end{pmatrix}, 
	\text{ and  }
	\bar{\mathbf{e}} := \begin{pmatrix}
		\mathbf{e_0} \\
		\mathbf{e_1}  \\
		\vdots\\
		\mathbf{e_n} 
	\end{pmatrix},
\end{equation}
we obtain the following system of  linear equations
\begin{equation}
	\bar{\mathbf{V}} =  \bar{\mathbf{A}}\theta_{\mathfrak{I}} + \bar{\mathbf{e}}\,,   
\end{equation}
for which we formulate  the augmented-gradient optimization problems  
\begin{equation}
	\label{e6}
	\tag{$AP_{\ell_2}$}
	\min_{\theta \in \mathbb{R}^q }\|  \bar{\mathbf{A}}\theta- \bar{\mathbf{V}} \|^2_2,
\end{equation}
and
\begin{equation}
	\label{e7}
	\tag{$AP_{\ell_1}$}
	\min_{\theta \in \mathbb{R}^q }\|  \bar{\mathbf{A}}\theta-\bar{\mathbf{V}} \|^2_2 + \lambda \| \theta\|_{1, \mathbf{w}}.
\end{equation}
Denoting by  $\bar{\theta}_{\ell_2} \in \mathbb{R}^q$ and  $\bar{\theta}_{\ell_1} \in \mathbb{R}^q$  the solutions to \eqref{e6} and \eqref{e6}, respectively,  we recover the following gradient-augmented approximations of $V(\bx)$
\begin{equation}
	\label{e10}
	\bar{V}_{\ell_2}(\bx) =  \sum_{ \mathbf{i} \in  \mathfrak{I}} (\bar{\theta}_{\ell_2})_{ \mathbf{i}}\Phi_{\mathbf{i}}(\bx) \quad  \text{ and } \quad  \bar{V}_{\ell_1}(\bx) = \sum_{ 
		\mathbf{i} \in  \mathfrak{I}} (\bar{\theta}_{\ell_1})_{ \mathbf{i}}\Phi_{\mathbf{i}}(\bx).  
\end{equation}
Similarly, we obtain the following representations for  $\nabla_x V(\bx) $,  where   $\nabla_x = ( \frac{\partial}{\partial x_1},  \dots, \frac{\partial}{\partial x_n})^{\top} $
\begin{equation}  
	\label{e11}
	\nabla_x \bar{V}_{\ell_2} (\bx)=  \sum_{ \mathbf{i} \in  \mathfrak{I}} (\bar{\theta}_{\ell_2})_{ \mathbf{i}}\nabla_x \Phi_{\mathbf{i}}(\bx) \quad  \text{ and } \quad \nabla_x \bar{V}_{\ell_1} (\bx)= \sum_{ \mathbf{i} \in  \mathfrak{I}} (\bar{\theta}_{\ell_1})_{ \mathbf{i}}\nabla_x\Phi_{\mathbf{i}}(\bx)\,\,
\end{equation}
recovering the optimal feedback laws
\begin{equation} 
	\label{eq:feedback}
	\bu_{\star}^*(\bx)= - \frac{1}{2\beta}\bg^{\top}(\bx)\nabla_x \bar V_{\star}(\bx)\,,\qquad \star\in\{\ell_1,\ell_2\}\,.
\end{equation}

\paragraph{On the numerical realization of the weighted LASSO regression.} The linear least squares problems \eqref{e4} and \eqref{e6} can be efficiently solved by using a preconditioned conjugate gradient method. The formulations \eqref{e5} and \eqref{e7} are instead convex, nonsmooth optimization problem which require a more elaborate treatment. In this work, we compute the solution of \eqref{e5} and \eqref{e7} by means of the  Alternating Direction Method of Multipliers (ADMM) ~\citep{boyd2011distributed}. To make matters precise, Algorithm \ref{ADMM} presents its implementation for problem \eqref{e5}.
\begin{algorithm}[htbp]
	\caption{ADMM for solving weighted LASSO}\label{ADMM}
	\begin{algorithmic}[1]
		\INPUT   Choose $\theta^0, z^0, h^0 \in \mathbb{R}^q$,  $\rho>0$, and tolerance $tol>0$.
		\State Set  $k =0$.
		\While {$\|\theta^k-z^k\| \geq tol$ and $\|\rho(h^k-h^{k-1})\| \geq tol$ }
		\State $\theta^{k+1} = \left(2\mathbf{A}\mathbf{A}^{\top}+\rho \mathbf{I} \right)^{-1}\left(2\mathbf{A}^{\top} \mathbf{V}+ \rho(z^k-h^k)\right)$   
		\State $z^{k+1} = \prox_{\frac{\lambda}{\rho}\| \cdot\|_{1,\mathbf{w}}}(\theta^{k+1}+h^k)$.
		\State $h^{k+1} = h^k+\theta^{k+1}-z^{k+1}$.
		\State Set  $k = k+1$ and go to Step 2.
		\EndWhile
	\end{algorithmic}
\end{algorithm}
In this algorithm, $\mathbf{I}$ stands for the  identity matrix and the proximal operator $\prox_{\frac{\lambda}{\rho}\| \cdot\|_{1,\mathbf{w}}}$
is explicitly given by a soft-thresholding type operator ~\citep[Chapter 6]{F-OMO}
\begin{equation}
	\prox_{\frac{\lambda}{\rho}\| \cdot\|_{1,\mathbf{w}}}(\mathbf{x}) = \left( [ |x_i| - \frac{\lambda w_i }{\rho} ] _+ \sgn(x_i)     \right)^q_{i = 1}    \text{ for  }   \mathbf{x} = (x_1, \dots, x_q)\,, 
\end{equation} 
where $[\cdot]_+$ denotes the positive part. The application of Algorithm \ref{ADMM} for problem \eqref{e7} is directly done by  replacing $\mathbf{A}$ and $\mathbf{V}$ by  $\bar{\mathbf{A}}$ and $\bar{\mathbf{V}}$, respectively.

\subsection{Recovering time-dependent feedback laws.} The present computational framework can be extended to recover time-dependent value functions and feedback laws. While we argue that the primary object of study in deterministic optimal control of physical systems is the synthesis of static feedback laws, there exist applications such as operations research and {\color{black}the stabilization to non-stationary trajectories ~\citep{sergio}} where the computation of time-dependent feedback controls is of great interest. As discussed at the end of Section \ref{intro}, the solution of the optimal control problem for a given initial condition $(t^j,\bx^j)$ generates data for $V(t,\by^*(t))$ and $\nabla V(t,\by^*(t))$ with $t\in(t^j,T)$, along the optimal trajectory departing from $\by^*(t^j)=\bx^j$. From an approximation viewpoint, this space-time data can be use to fit a model for $V(t,\bx)$. The simplest option is to approximate $V(t,\bx)$ along the space-time cylinder treating time in the same way as $\bx$, that is
\begin{align}
	V(t,\bx)\approx V_{\theta}(\tilde\bx)=\sum_{\veci\in\basis}\theta_\veci\Phi_\veci(\tilde\bx)\,,\qquad \tilde\bx=(t,\bx)\in\R^{n+1}\,.
\end{align}
The computational cost of this augmented representation is related to the new polynomial basis and the increase of $|\basis(s)|$ in $\R^n$ to $\R^{n+1}$. An alternative to this treatment is to establish a time-marching structure for $t$, and embed this time dependence in $\theta$, similar to the method of lines for parabolic PDEs ~\citep{MOL}. Formally, we write
\begin{align}
	V(t,\bx)\approx V_{\theta}(t,\bx)=\sum_{\veci\in\basis}\theta_\veci(t)\Phi_\veci(\bx)\,,
\end{align}
where an additional approximation for $\theta$ is required. It is reasonable to assume that the artificial dataset from the numerical optimal control solutions will be provided as a time series with a uniform time discretization parameter $\tau$, so we use a piecewise constant approximation for $\theta_{\veci}(t)$
\begin{align}
	\theta_{\veci}(t)\approx \theta_{k,\veci}\,,\qquad \text{ for }\, t\in [(k-1)\tau,k\tau)\,,\qquad k=1,\ldots,N_T
\end{align}
where $N_T=T/\tau$, and the space-time approximation of $V(t,\bx)$ becomes
\begin{align}
	V_{\theta}(t,\bx)=\sum\limits_{k=1}^{N_T}\sum_{\veci\in\basis}\theta_{k,\veci}\Phi_\veci(\bx)\,,
\end{align}
which can be further simplified in the absence of terminal penalties in the cost functional since in this case $V(T,\bx)=0$. Thus, the computational increase is linear with respect to the cost associated to the static feedback law. A high-order discretization in time can be used to reduce the number of time nodes, however, it is necessary to always maintain a linear structure in $\theta$.

\section{Numerical  Tests}
In this section we assess the proposed methodology for recovering optimal feedback laws in three different tests. After presenting the practical aspects of our numerical implementation, we study the control of a nonlinear, low-dimensional oscillator. Then, we study large-dimensional dynamics arising in optimal control of nonlinear parabolic PDEs and non-local agent-based dynamics. In these tests, we focus on studying the effects of the sparse regression and the selection of weights in the $\ell_1$ penalty, the gradient-augmented recovery, the selection of a suitable polynomial basis, and the effectiveness of the recovered control law. Both sampling and regression algorithms were implemented in MATLAB R2014b, and the numerical tests were run in a  MacBook Pro with $2.9$ GHz Dual-Core Intel Core i5 and memory  $16$ GB $1867$ MHz DDR$3$.
\subsection{Practical aspects}\label{sec:pract} 
\paragraph{Generating the samples.} For each test we fixed an $n-dimensional$ hyperrectangle as the domain for sampling initial condition vectors $\{ \mathbf{x}^j\}^{N}_{j =1} \in \mathbb{R}^n$.  These initial vectors were generated using Halton quasi-random sequences\footnote{\url{https://www.mathworks.com/help/stats/generating-quasi-random-numbers.html}} in dimension $n$. Then for every $i\in \{1,\dots,N\}$, we compute the value function $V^j = V(0,\mathbf{x}^j)$ by solving the open-loop optimal control problem \eqref{eq:optcost}-\eqref{eq:state}. Every optimal control problem was  solved in the reduced form by  using Algorithm \ref{BBa} with $tol =10^{-5}$ as discussed in Section \ref{sec:rdg}. Note that the computational burden associated to solving an optimal control problem for each initial condition of the ensemble can be alleviated by directly parallelising this task.  Further, for each control problem, the gradient of the value function was  obtained by evaluating  the solution $\vp^*$ of the adjoint equation \eqref{e22} at initial time $t_0 = 0$, so that  $\nabla V_x^j(\bx^j)=\nabla V(0,\bx^j)=\vp^*(0)$. This quantity is obtained as a by-product of solving the optimal control problem at no additional cost.

\paragraph{Training and validation.} We split the  sampling dataset $\{ \mathbf{x}^j, V^j,  V_x^j \}^{N}_{j = 1}$ into two sets:  a set of training indices $\mathcal{I}_{tr}$ which is used for regression, and a set of validation indices  $\mathcal{I}_{val}$, with  $\mathcal{I}_{val} \cup \mathcal{I}_{tr}  = \{ 1,\dots,N\}$. Without loss of generality, we assume that $\mathcal{I}_{tr}  = \{ 1,\dots,N_d\}$ and   $\mathcal{I}_{val}  = \{ N_d+1,\dots,N\}$ for $N \in \mathbb{N}$ with  $ N_d < N$. 

The linear least square problems \eqref{e4} and \eqref{e6} were solved using a preconditioned
conjugate gradient method, and the algorithm was terminated when the norm of residual was less than $10^{-8}$. For the LASSO regressions \eqref{e5} and \eqref{e7}, we employed Algorithm \ref{ADMM} with $tol = 10^{-5}$. 

To analyse the generalization error of the approximated value function $\bar{V}(\bx)$ with respect to the exact value $V(0,\bx)$, we use the following relative errors:
\begin{equation}
	\label{eq:kk1}
	\begin{split}
		Err_{L^2}(\bar{V})  &=\left(\frac{\sum\limits_{j \in \mathcal{I}_{val}}| \bar{V}(\bx^j)-V(0,\bx^j)|^2}{ \sum\limits_{j \in \mathcal{I}_{val}}| V(0,\bx^j)|^2}\right)^{\frac12}\,, \\
		Err_{H^1}(\bar{V})  &=\left(\frac{\sum\limits_{j \in \mathcal{I}_{val}}\left( | \bar{V}(\mathbf{x}^j)-V(0,\mathbf{x}^j)|^2 +\sum^n_{i=1}\left| \frac{\partial \bar{V}(\mathbf{x}^j)}{\partial x_i}  -\frac{\partial V(0,\mathbf{x}^j)}{\partial x_i}\right|^2    \right)}{\sum\limits_{j \in \mathcal{I}_{val}}\left( |V(0,\mathbf{x}^j)|^2 +\sum^n_{i=1}\left| \frac{\partial V(0,\mathbf{x}^j)}{\partial x_i}\right|^2    \right)}\right)^{\frac12}\ .
	\end{split}
\end{equation}
\paragraph{Weights in the $\ell_1$ norm.} Regarding the selection of weights for the $\|\cdot\|_{1,\bw}$ norm, we consider expressions of the form
\begin{equation}\label{alpha}
	w_i =v^{\alpha}_i \;  \text{ for }   \alpha>0, 
\end{equation}
where the terms $v_{\mathbf{i}}$ depend on the polynomial basis chosen for regression. In the case of Legendre and Chebyshev polynomials, we proceed as in ~\citep{MR3440176}, which we summarize in the following. We consider tensorized Legendre polynomials on $\mathcal{D} =[-1,1]^n$ of the from
\begin{equation}
	\mathbf{L}_{\mathbf{i}}(\mathbf{x}) := \prod^n_{j =1} {L}_{i_j}(x_i)  \text{ with } \quad  \mathbf{i} =(i_1,i_2,\dots,i_n) \in \mathbb{N}^n,  \quad \mathbf{x}=(x_1,x_2,\dots,x_n),  
\end{equation}
with $L_k$ defined as the univariate orthonormal Legendre polynomials of degree $k$. In this case, the Legendre polynomials form a basis for the real algebraic polynomials on $\mathcal{D}$ and are orthogonal with respect to the tensorized uniform measure on $\mathcal{D}$. Moreover, due to the fact that $\|L_k\|_{L^{\infty}(-1,1)} \leq \sqrt{k}$, we can write that 
\begin{equation}\label{wleg}
	\|\mathbf{L}_{\mathbf{i}}\|_{L^{\infty}(\mathcal{D})} \leq \prod^n_{j =1}(1+ i_j)^{\frac{1}{2}},
\end{equation}    
from where we take the hyperbolic  cross  weights $v_{\mathbf{i}}= \prod^n_{j =1}(1+ i_j)^{\frac{1}{2}}$. 
\noindent For tensorized Chebyshev polynomials on $\mathcal{D}$
\begin{equation}
	\mathbf{C}_{\mathbf{i}}(\mathbf{x}) := \prod^n_{j =1} {C}_{i_j}(x_i)  \text{ with } \quad  \mathbf{i} =(i_1,i_2,\dots,i_n) \in \mathbb{N}^n,  \quad \mathbf{x}=(x_1,x_2,\dots,x_n),  
\end{equation}
with $C_k(x)= \sqrt{2}\cos((k-1)\arccos(x))$, the uniform bound $\|C_k\|_{L^{\infty}(-1,1)} \leq \sqrt{2}$ holds, leading to $\|\mathbf{C}_{\mathbf{i}}\|_{L^{\infty}(\mathcal{D})} \leq 2^{\frac{\|\mathbf{i}\|_0}{2}}$. The latter is a valid alternative for setting $v_{\mathbf{i}}$, however we note that the bound $\eqref{wleg}$ also holds in this case, so we choose $v_{\mathbf{i}}= \prod^n_{j =1}(1+ i_j)^{\frac{1}{2}}$. To fit with these settings, in our numerical experiments we rescale the sampling set of initial conditions $\mathcal{D}$ to the unit hypercube $[-1,1]^n$.

\subsection{Test 1: Van der Pol oscillator}\label{exp1}
We consider the optimal control of the Van der Pol oscillator expressed as
\begin{equation}
	\label{e26}
	\min_{\mathbf{u}\in L^2(0,T;\mathbb{R})} \int^{T}_{0} y^2_1(t)+y^2_2(t)+\beta \mathbf{u}^2(t)dt
\end{equation}
subject to
\begin{equation}
	\label{e27}
	\begin{cases}
		\partial_t y_1 =  y_2,\\
		\partial_t y_2 = -y_1+y_2(1-y^2_1)+u, \\
		(y_1(0),y_2(0))=(x_1,x_2),
	\end{cases}
\end{equation}
where we set  $\bx := (x_1,x_2)$,    $\beta =0.1$, and $T = 3$. A dataset  $\{\mathbf{x}^j,V^j, V^j_x \}^{N}_{j =1}$  with $N=2000$ is prepared by solving open-loop problems for different values of quasi-randomly chosen initial vectors from the domain $\mathcal{D}=[-3,3]^2$. The temporal discretization is done by the Crank-Nicolson time stepping method with step-size $\Delta t =10^{-4}$.
Here we set $s = 16$ in the Hyperbolic cross index set  $\mathfrak{I}(s)$  given in \eqref{hcset}. In this case  we have $q\equiv|\mathfrak{I}(s)| =52$. That is, we use $52$  polynomial  basis functions to approximate the value function $V(\bx): = V(0,\bx)$ corresponding to the optimal control problem \eqref{e26}-\eqref{e27}.

We computed the solutions $\theta_{\ell_2}$, $\bar{\theta}_{\ell_2}$, $\theta_{\ell_1}$ and $\bar{\theta}_{\ell_1}$ to the problems \eqref{e4} , \eqref{e6}, \eqref{e5}, and  \eqref{e7}, respectively, analysing different sizes for the training dataset $N_d$, the choice of  $\ell_1$ weights encoded through $\alpha$ in \eqref{alpha}, and polynomial bases. A compact summary of these results is given in Table
\ref{table1}.  The first column gives the errors of the value function and the nonzero components in its expansion without relying on gradient information and without sparsification. In the second column gradient information is added and the errors decrease for the same number of training samples ($N_d=40$). For the third column the sparsity enhancing functional is added and approximately the same errors are obtained with significantly fewer nonzero components in the expansion.

\begin{table}[htbp!]
	
	\begin{center}
		Training Errors\\\vskip 2mm
		\begin{tabular}{cccc}
			\hline\hline
			&  $Err_{L^2} $& $Err_{H^1} $&  Nonzero components  \\
			\cmidrule(lr){2-2}\cmidrule(lr){3-3}\cmidrule(lr){4-4} \morecmidrules\cmidrule(lr){2-2}\cmidrule(lr){3-3}\cmidrule(lr){4-4}
			$V_{\ell_2}$  for  $N_d =40$ & $9.93 \times 10^{-9}$ &$ 4.83  \times 10^{-1}
			$&$52/52$  \\
			$\bar{V}_{\ell_2}$ for $N_d =40$ &$1.38 \times 10^{-3}$ &$ 2.34 \times 10^{-3}$&$52/52$ \\
			$\bar{V}_{\ell_1}$ for $N_d =40$, $\lambda = 0.01$  &$ 1.05 \times 10^{-2}$ &$1.45\times10^{-2}$&$19/52$ \\
			\hline\hline
		\end{tabular}\\\vskip 5mm
		
		Validation Errors\\
		\vskip 2mm
		\begin{tabular}{cccc}			
			\hline\hline
			&  $Err_{L^2} $& $Err_{H^1} $&  Nonzero components  \\
			\cmidrule(lr){2-2}\cmidrule(lr){3-3}\cmidrule(lr){4-4} \morecmidrules\cmidrule(lr){2-2}\cmidrule(lr){3-3}\cmidrule(lr){4-4}
			$V_{\ell_2}$  for  $N_d =40$ & $1.46 \times 10^{-1}$ &$1.17$&$52/52$  \\
			$\bar{V}_{\ell_2}$ for $N_d =40$ &$ 9.38 \times 10^{-3}$ &$ 3.25\times10^{-2}$&$52/52$ \\
			$\bar{V}_{\ell_1}$ for $N_d =40$, $\lambda = 0.01$  &$ 1.20  \times 10^{-2}$  &$2.05\times10^{-2}$&$19/52$ \\
			\hline\hline
		\end{tabular}
		
	\end{center}
	\label{table1b}
	\caption{Test 1. Numerical Results for Legendre polynomial basis and $\alpha = 1$ for $\bar{V}_{\ell_1}$. Including gradient information and sparsification leads to less error and fewer components in the expansion with a reduced number of training samples.}
	\label{table1}
\end{table}

\begin{table}[htbp!]
	
\end{table}

\begin{figure}[htbp!]
	\centering
	\subfigure[Scatter plot of sampling data]
	{
		\label{Fig5}
		\includegraphics[height=3.5cm,width=3.5cm]{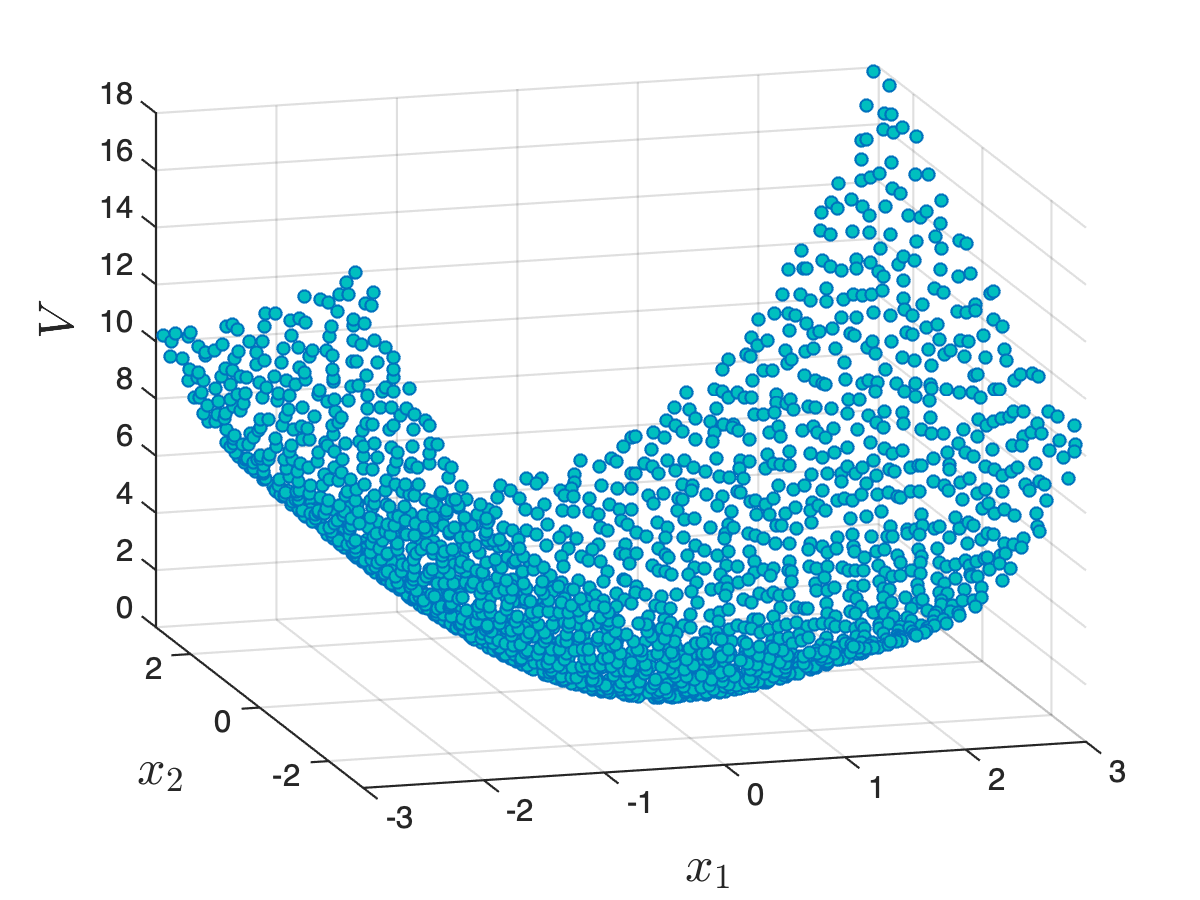}
	}
	\subfigure[$\bar{V}_{\ell_1} $ for $N_d =40$]
	{
		\label{Fig6}
		\includegraphics[height=3.5cm,width=3.5cm]{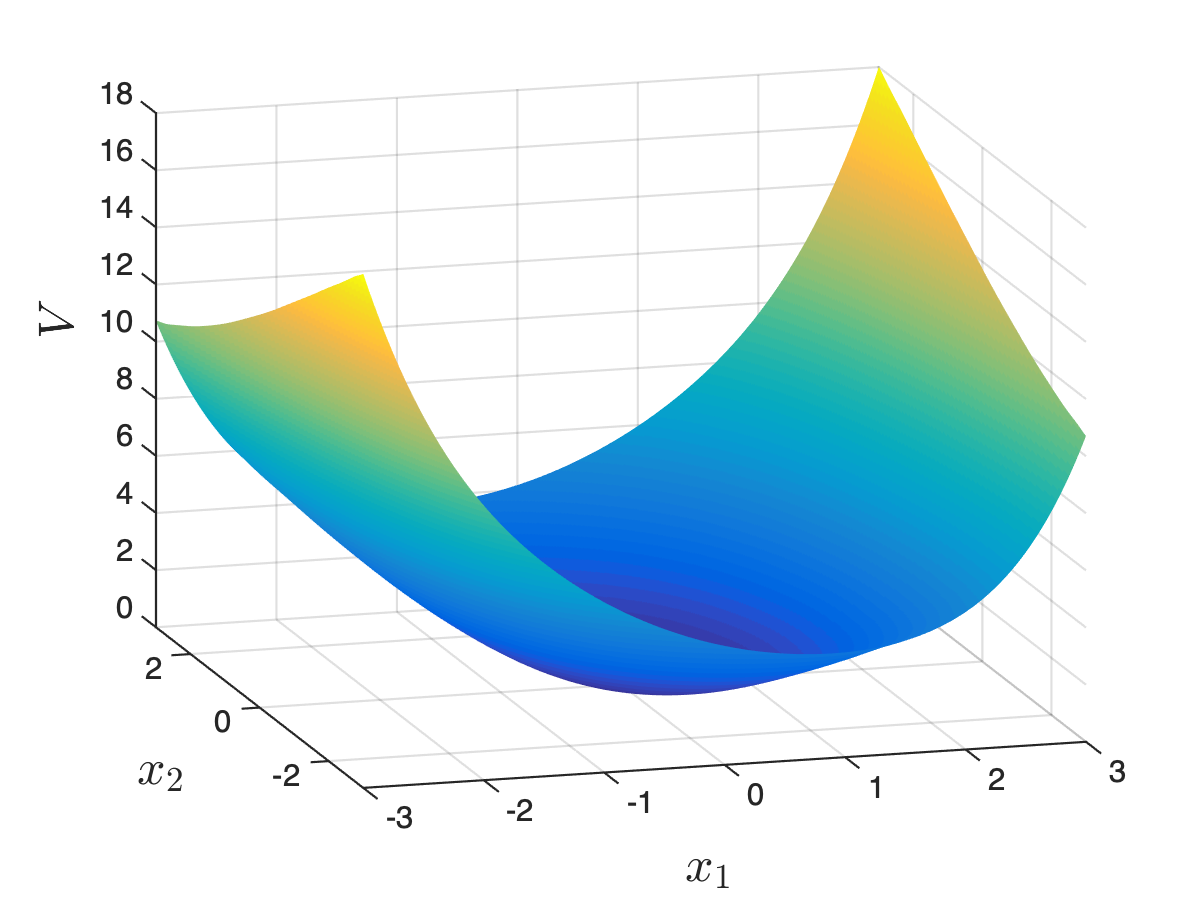}
	}
	\subfigure[${V}_{\ell_2} $ for  $N_d =40$]
	{
		\label{Fig8}
		\includegraphics[height=3.5cm,width=3.5cm]{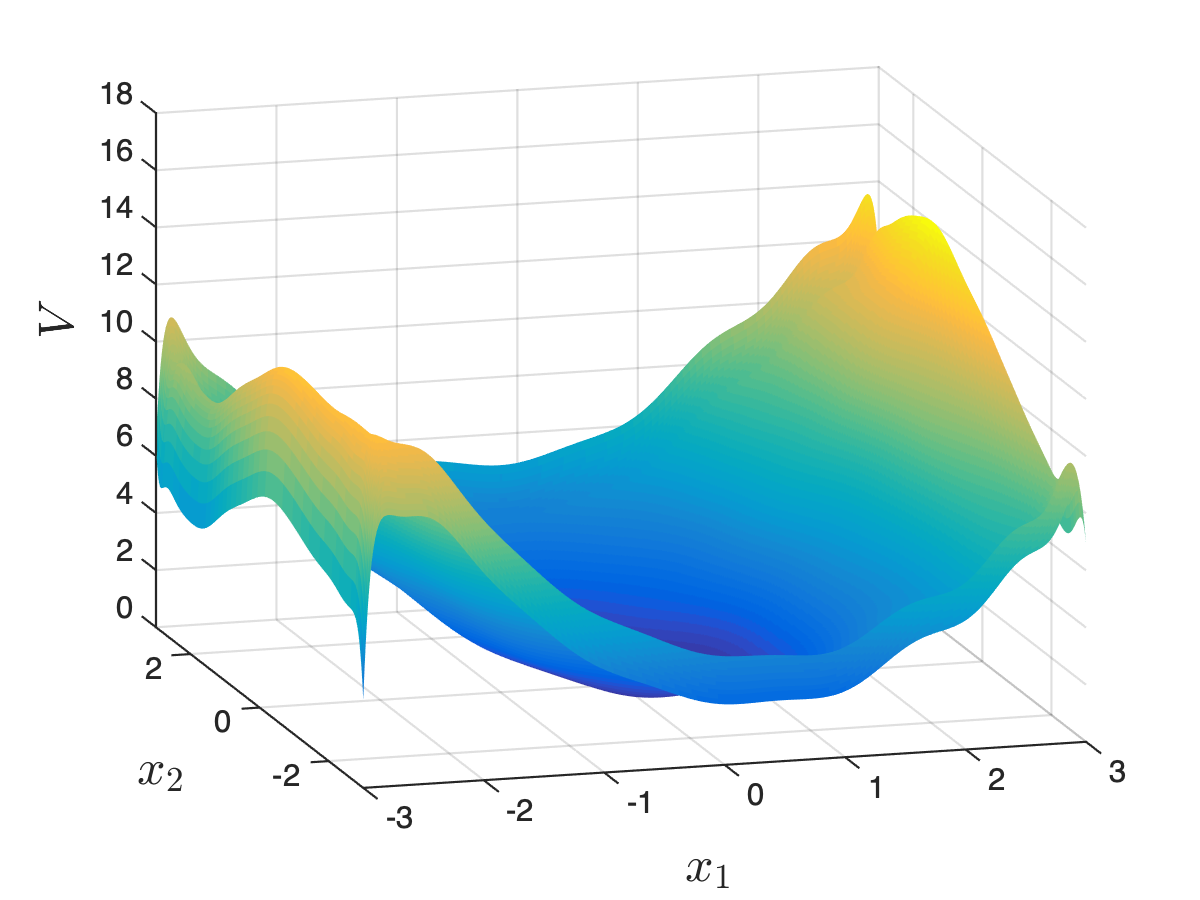}
	}
	\subfigure[$V_{\ell_2} $ for $N_d =120$]
	{
		\label{Fig9}
		\includegraphics[height=3.5cm,width=3.5cm]{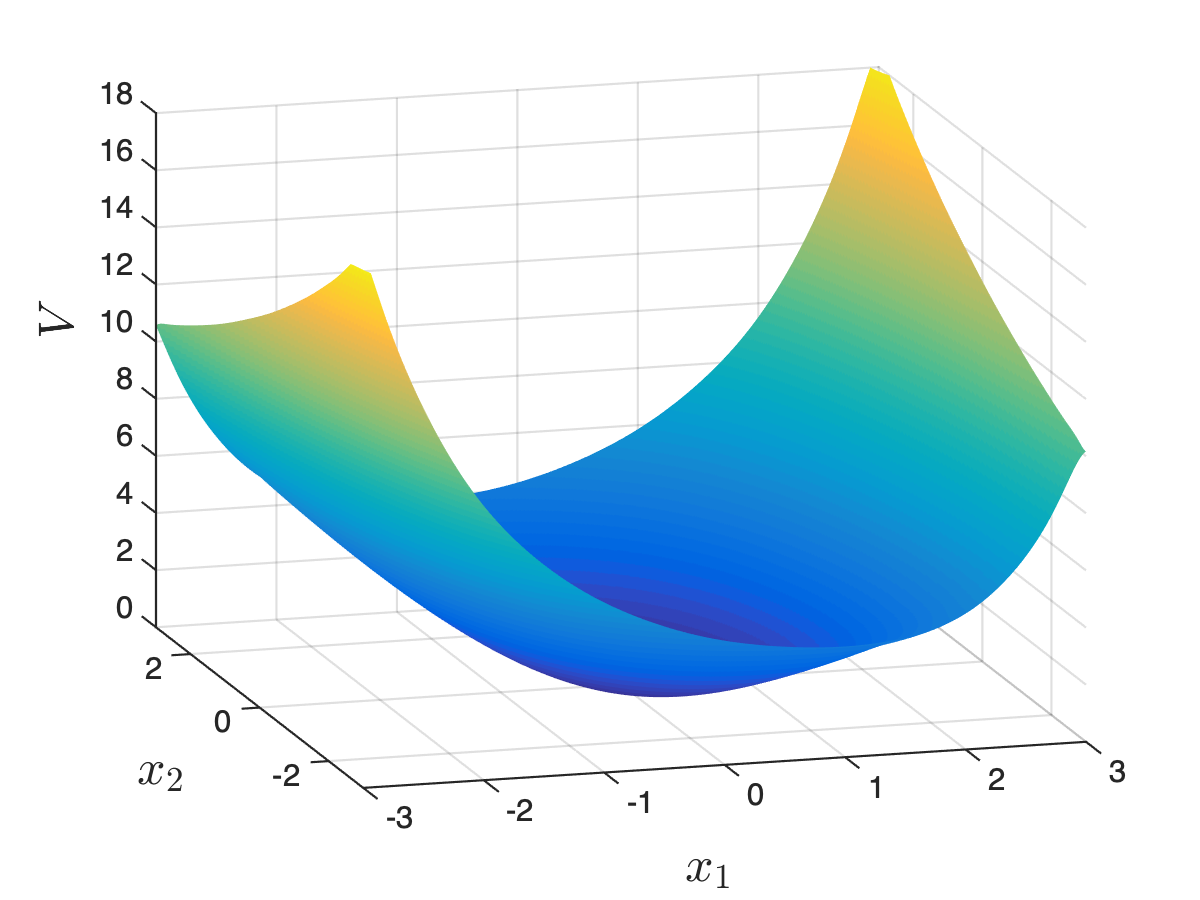}
	}
	\caption{Test 1. (a) Training and validation dataset (b) Sparse regression with gradient-augmented information and $N_d=40$ training points (c) Linear least squares without gradient-augmented information, $N_d=40$ (d) Linear least squares without gradient information with $N_d=120$. Sparse regression with gradient-augmented information provides an accurate reduced complexity approximation with fewer training samples.}
\end{figure}

To illustrate the approximation of the value function  $ V(\bx)$ by the different regression formulations,  we consider Figure \ref{Fig5}. This figure displays the scatter plot associated to  the training and validation data $\{\mathbf{x}^j,V^j\}^{N}_{j =1}$  with $ N =2000$.   Figure \ref{Fig6} shows the approximation on the bases of \eqref{e7} for the Legendre  polynomial basis, with $\lambda = 0.01$, $\alpha = 1$, and  $N_d = 40$. It clearly outperforms the approximation  on the basis of \eqref{e4} given in Figure \ref{Fig8}, again with $N_d=40$. To achive a similar result without including gradient information would require to increase the size of the training set to $N_d=120$, as shown in Figure \ref{Fig9}. Sparse regression with gradient-augmented information provides an accurate reduced complexity approximation with fewer training samples.

\begin{figure}[ht!]
	\centering
	\centering
	\subfigure[$\lambda =0.002$]
	{
		\includegraphics[height=4cm,width=4cm]{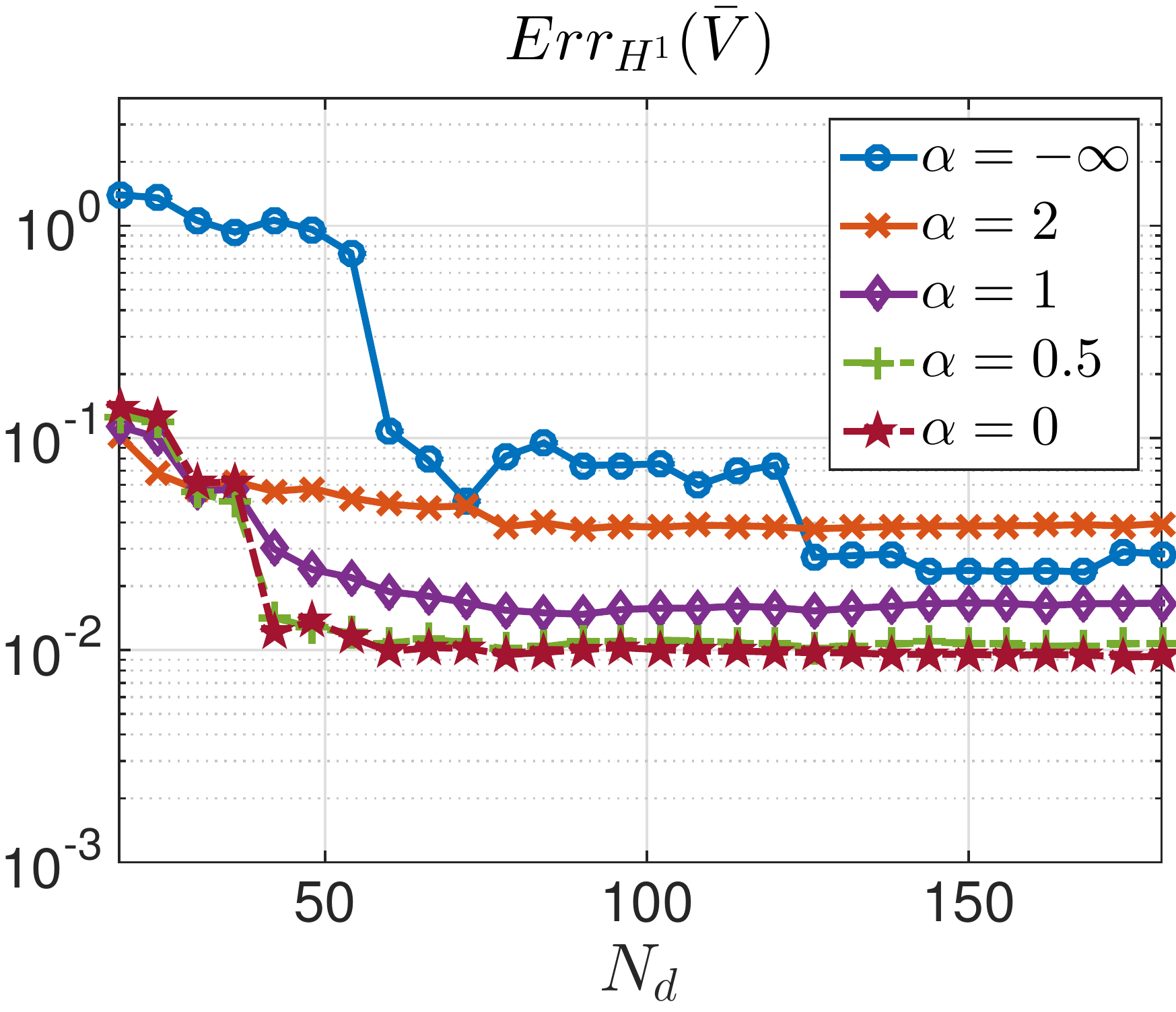}
	}
	\subfigure[$\lambda = 0.002$]
	{
		\includegraphics[height=4cm,width=4cm]{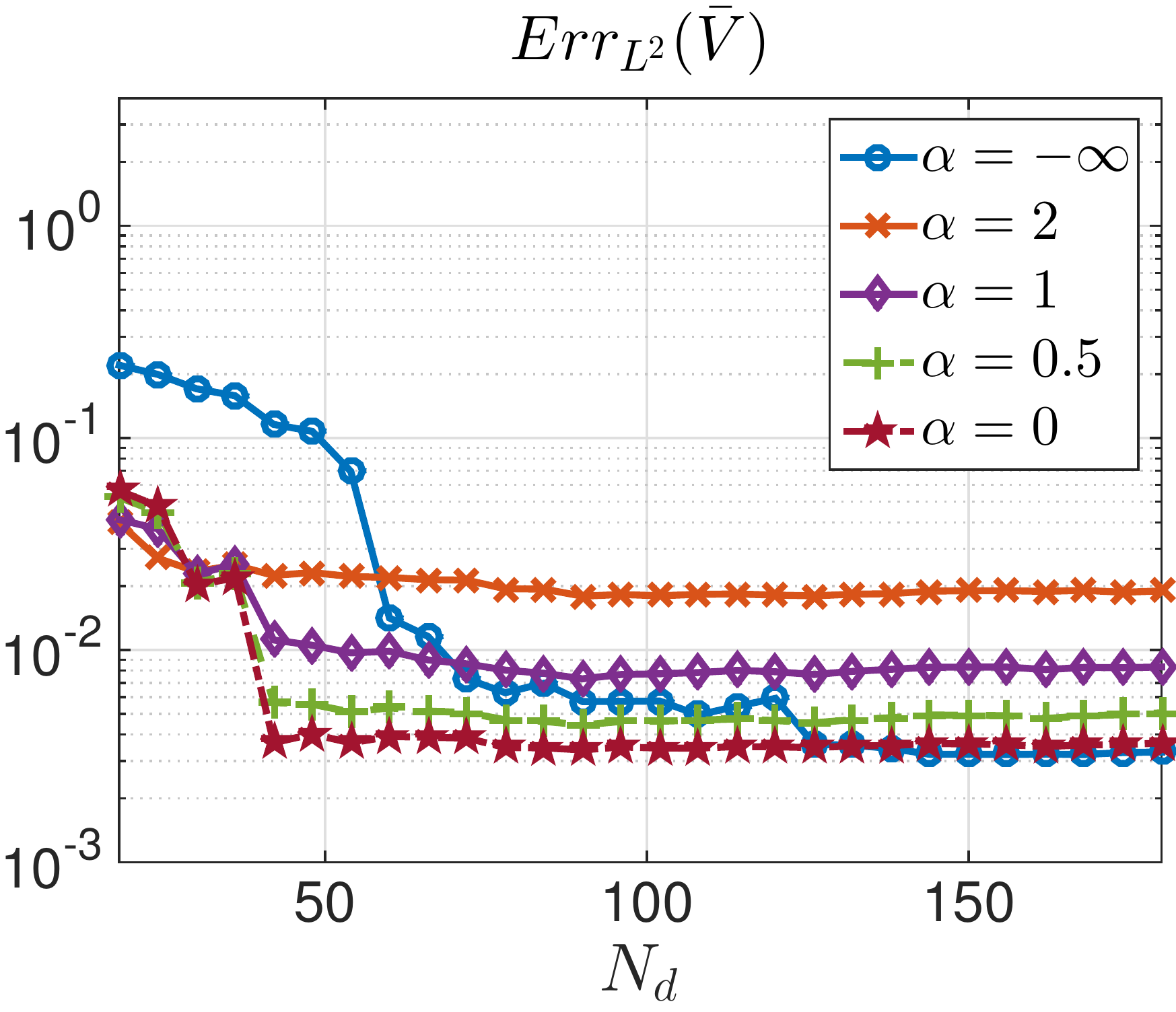}
	}
	\subfigure[$\lambda = 0.002$]
	{
		\includegraphics[height=4cm,width=4cm]{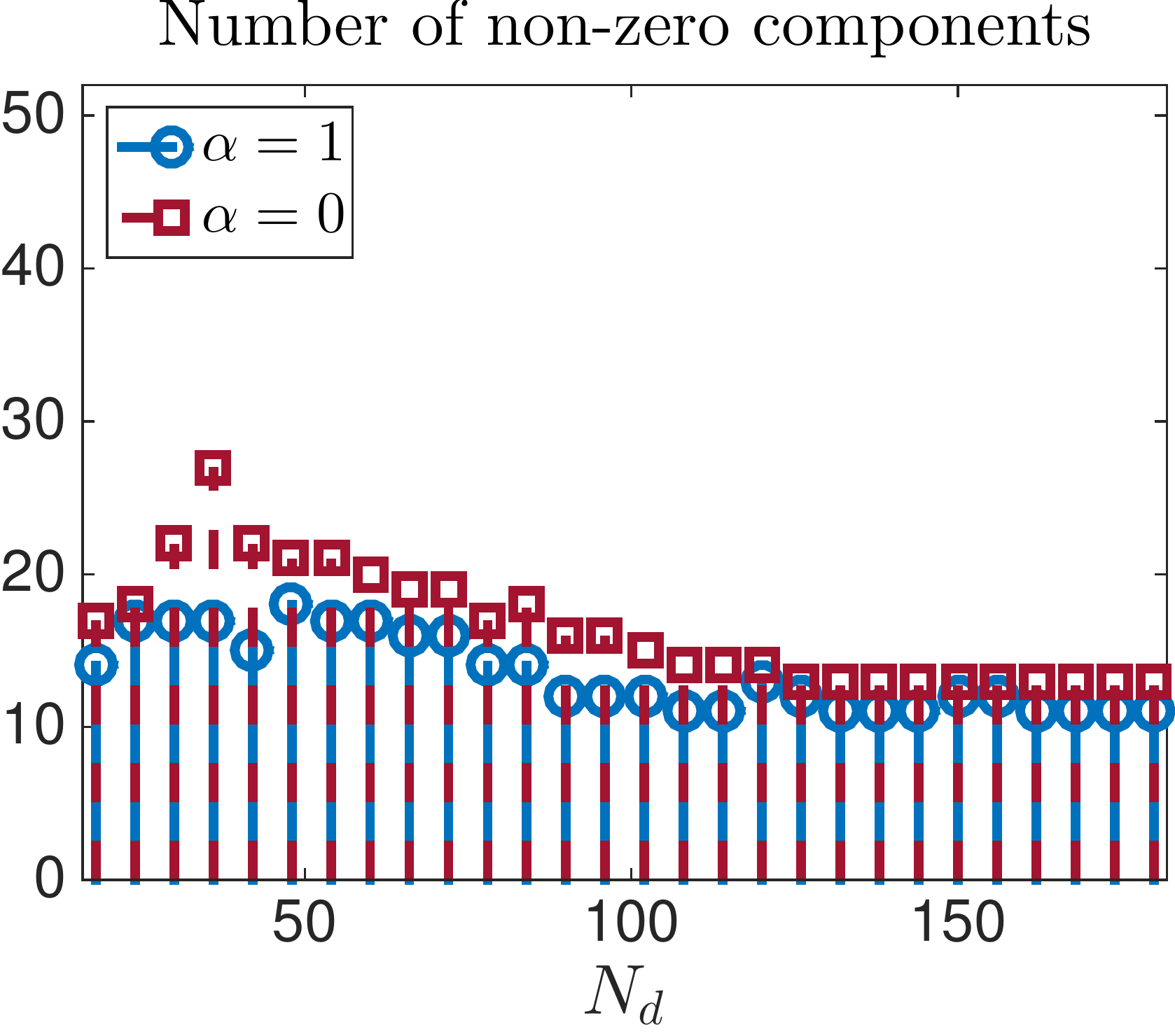}
	}
	\caption{Test 1. Numerical  results for the polynomial approximation without gradient information using the  Legendre  polynomial basis. {\color{black} Both $H^1$ and $L^2$ validation errors decrease as the number of training samples increases for different weights in the $\ell_1$ norm penalty. The sparse regression reduces the number of non-zero coefficients in the feedback law expansion.}}
	\label{Fig3}
\end{figure}

\begin{figure}[ht!]
	\centering
	\centering
	\subfigure[$\lambda =0.002$]
	{
		\includegraphics[height=4cm,width=4cm]{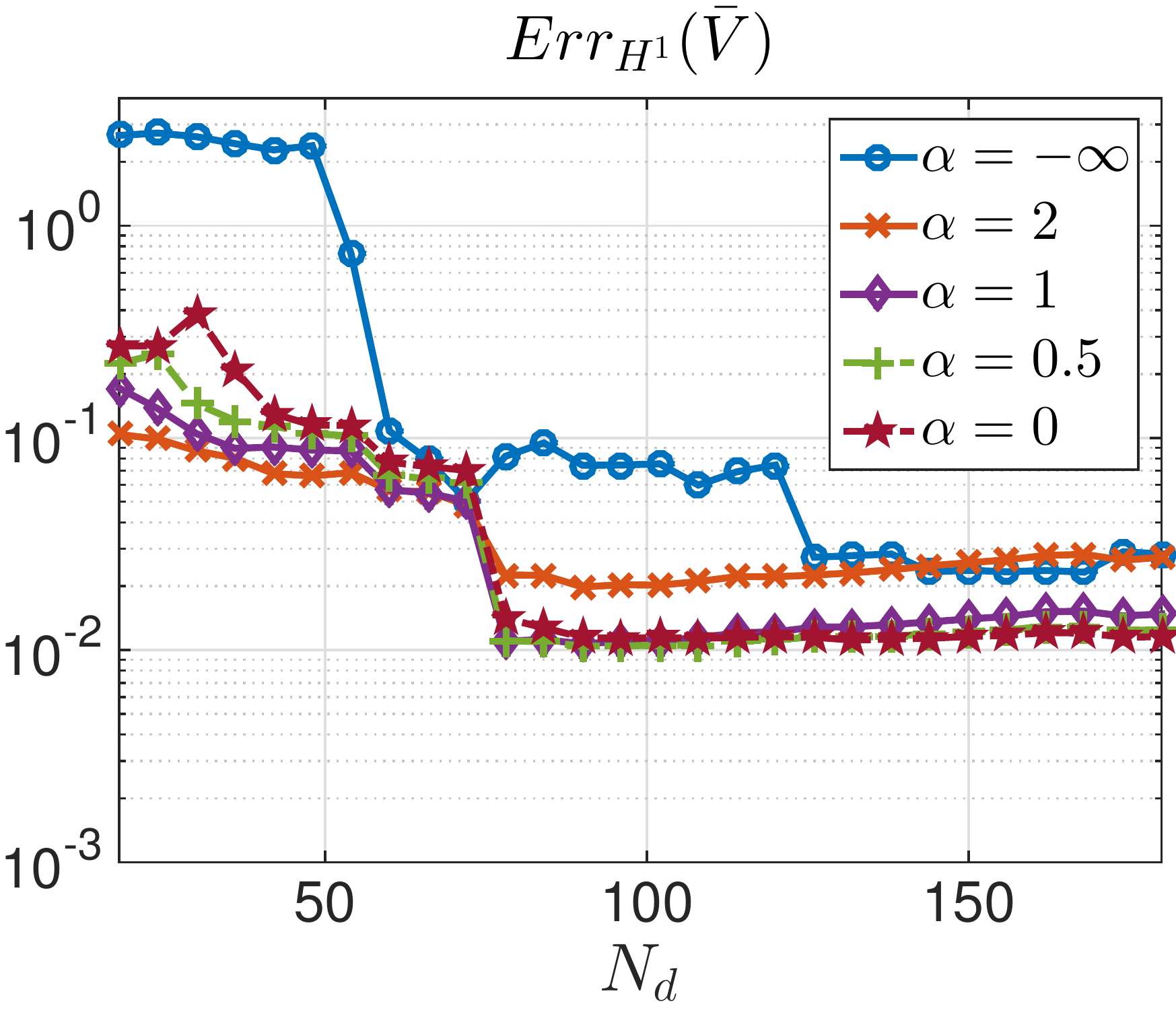}
	}
	\subfigure[$\lambda = 0.002$]
	{
		\includegraphics[height=4cm,width=4cm]{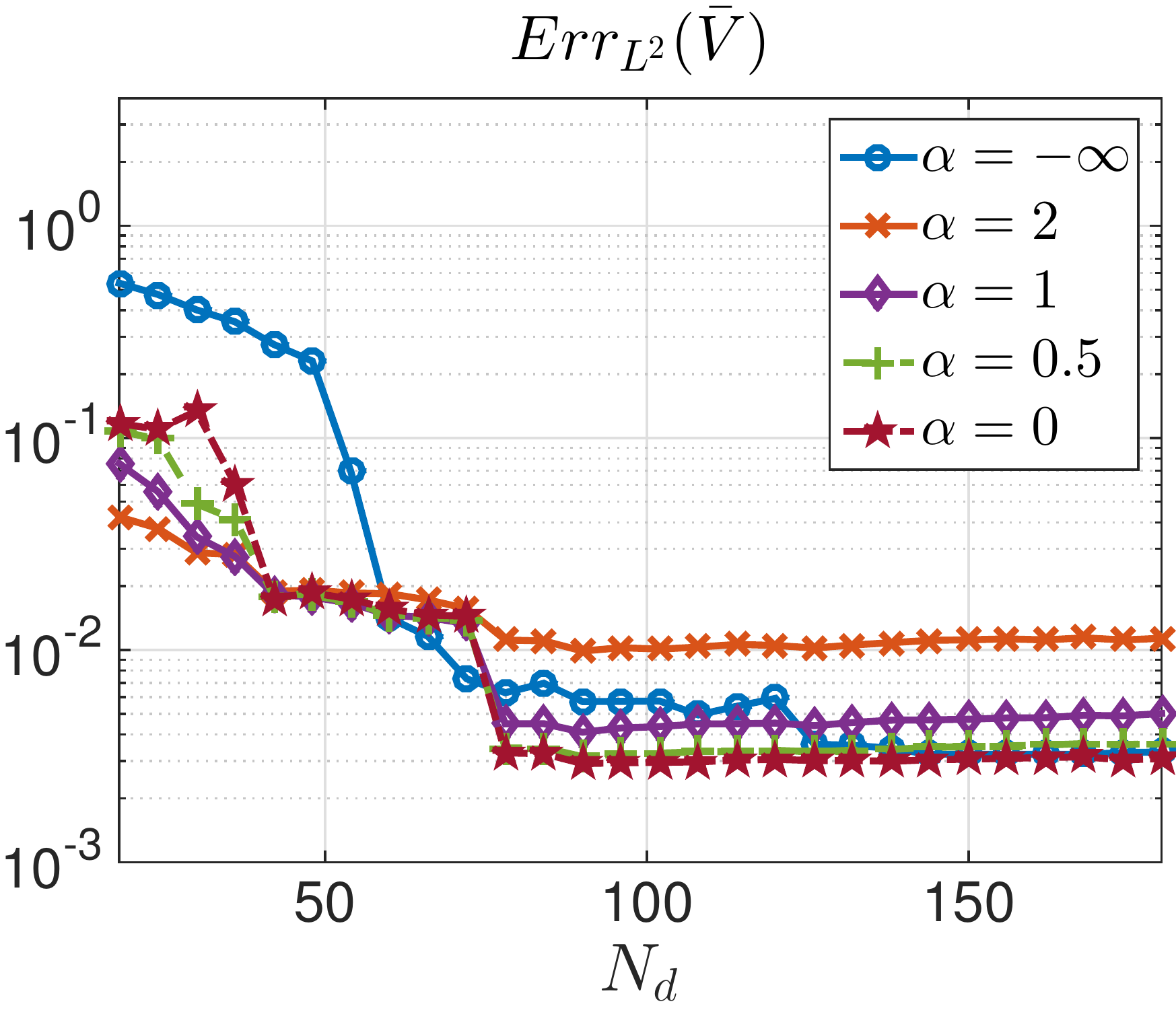}
	}
	\subfigure[$\lambda = 0.002$]
	{
		\includegraphics[height=4cm,width=4cm]{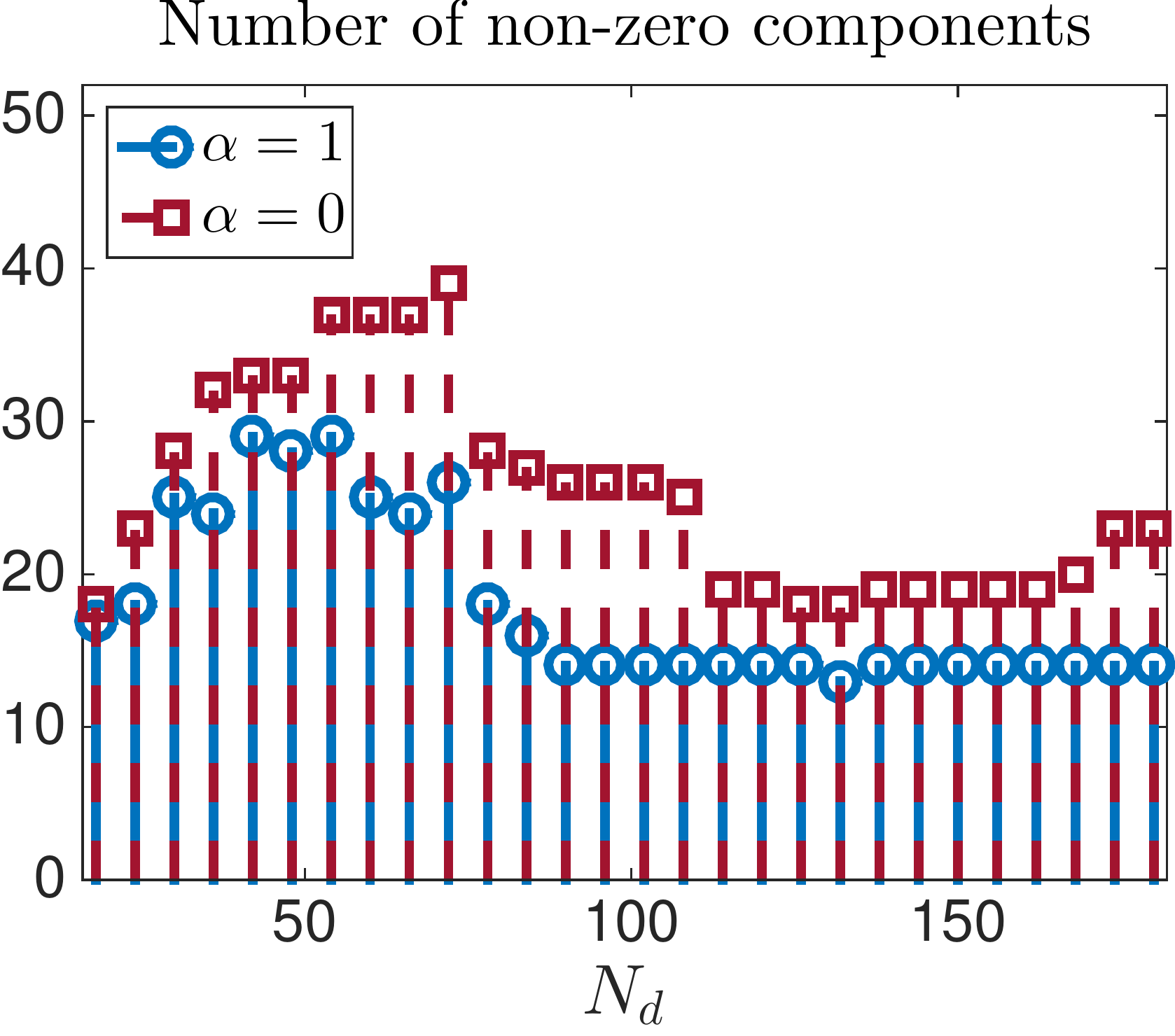}
	}
	\caption{Test 1. Numerical  results for the polynomial approximation without gradient informations using the  Chebyshev polynomial basis. {\color{black} Results are qualitatively similar to those in Figure \ref{Fig3} for a Legendre basis.}}
	\label{Fig4}
\end{figure}

\begin{figure}[ht!]
	\centering
	\subfigure[$\lambda =0.01$]
	{
		\includegraphics[height=4cm,width=4cm]{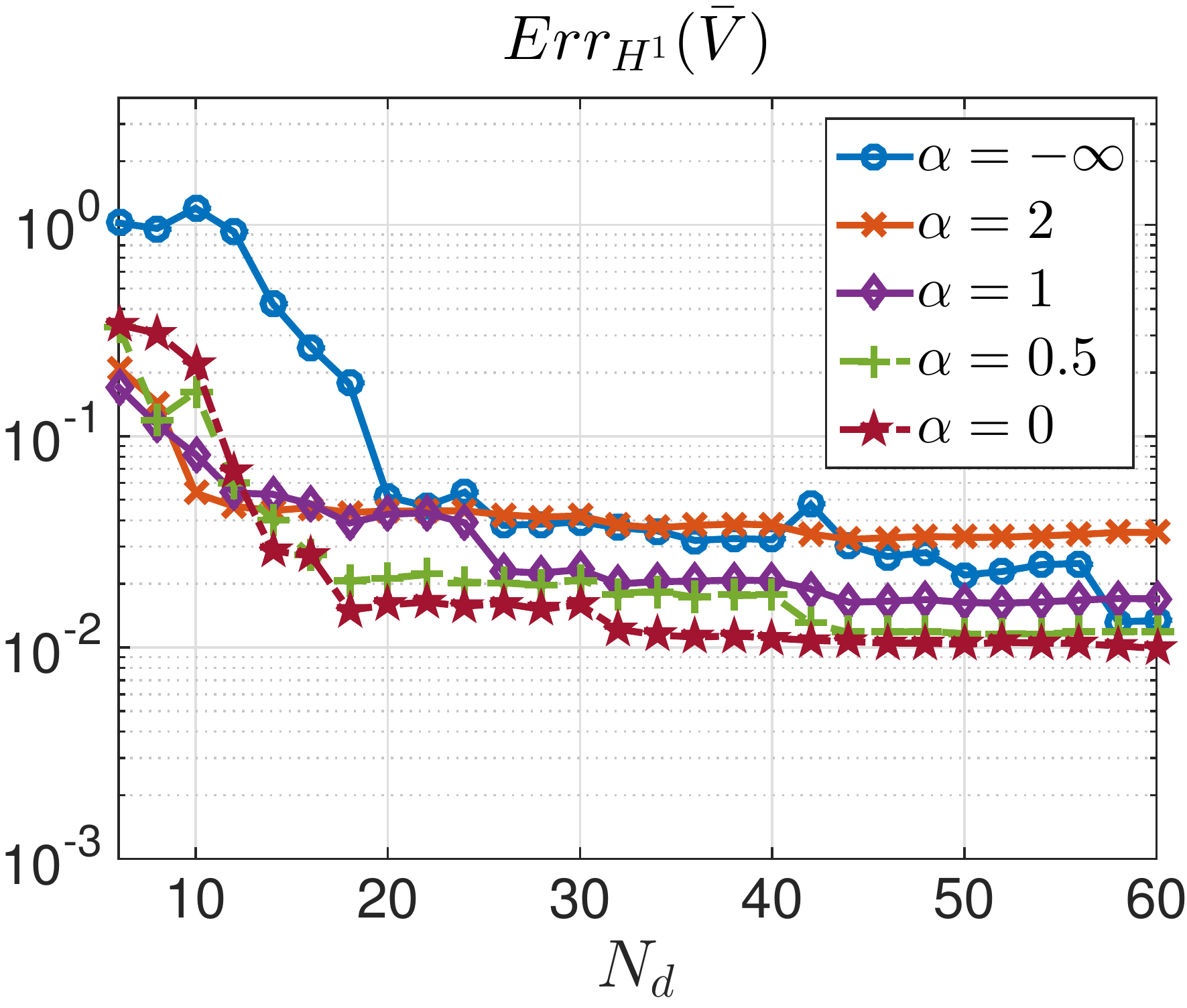}
	}
	\subfigure[$\lambda = 0.01$]
	{
		\includegraphics[height=4cm,width=4cm]{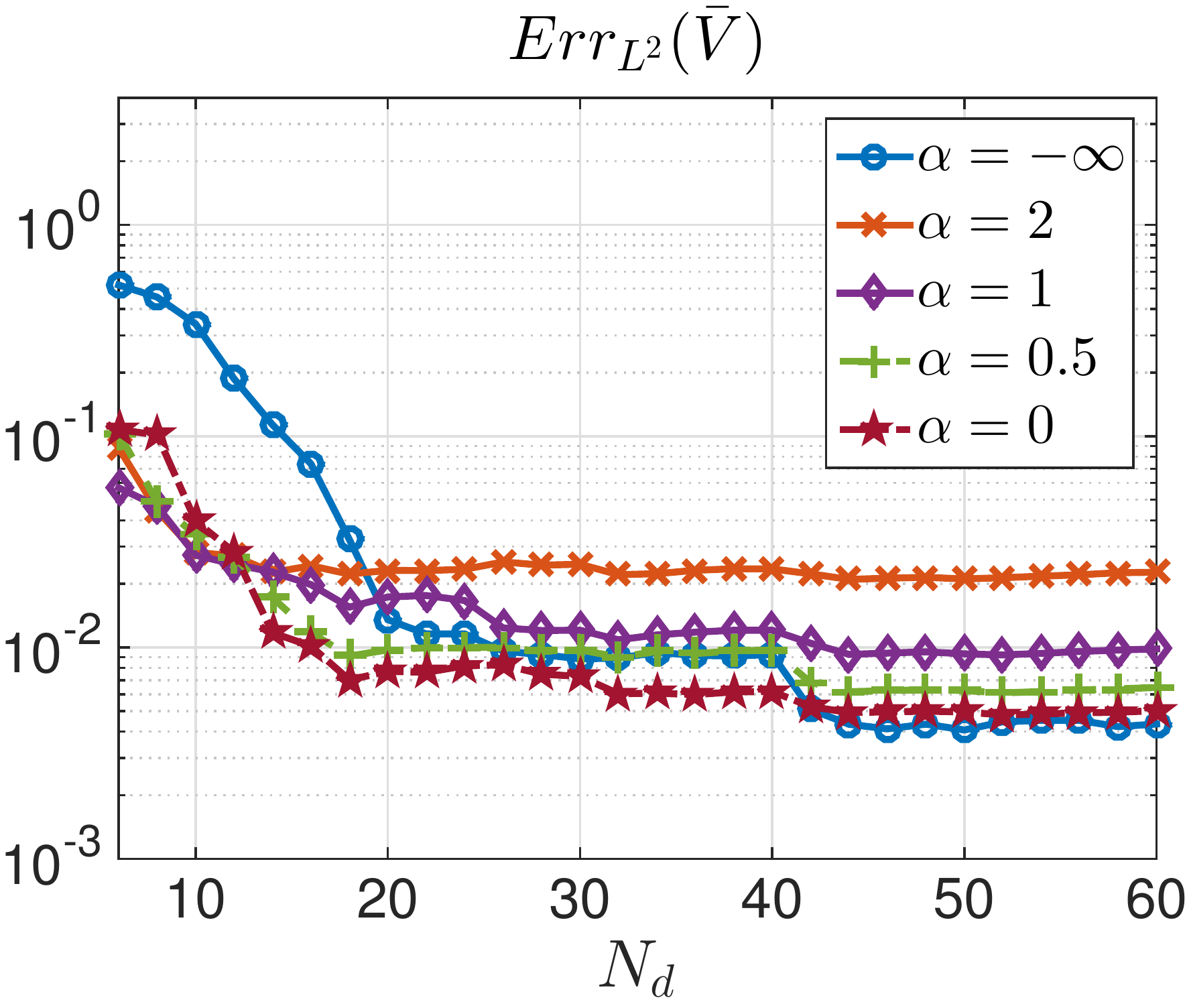}
	}
	\subfigure[$\lambda = 0.01$]
	{
		\includegraphics[height=4cm,width=4cm]{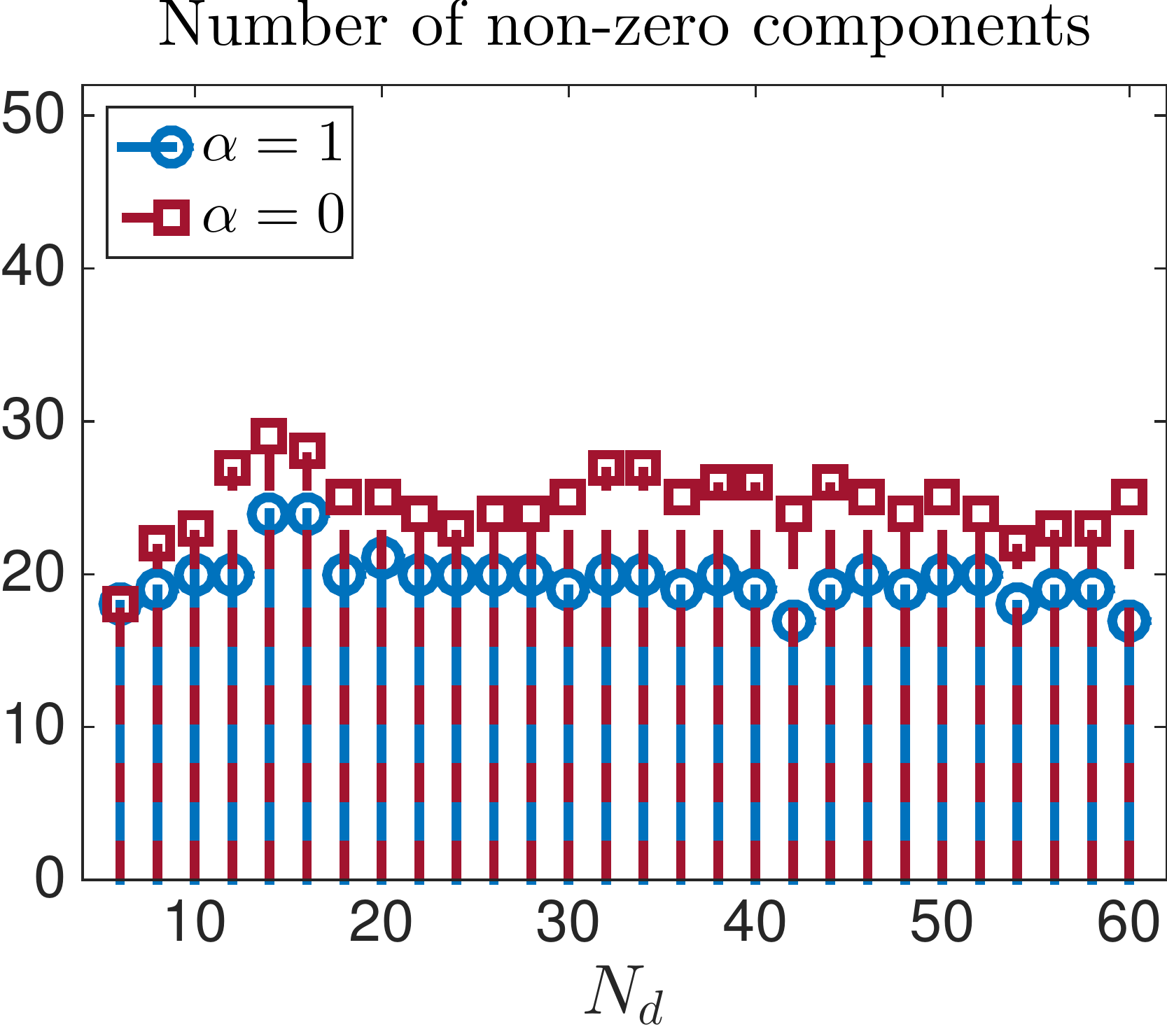}
	}
	\subfigure[$\lambda =0.02$]
	{
		\includegraphics[height=4cm,width=4cm]{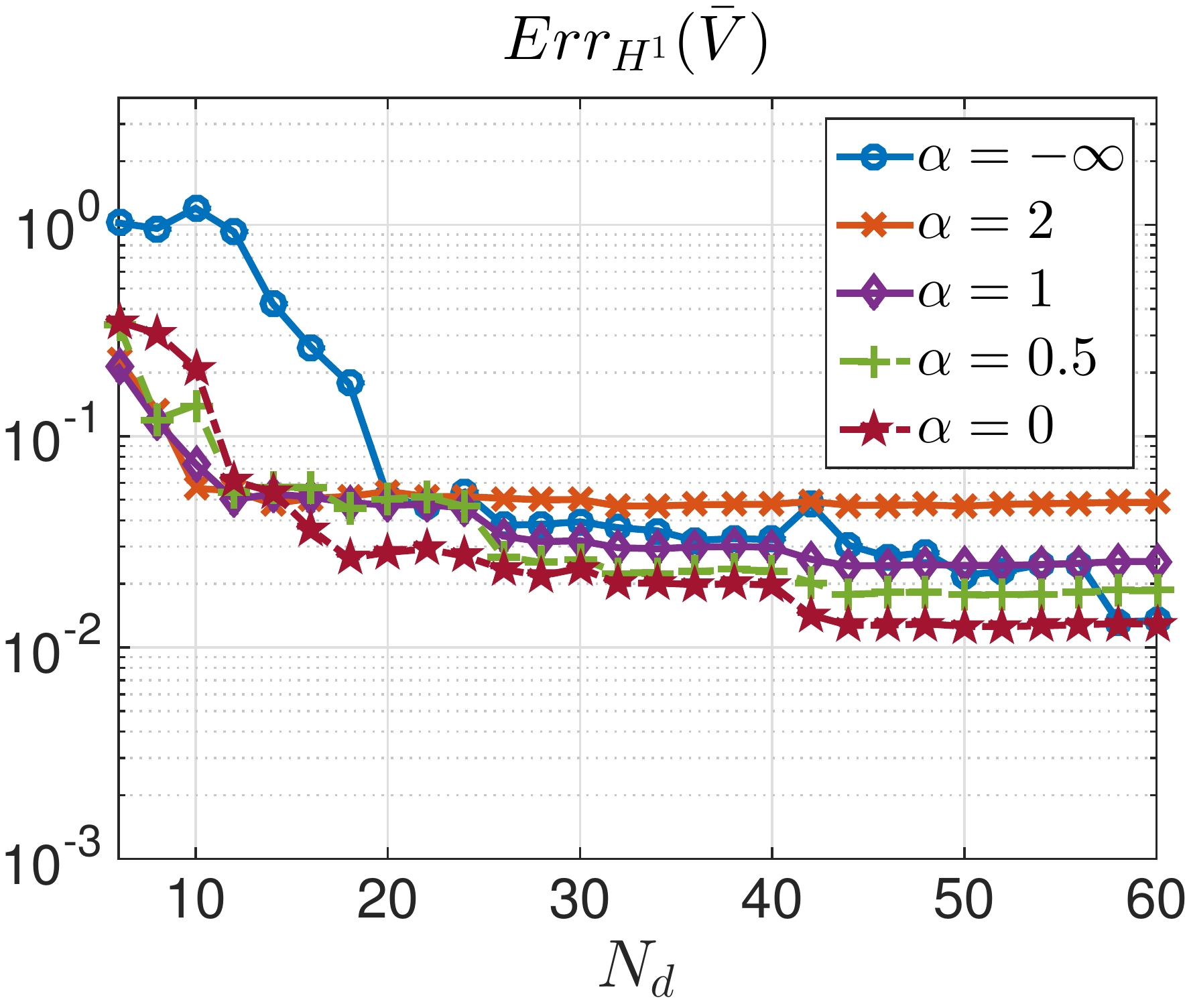}
	}
	\subfigure[$\lambda = 0.02$]
	{
		\includegraphics[height=4cm,width=4cm]{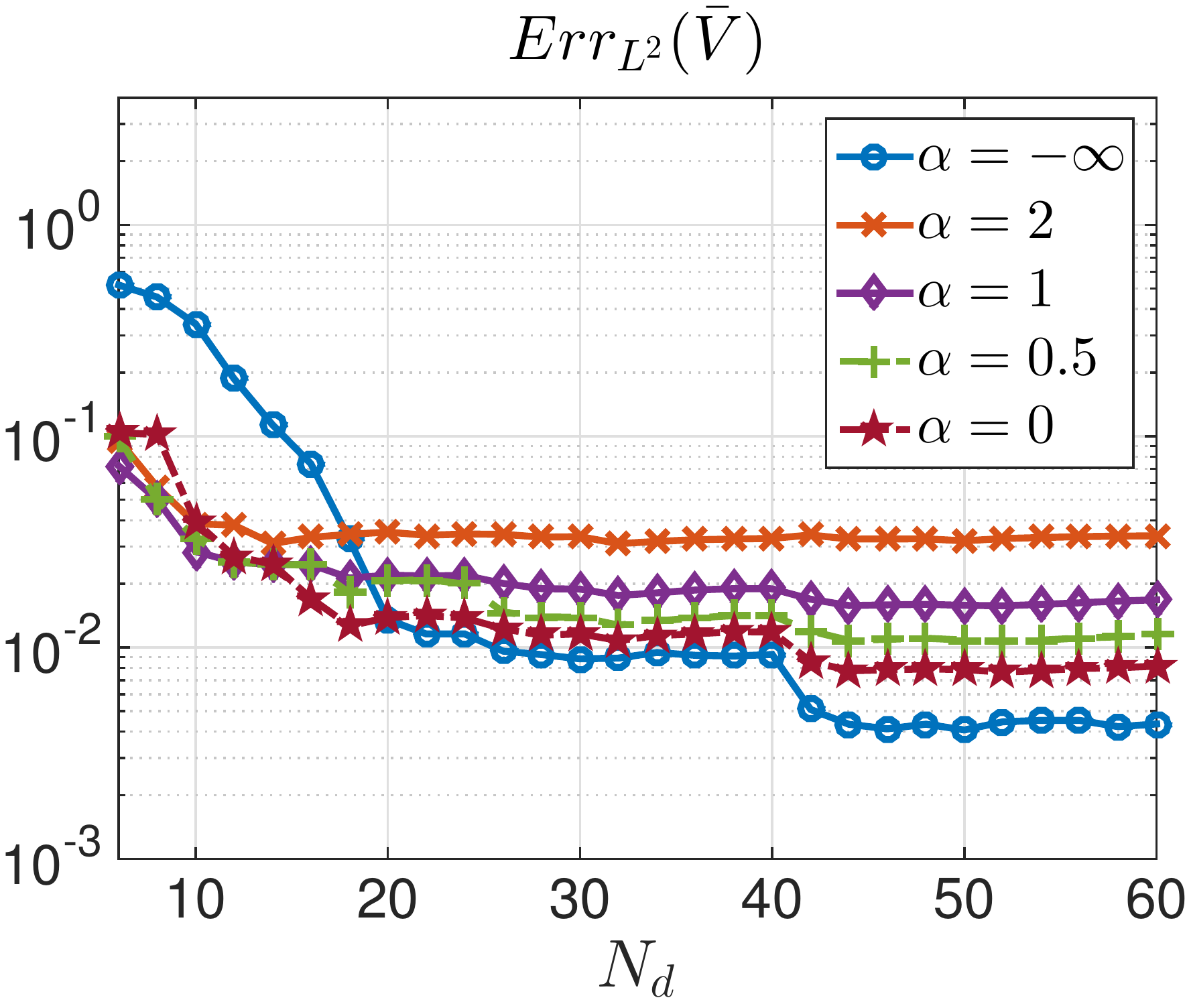}
	}
	\subfigure[$\lambda = 0.02$]
	{
		\includegraphics[height=4cm,width=4cm]{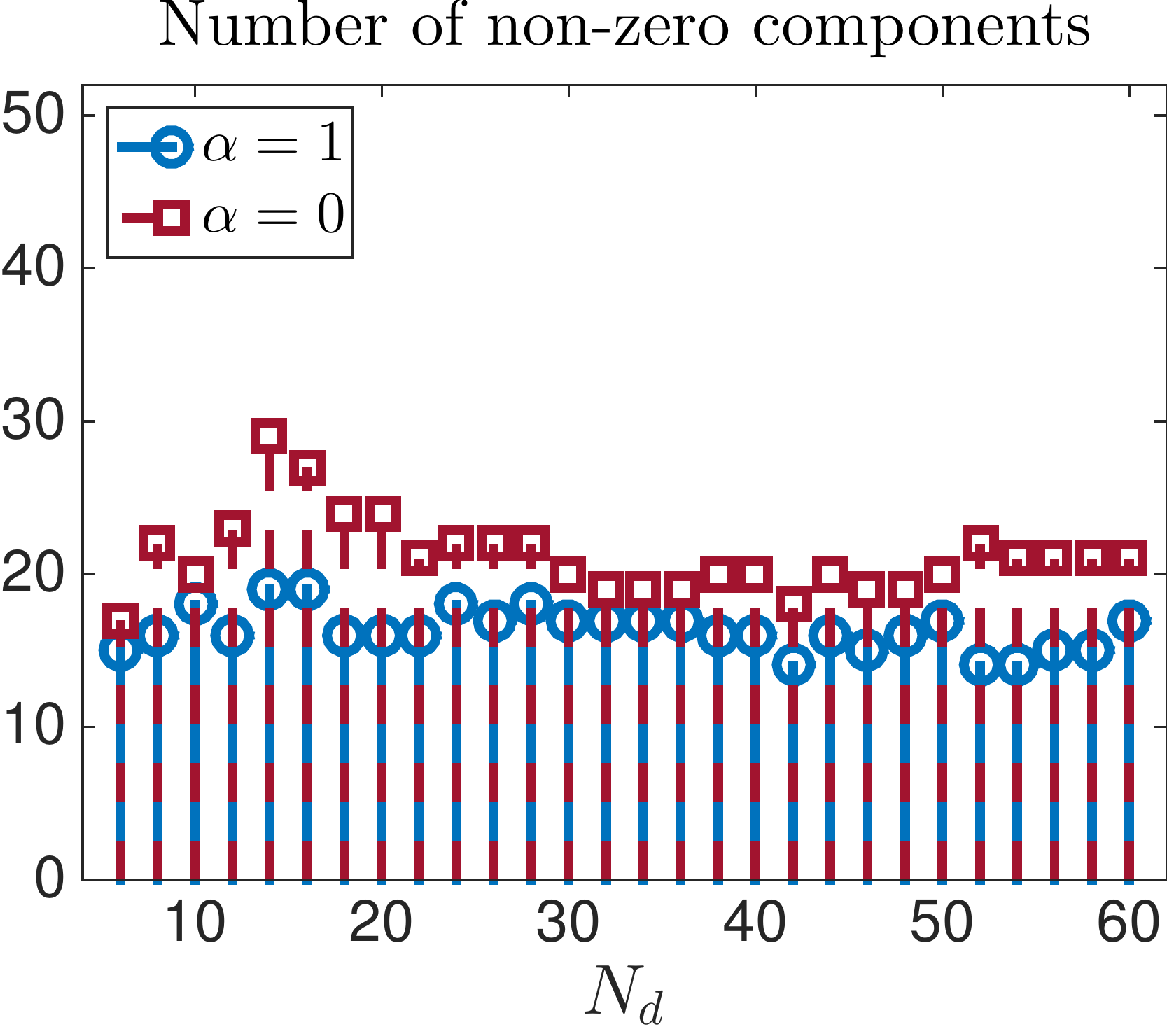}
	}
	\caption{Test 1. Numerical  results for the gradient-augmented polynomial approximation using the Legendre  polynomial basis. {\color{black} Similarly to Figures \ref{Fig3} and \ref{Fig4}, $H^1$ and $L^2$ validation errors decrease as the number of training samples increases. However, lower errors are reached with fewer samples due to the inclusion of gradient information, and without affecting the sparsity of the feedback law expansion.}}
	\label{Fig1}
\end{figure}

\begin{figure}[ht!]
	\centering
	\subfigure[$\lambda =0.01$]
	{
		\includegraphics[height=4cm,width=4cm]{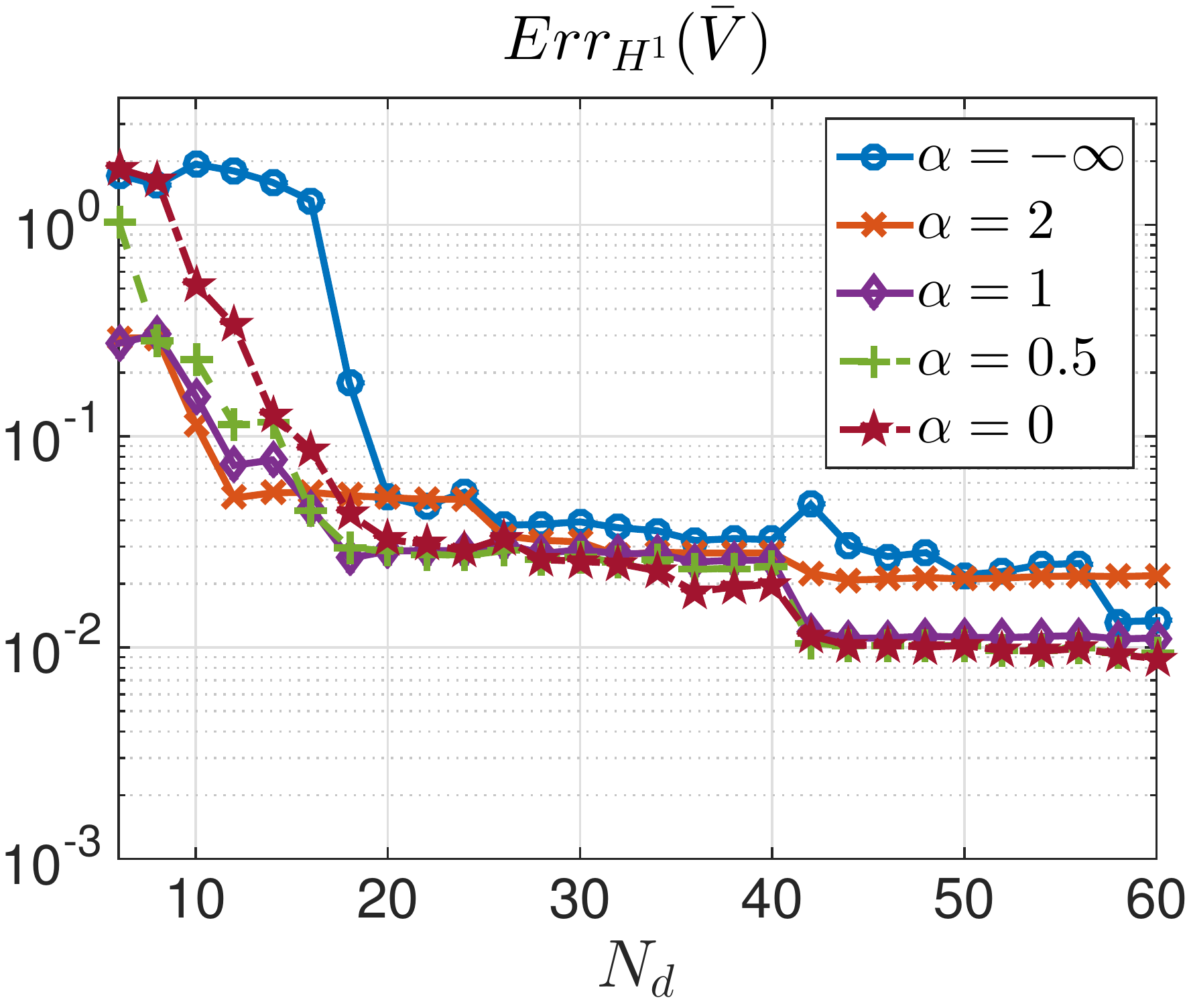}
	}
	\subfigure[$\lambda = 0.01$]
	{
		\includegraphics[height=4cm,width=4cm]{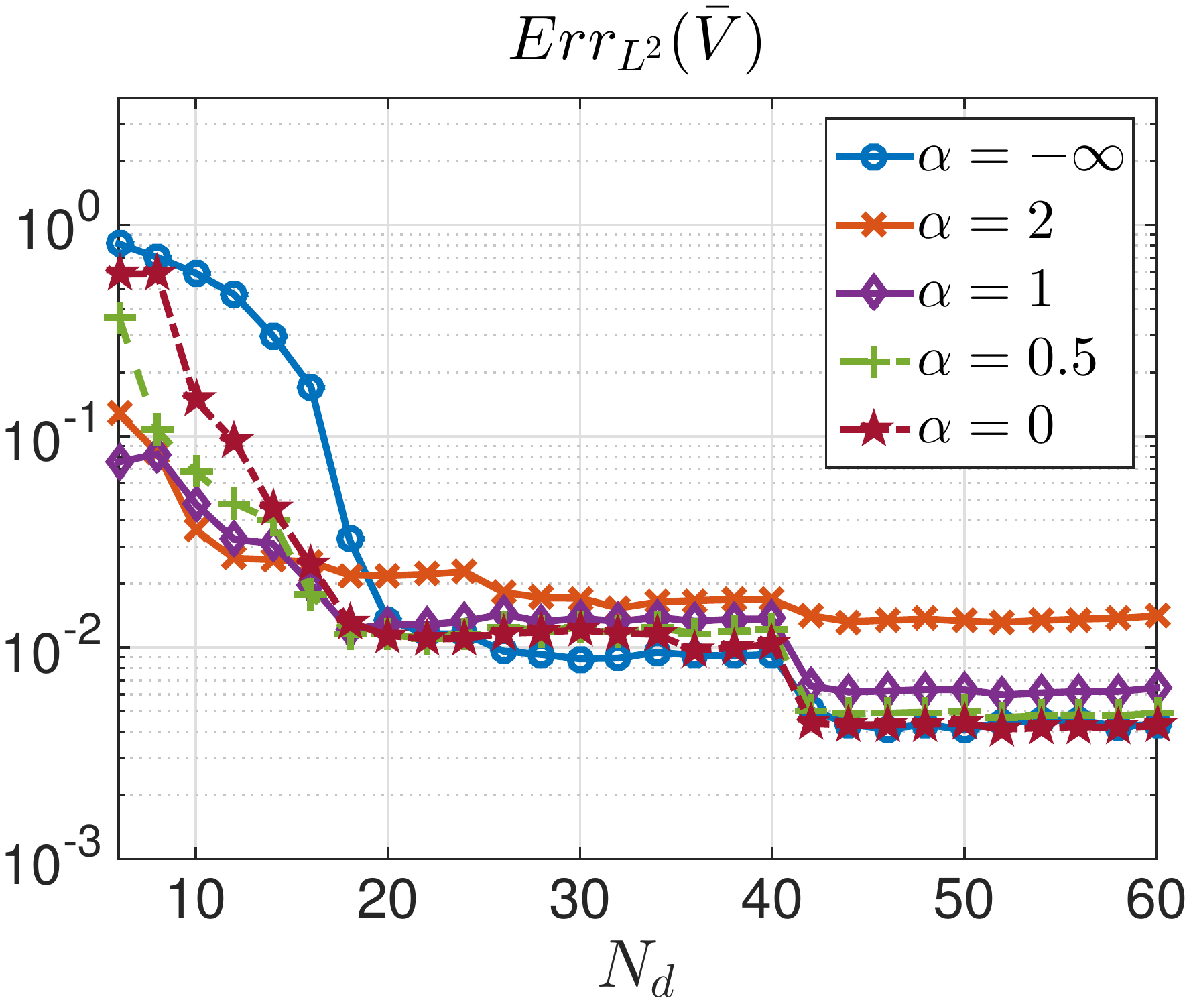}
	}
	\subfigure[$\lambda = 0.01$]
	{
		\includegraphics[height=4cm,width=4cm]{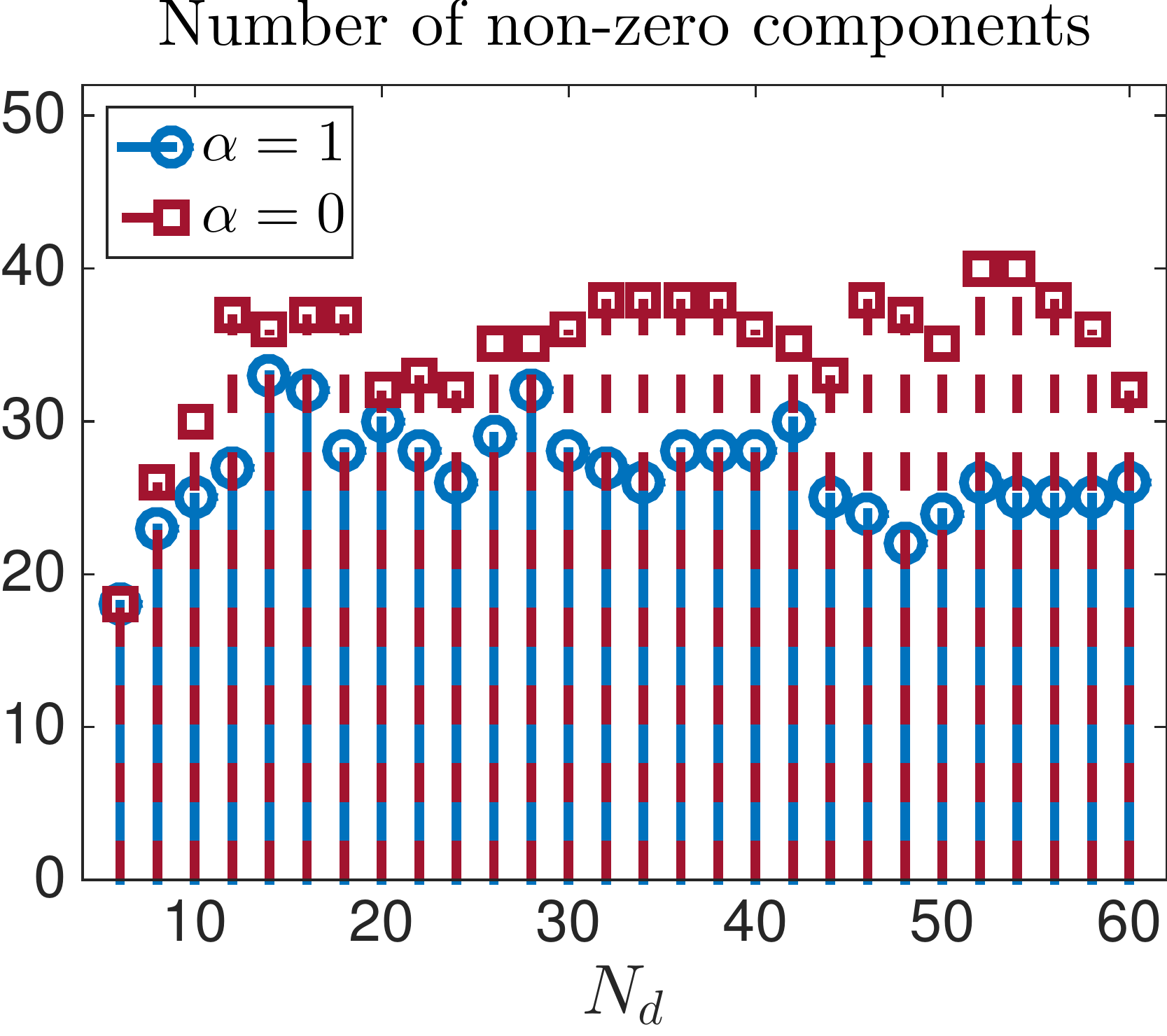}
	}
	\subfigure[$\lambda =0.02$]
	{
		\includegraphics[height=4cm,width=4cm]{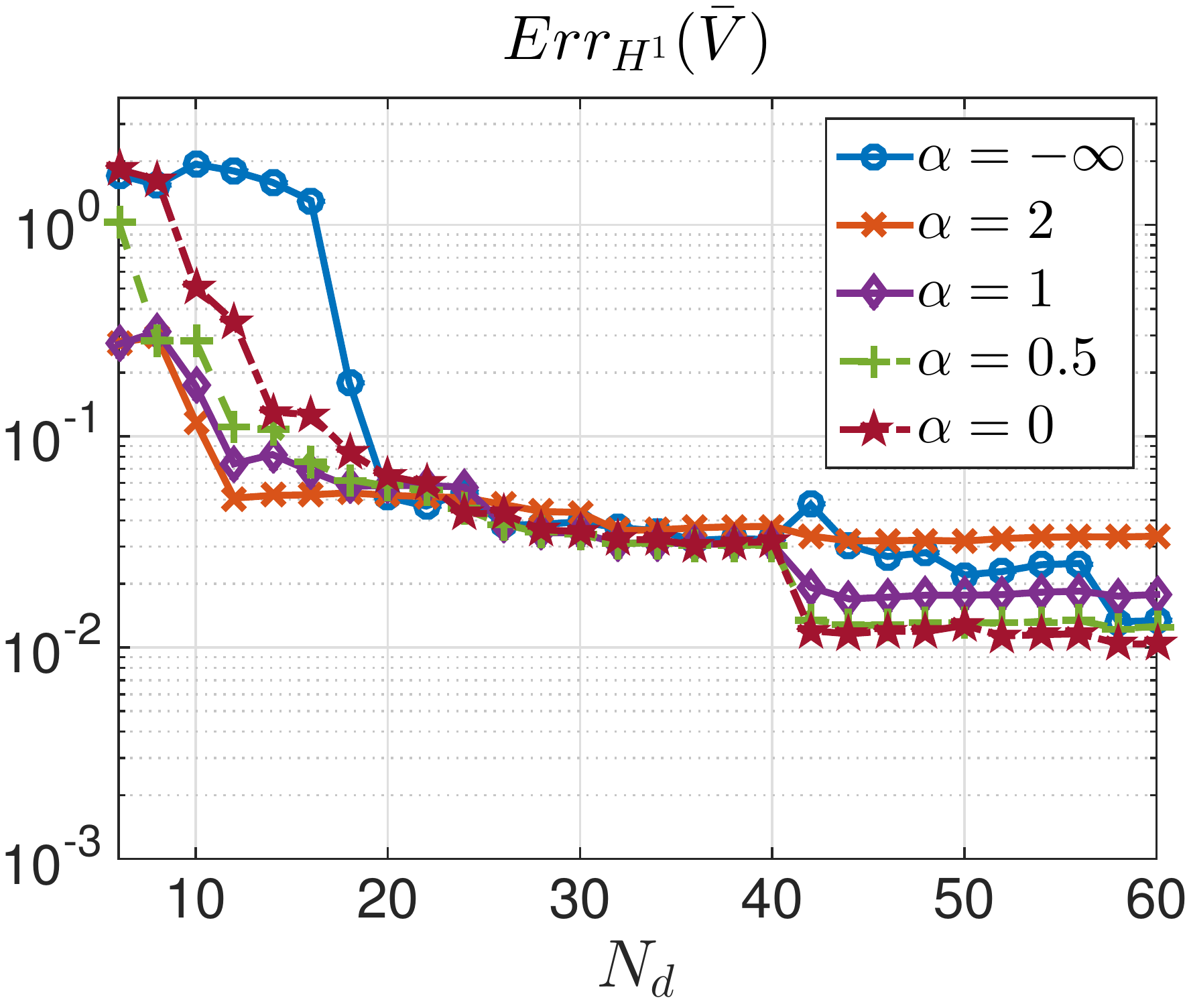}
	}
	\subfigure[$\lambda = 0.02$]
	{
		\includegraphics[height=4cm,width=4cm]{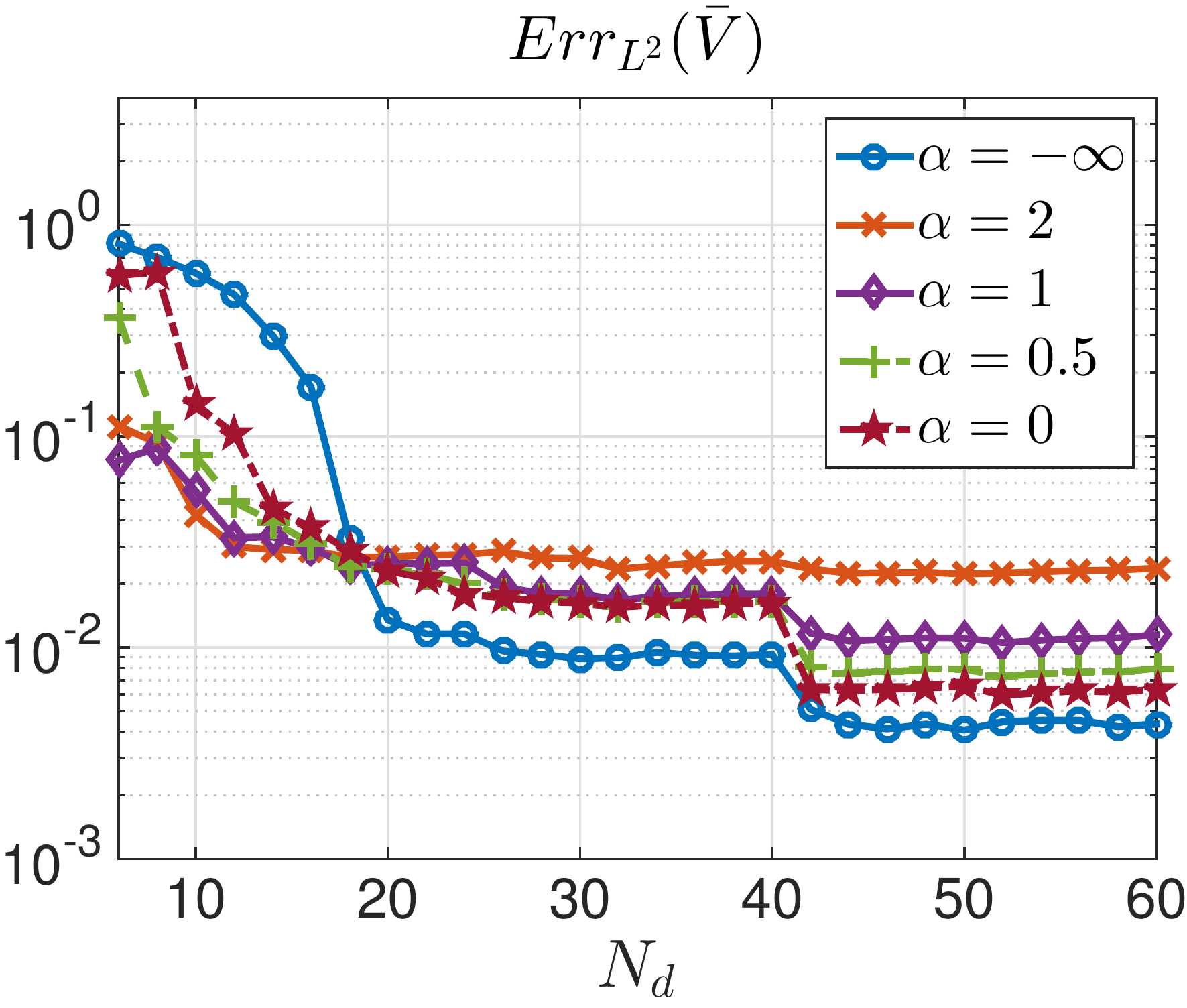}
	}
	\subfigure[$\lambda = 0.02$]
	{
		\includegraphics[height=4cm,width=4cm]{VanderPol_SLAlpha01Pol2-eps-converted-to}
	}
	\caption{Test 1. Numerical  results for the gradient-augmented polynomial approximation using the Chebyshev polynomial basis. {\color{black} Results are qualitatively similar to those in Figure \ref{Fig1} for the Legendre basis.}}
	\label{Fig2}
\end{figure}

We next turn to Figures \ref{Fig3}--\ref{Fig2} where the errors according to \eqref{eq:kk1} are plotted on a logarithmic scale ( $\log_{10} $) with respect to the number of samples $N_d$ used for training. Here  $|\mathcal{I}_{val}| = 1800$ validation samples were used.
For problem \eqref{e5} and  \eqref{e7} we chose the sparse penalty parameter $\lambda =  0.002$ and $\lambda = \{0.01, 0.02\}$,  respectively. Choosing $\lambda$ larger for    \eqref{e7} than for \eqref{e5} allows to approximately balance the contributions for the data and the regularization terms in the cost functionals of these two problems.
The cardinality of the non-zero coefficients of  $\theta_{\ell_1}$ and $\bar{\theta}_{\ell_1}$ is determined by defining  components as nonzero if its absolute value is bigger than double machine extended precision $10^{-20}$.  Let us next make some observations on these results.

As expected, the error decreases with the training size $N_d$, up to a certain threshold.
The best possible fit for the chosen order $s=16$ ($q=52$) of the polynomial approximation is reached at about $N_d=50$ and $N_d=75$ for Legendre and Chebyshev polynomials, respectively, without the use of gradient information, see Figures \ref{Fig3} and \ref{Fig4}.  These error levels are reached much earlier when we include gradient information, as shown in Figures \ref{Fig1} and \ref{Fig2}. The influence of the $\ell_1$ weights, expressed in terms of $\alpha$, is not very pronounced. Note that $\alpha=-\infty$ corresponds to no regularisation, whereas $\alpha=0$ corresponds to a constant weight. In the case that  $N_d \lll 52$, the system is highly under-determined for $\alpha=-\infty$, which goes along with a large error. For small $N_d$, the choice $\alpha=2$ can be favoured over the choice $\alpha=0$, with the latter giving best results for $N_d$ sufficiently large.

In the last column of these plots the cardinality of nonzero coefficients for $\bar \theta_{\ell_1}$ is depicted. It typically increases with $N_d$ up to a certain threshold, and roughly stays constant thereafter, at less than 50 percent of the total number of free coefficients. Increasing $\lambda$
promotes sparsity, as expected, see Figures \ref{Fig1} -\ref{Fig2}, (c) and (f).

Figures \ref{Fig3}  and \ref{Fig1}  correspond to the  Legendre  polynomial basis, while  Figures \ref{Fig4}  and \ref{Fig2} are obtained with  Chebyshev polynomials. The results are quite
similar  in terms of the asymptotic (with respect to $N_d$) behavior of the errors. The number of non-zero components is higher for the Chebyshev than for the Legrendre polynomial expansion.
By comparing Figures \ref{Fig1} (resp. Figure \ref{Fig2} )  with \ref{Fig3}
(resp. Figure \ref{Fig4} ), we can see that for obtaining the same precision of approximation we need to consider more samples.  As expected, for the  case of gradient-augmented approximation we obtained better results for the $H_1$-errors for small $N_d$.

\begin{figure}[ht!]
	\centering
	\subfigure[State]
	{
		\label{Fig10}
		\includegraphics[height=4cm,width=4cm]{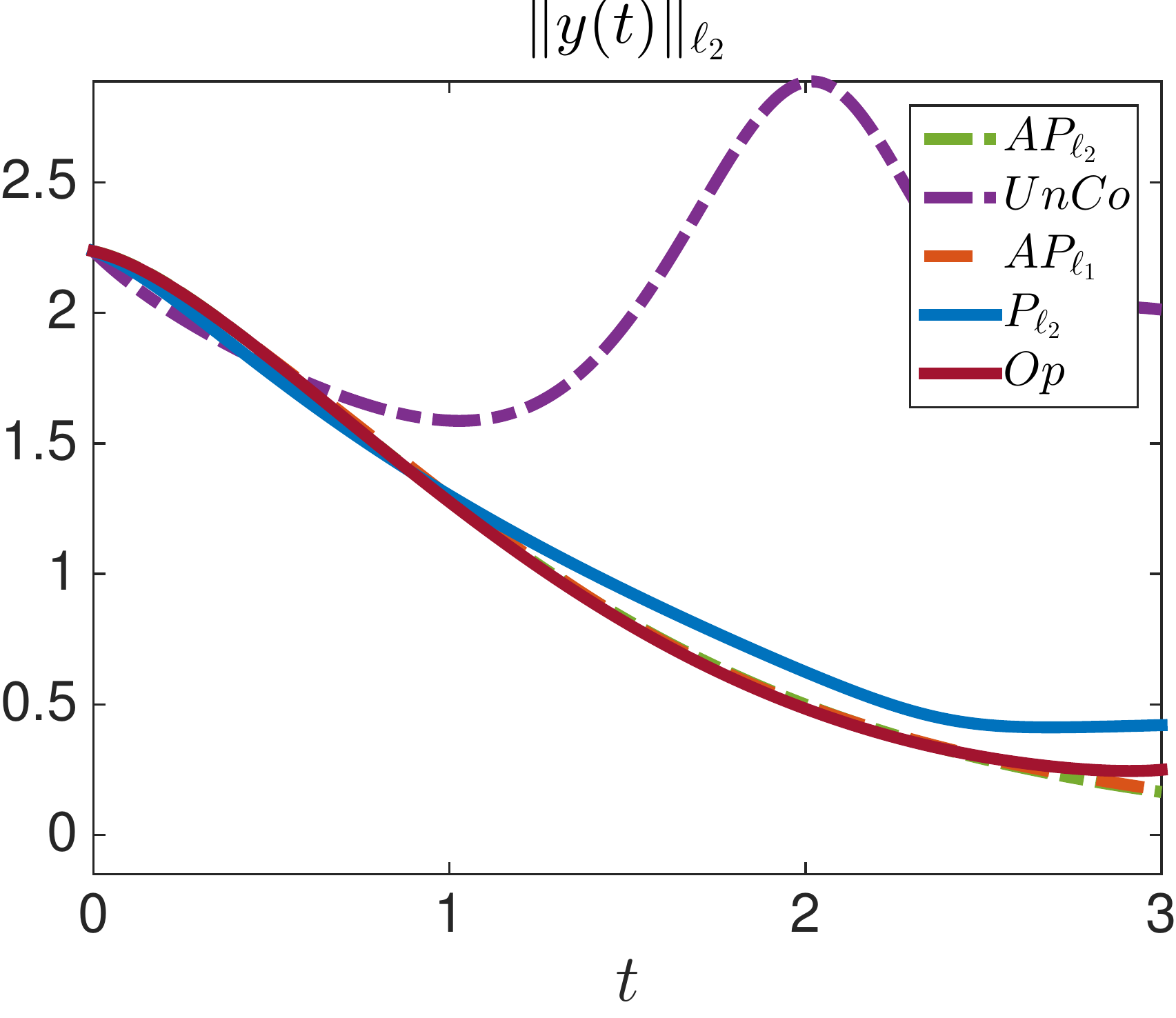}
	}
	\subfigure[Control]
	{
		\label{Fig12}
		\includegraphics[height=4cm,width=4cm]{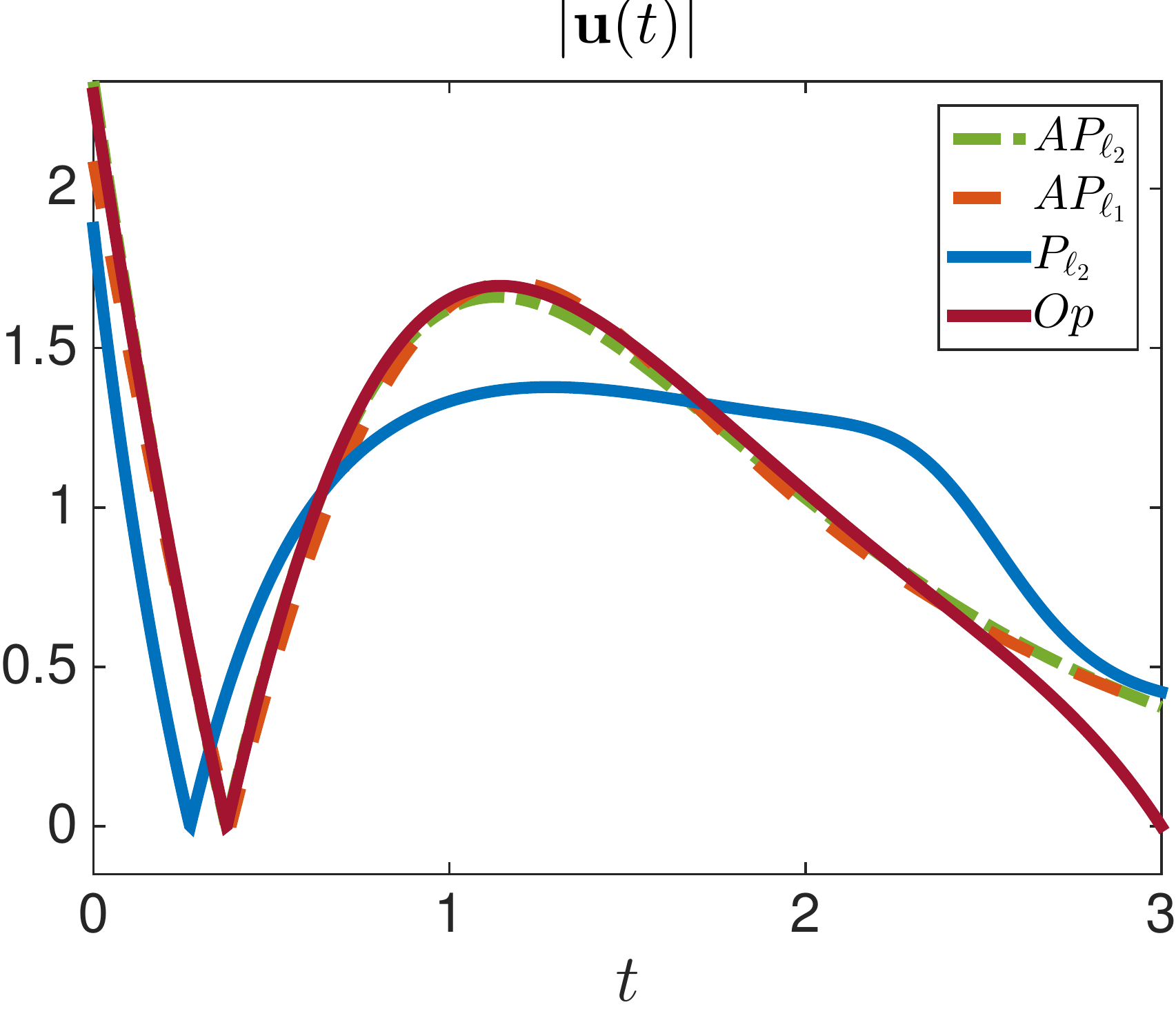}
	}
	\caption{Test 1. Evolution of  $\| y(t) \|_{\ell_2}$ and  $|\mathbf{u}(t)|$ for $\bx = (2,-1)$. Here \textsl{UnCo} stands for the uncontrolled trajectory, and \textsl{Op} refers to the exact optimal trajectory. {\color{black} We observe that the optimal feedback law obtained from the gradient-augmented sparse polynomial regression follows the true optimal trajectory, unlike the feedback law computed without gradient information.}}
\end{figure}

\begin{figure}[ht!]
	\centering
	\subfigure[State]
	{
		\label{Fig11}
		\includegraphics[height=4cm,width=4cm]{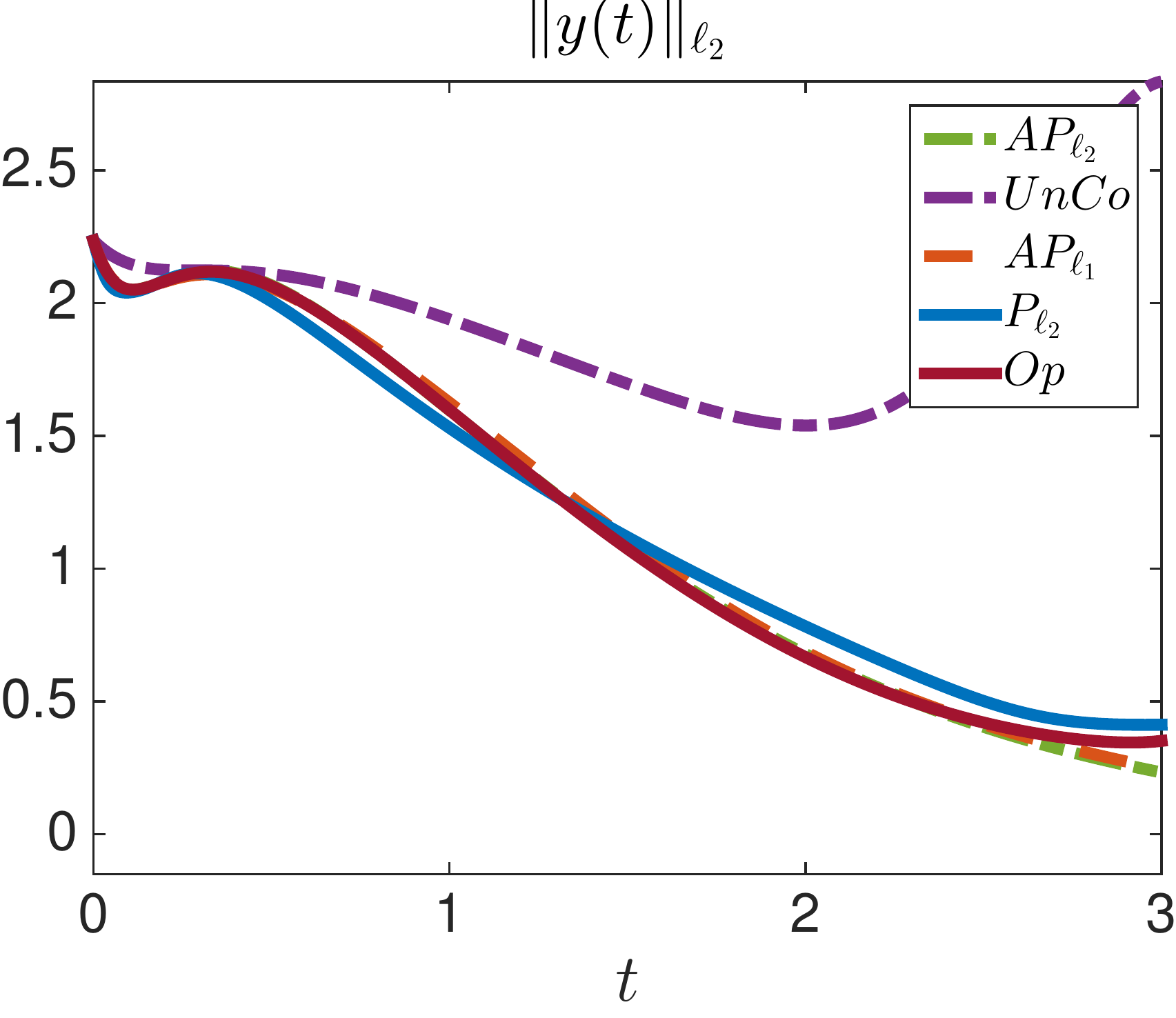}
	}
	\subfigure[Control]
	{
		\label{Fig13}
		\includegraphics[height=4cm,width=4cm]{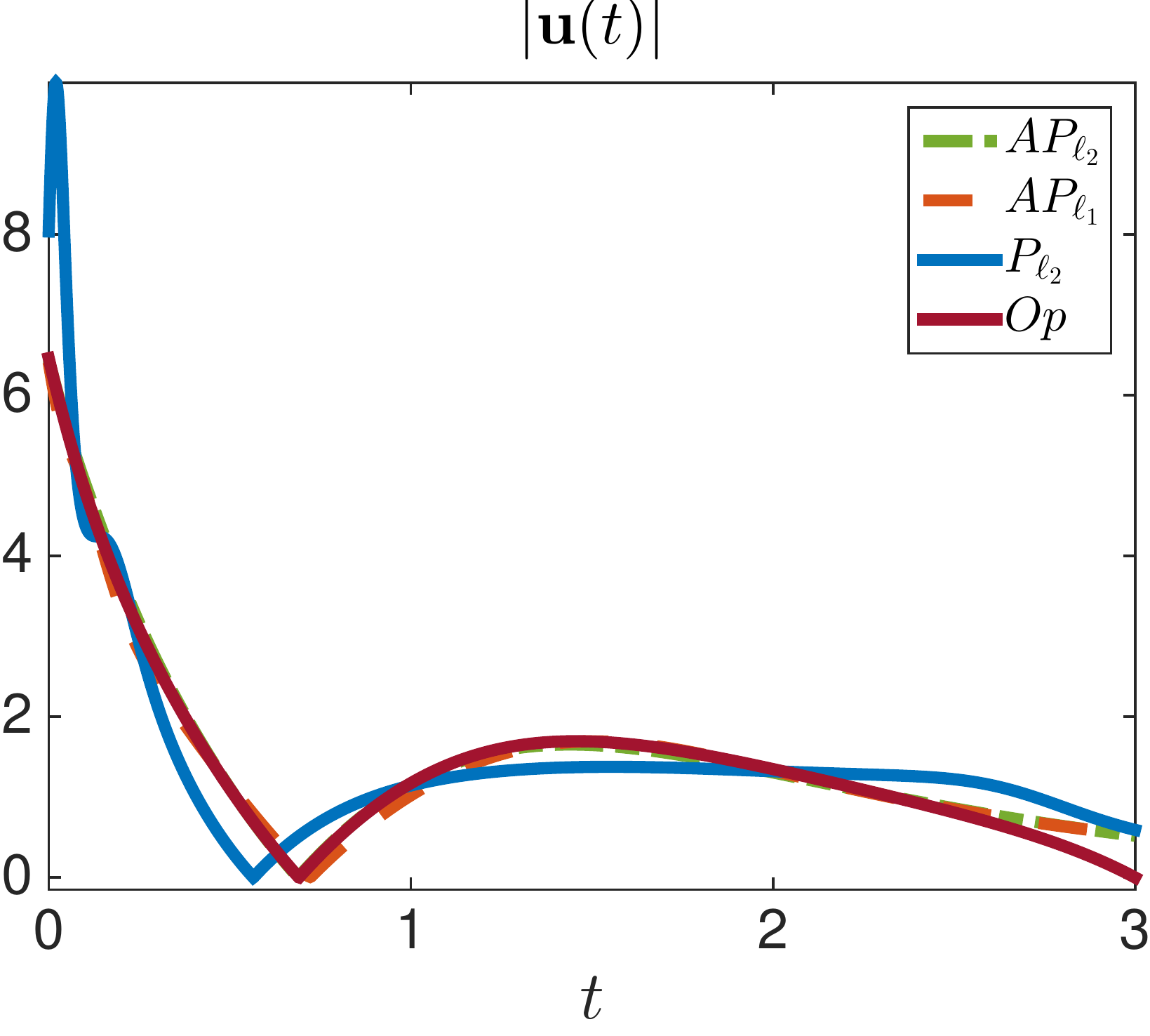}
	}
	\caption{Test 1. Evolution of  $\| \by(t) \|_{\ell_2}$ and  $|\mathbf{u}(t)|$ for $\bx = (2, 1)$. Here \textsl{UnCo} stands for the uncontrolled trajectory, and \textsl{Op} refers to the exact optimal trajectory. {\color{black} Results are qualitatively similar to Figure \ref{Fig12}}.}
\end{figure}

Moreover, we approximated the optimal control by using \eqref{e11} and \eqref{eq:feedback}. To be more precise, we compute the following feedback law:
\begin{equation}\label{eqkk3}
	\mathbf{u}_{\theta}(\bx)= - \frac{1}{2\beta}\sum_{ \mathbf{i} \in  \mathfrak{I}} \theta_{ \mathbf{i}}\nabla_x \Phi_{\mathbf{i}}(\bx).
\end{equation}
where $\nabla_x \Phi_{\mathbf{i}}$ stands for the gradient of the polynomial basis. We applied this feedback law for the choices  $\theta_{\ell_2}$, $\bar{\theta}_{\ell_2}$, and $\bar{\theta}_{\ell_1}$ for two initial vectors $(2,1)$ and $(2,-1)$. The evolution of the norm for the states controlled by these feedback laws,  compared to the optimal state, and the uncontrolled state is illustrated in Figures \ref{Fig10} and \ref{Fig11}. Figures  \ref{Fig12} and \ref{Fig13} depict the  evolution of the absolute value of the controls. Clearly,  the controls $\mathbf{u}_{ \bar{\theta}_{\ell_2}}$ and  $\mathbf{u}_{ \bar{\theta}_{\ell_1}}$ on the basis of
$\bar{V}_{\ell_2}$ and $\bar{V}_{\ell_1}$ approximate well the challenging behaviour of the optimal control, and they  outperform  $\mathbf{u}_{ {\theta}_{\ell_2}}$ obtained by  $V_{\ell_2}$.

\noindent Another advantage of using a sparse regression is the synthesis of a feedback law of reduced complexity. This is particularly relevant for the implementation of feedback laws in a real-time environment, where the number of calculations in the control loop needs to be minimized. In general, a feedback law expressed in the form
\begin{equation} 
	\bu_{\theta}(\bx)= - \frac{1}{2\beta}\bg^{\top}(\bx)\sum_{ \mathbf{i} \in  \mathfrak{I}} \theta_{ \mathbf{i}}\nabla_x \Phi_{\mathbf{i}}(\bx)\,,\qquad\text{with }\; g(\bx)\in\R^{n\times m}\,,\nabla\Phi_i(\bx)\in\R^n\,,
\end{equation}
requires $\mathcal{O}((mn^2+n)q)$ floating-point operations, where $q$ is the number of non-zero components in the expansion. Thus, the operation count decreases linearly with the level of sparsity. Going back to Table \ref{table1}, this implies a reduction of 63\% in the number of operations with respect to an $\ell_2$-based controller.


\subsection{Test 2: Controlled Allen-Cahn equation}\label{exp2}
In the following test, we consider  the PDE-constrained optimal control problem
\begin{align}
	\label{e13}
	&\min_{\mathbf{u}\in L^2(0,T;\mathbb{R}^2)} \int^{T}_{0}(\| y(t)\|^2_{L^2(0,1)}+\beta \|\mathbf{u}(t)\|^2_2)dt
\end{align}subject to
\begin{equation}\label{pdecontrol}
	\begin{cases}
		\partial_t y-\nu \partial^2_{x}y -y(1-y^2)   =  \sum^3_{i =1} u_i(t)\mathbf{1}_{\omega_i}    &\text{ in } (0,T)\times \Omega,\\
		\partial_x y (t,1) = \partial_xy(t,0) =0 &\text{ in } (0,T),  \\
		y(0,y)=y_0 &\text{ in } \Omega,
	\end{cases}
\end{equation}
with the 3-d control vector $\mathbf{u}(t):=[ u_1(t),u_2(t), u_3(t)]  \in L^2(0,4;\mathbb{R}^3)$, $\Omega  =(-1,1)$, $\nu = 0.1$, $\beta = 0.01$ and $T=4$. The control signals act through  $\mathbf{1}_{\omega_i}=\mathbf{1}_{\omega_i}(x)$, which denote the indicator functions with supports $\omega_1  = (-0.7,-0.4)$, $\omega_2  = (-0.2,0.2)$, and  $\omega_3  = (0.4,0.7)$. Due to the infinite-dimensional nature of the state equation this problem does not fall directly into the  optimal control setting considered in this paper. We first perform an approximation of \eqref{pdecontrol} in space by doing a pseudospectral collocation using Chebyshev spectral elements  with 18 degrees of freedom as in ~\citep{poli2}. This approximates the PDE control dynamics as an 18-dimensional nonlinear dynamical system. The resulting ODE system was treated numerically  by the Crank-Nicolson time stepping method with step-size $\Delta t = 0.005$.
Subsequently, a  dataset  $\{\mathbf{x}^j,V^j, V^j_x \}^{N}_{j =1}$ with $N=9000$ (including samples for training and validation) was generated by solving open-loop problem \eqref{e13}-\eqref{pdecontrol} for different values of quasi-randomly chosen initial vectors from the hypercube $[-10,10 ]^{18}$.

\begin{figure}[ht!]
	\centering
	\subfigure[$s =4$, $\lambda = 0.04$]
	{
		\includegraphics[height=4cm,width=4cm]{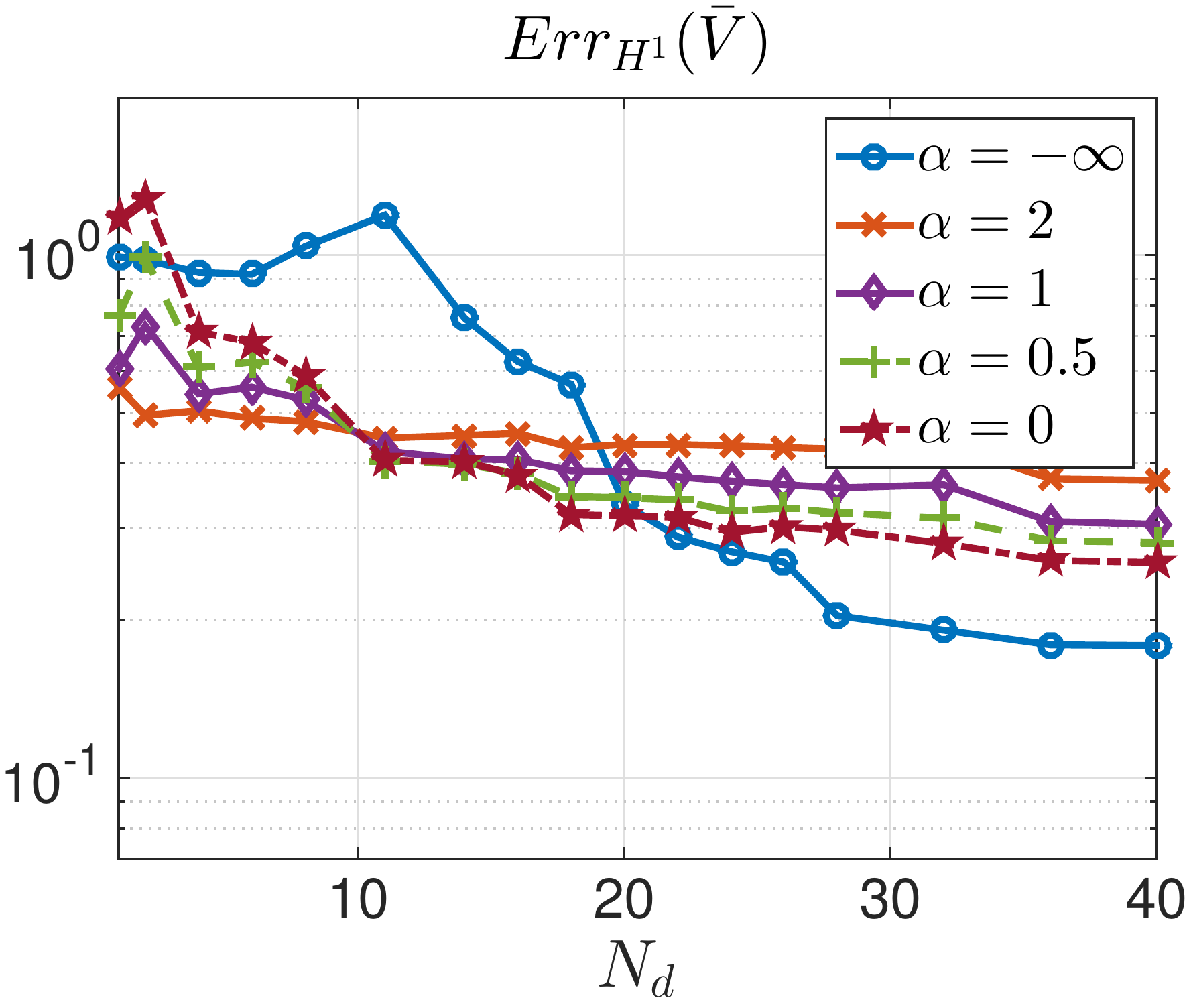}
	}
	\subfigure[ $s =4$, $\lambda = 0.04$]
	{
		\includegraphics[height=4cm,width=4cm]{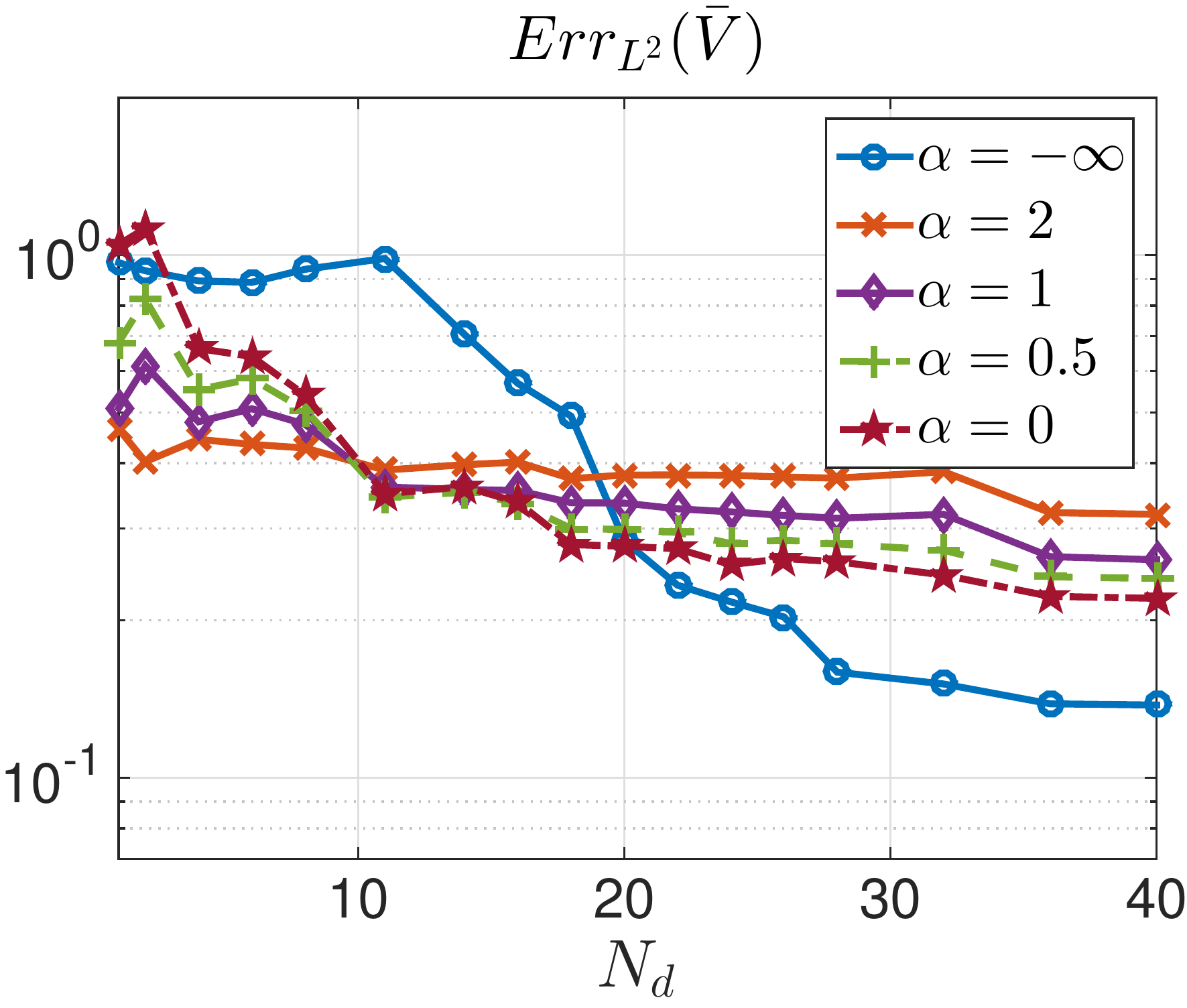}
	}
	\subfigure[$s =4$, $\lambda = 0.04$]
	{
		\includegraphics[height=4cm,width=4cm]{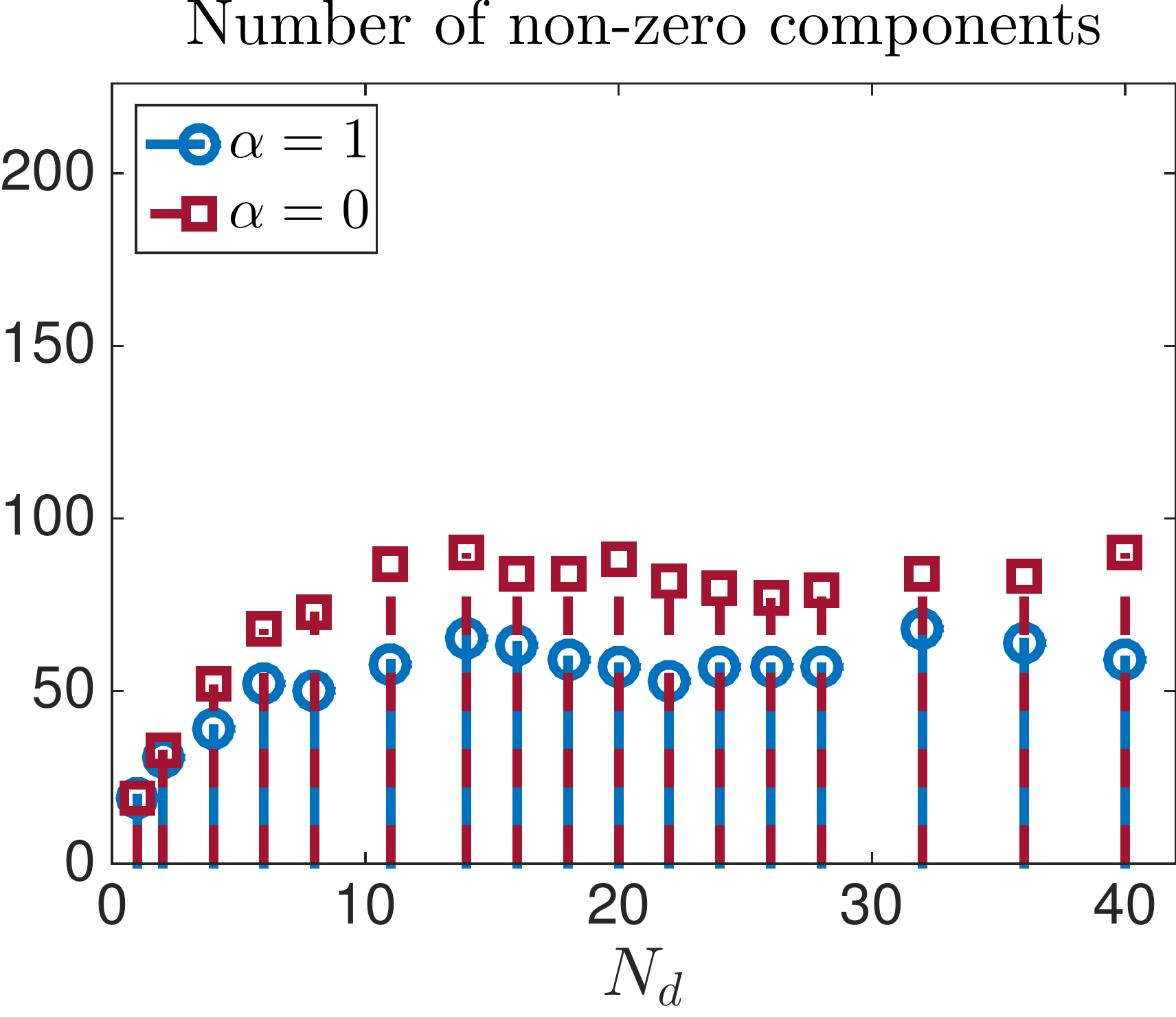}
	}
	\subfigure[ $s =4$, $\lambda = 0.008$]
	{
		\includegraphics[height=4cm,width=4cm]{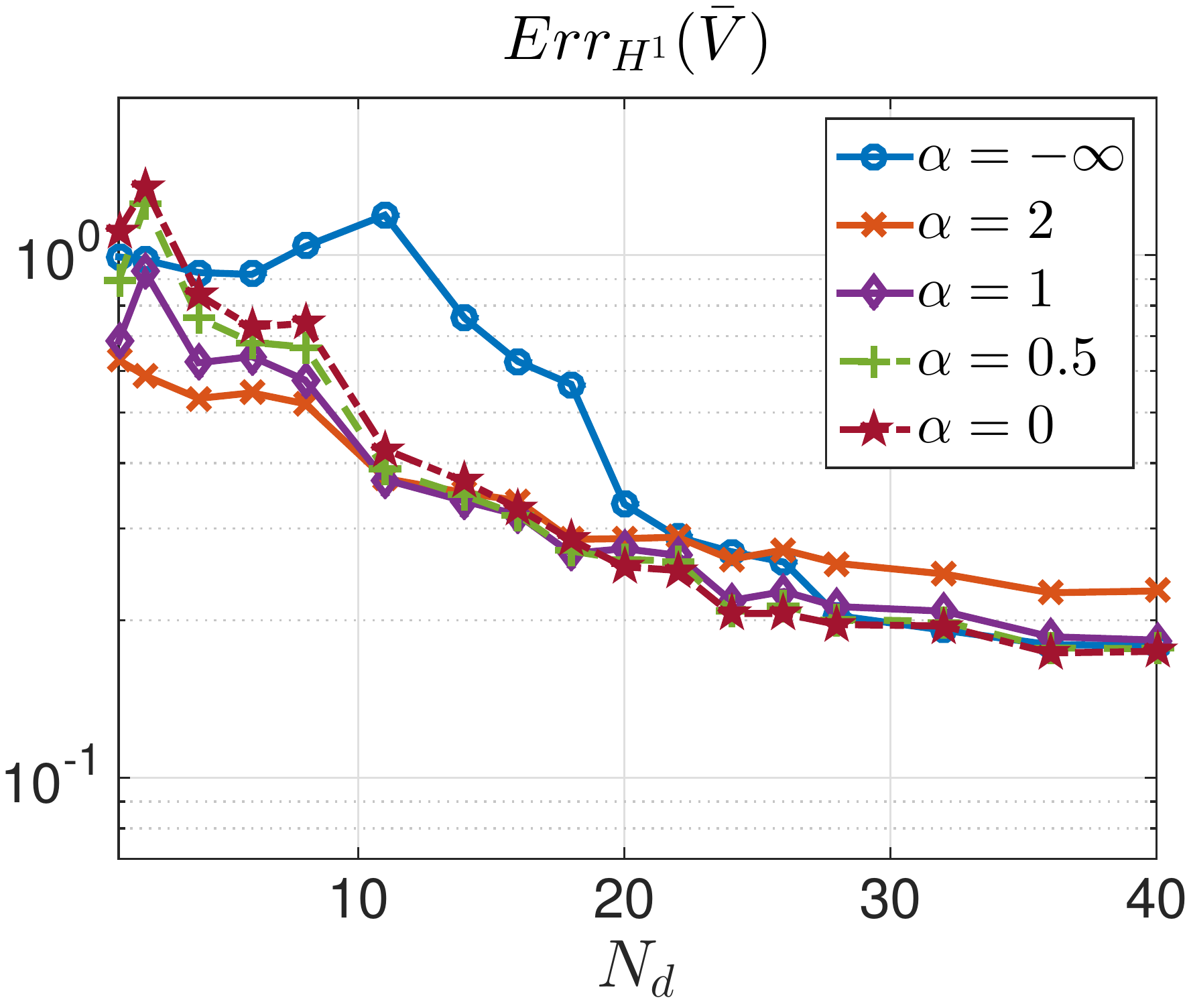}
	}
	\subfigure[$s =4$, $\lambda = 0.008$ ]
	{
		\includegraphics[height=4cm,width=4cm]{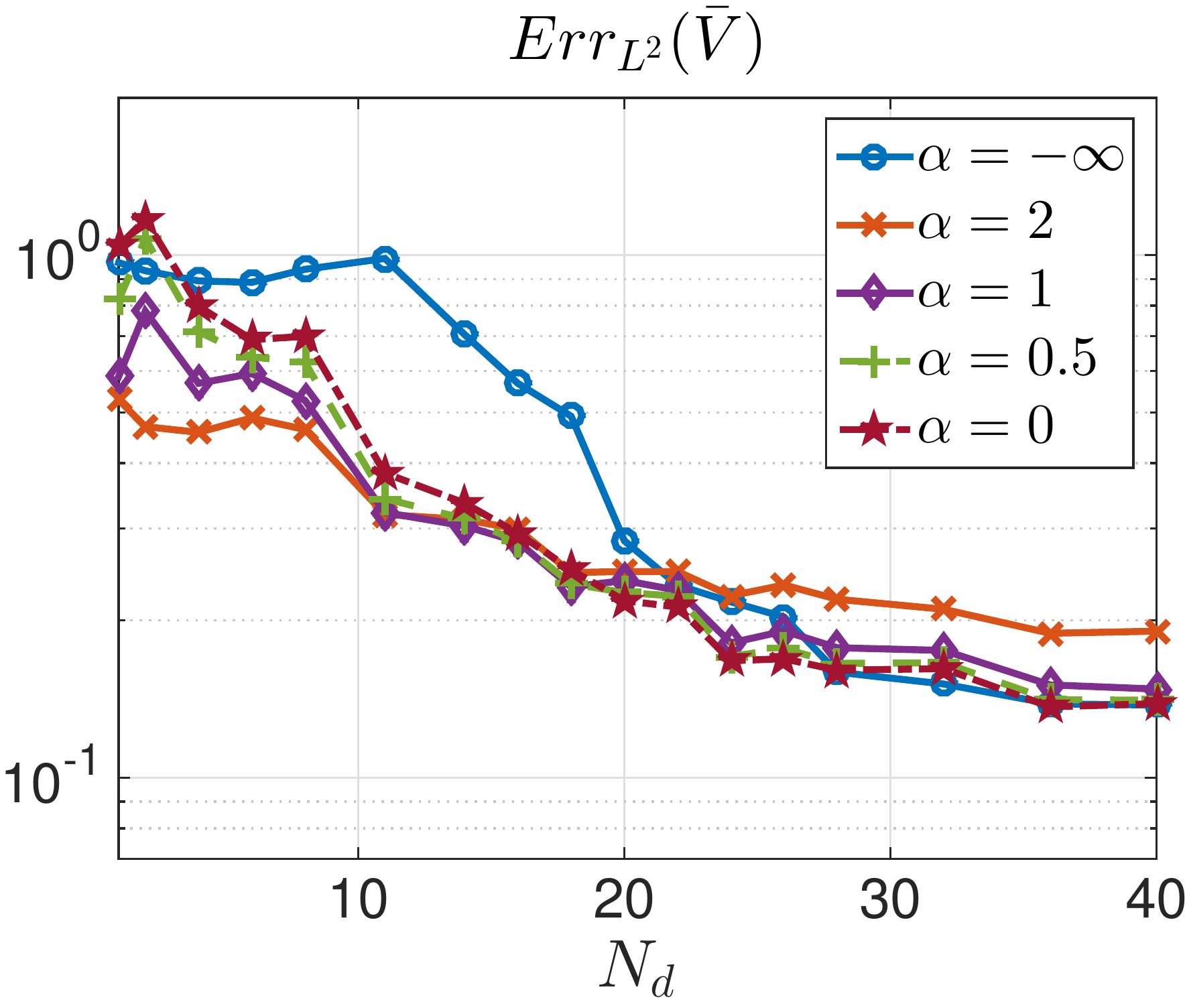}
	}
	\subfigure[$s =4$, $\lambda = 0.008$]
	{
		\includegraphics[height=4cm,width=4cm]{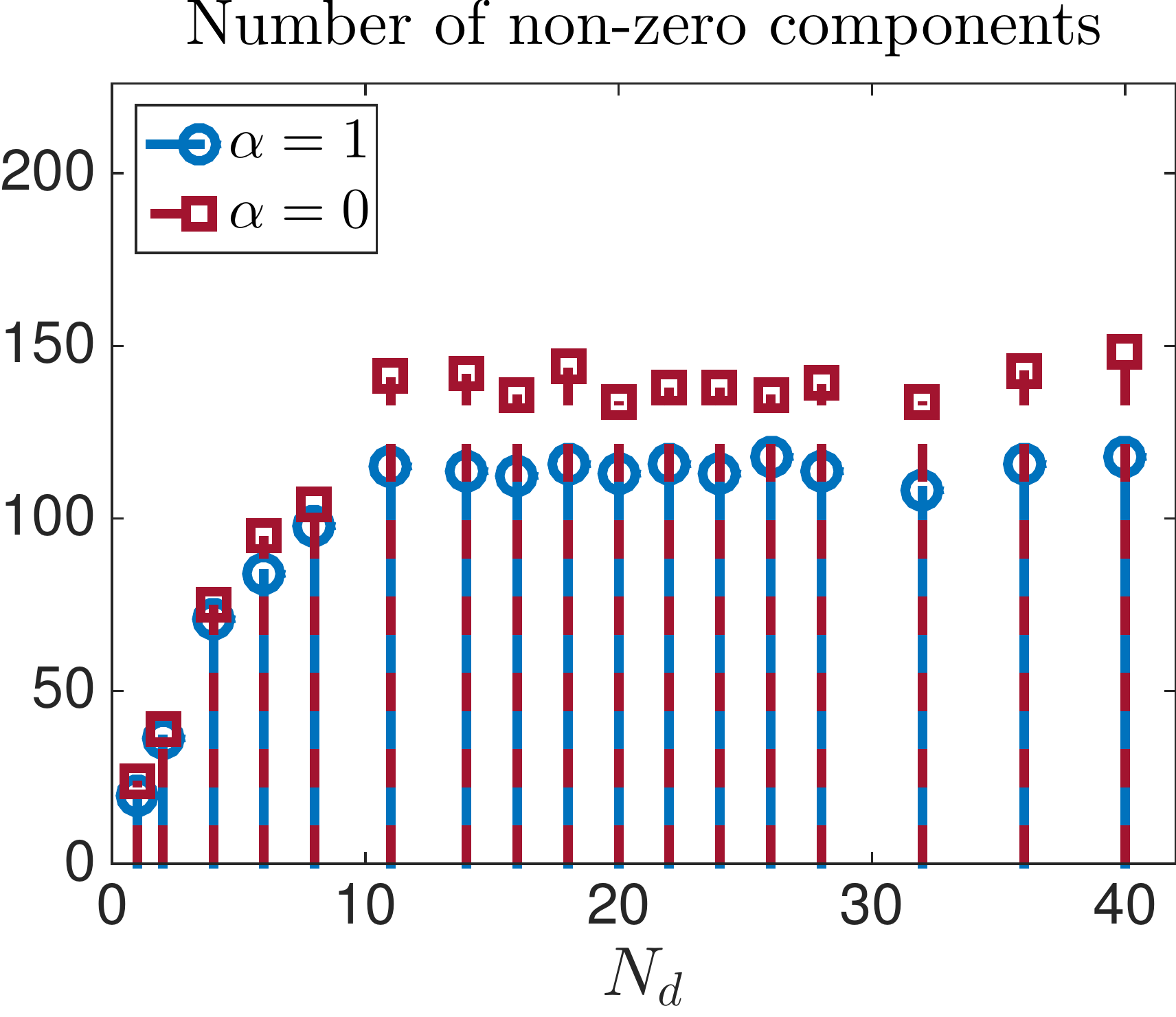}
	}
	\subfigure[$s =8$, $\lambda =0.008$ ]
	{
		\includegraphics[height=4cm,width=4cm]{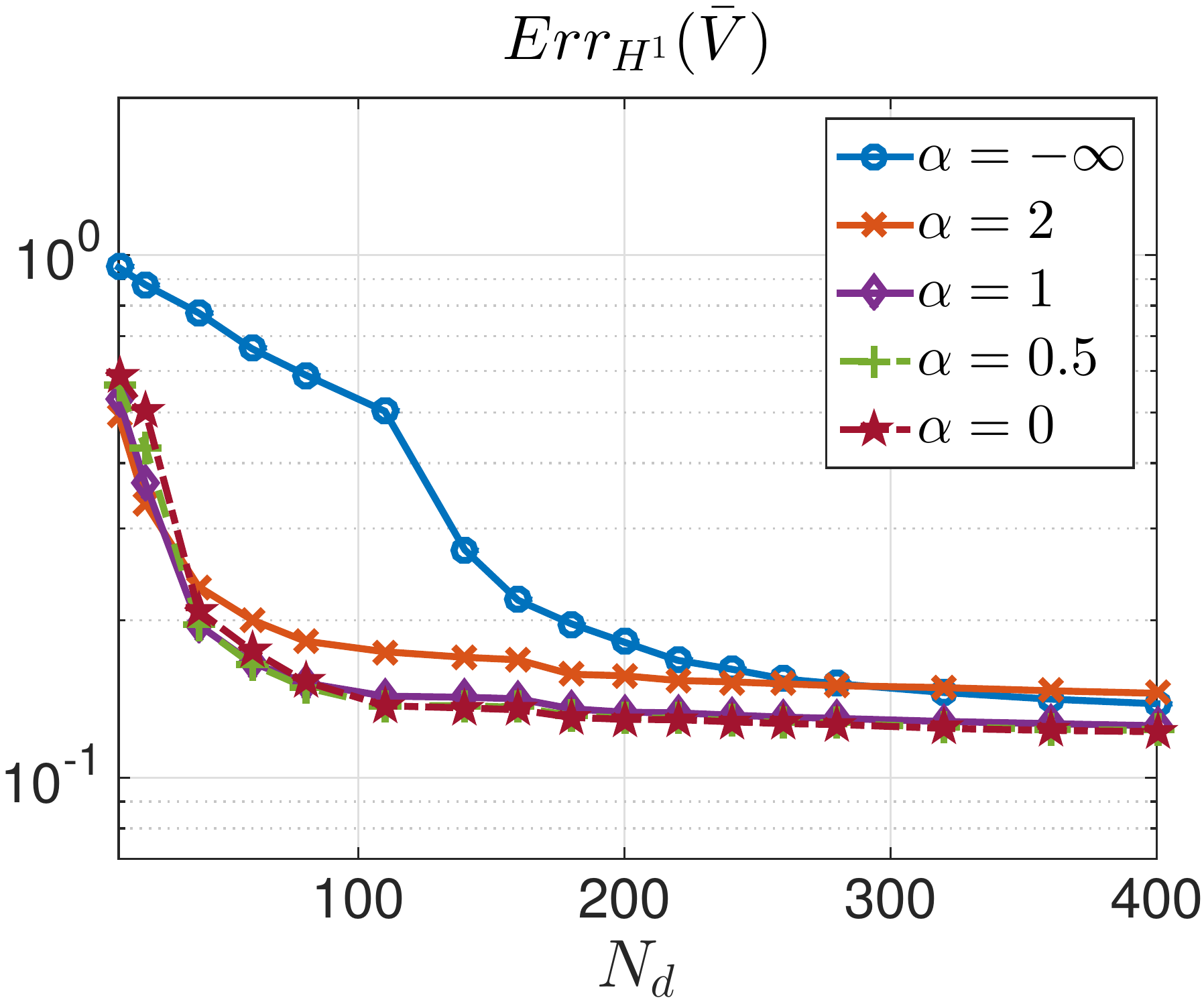}
	}
	\subfigure[ $s =8$, $\lambda =0.008$]
	{
		\includegraphics[height=4cm,width=4cm]{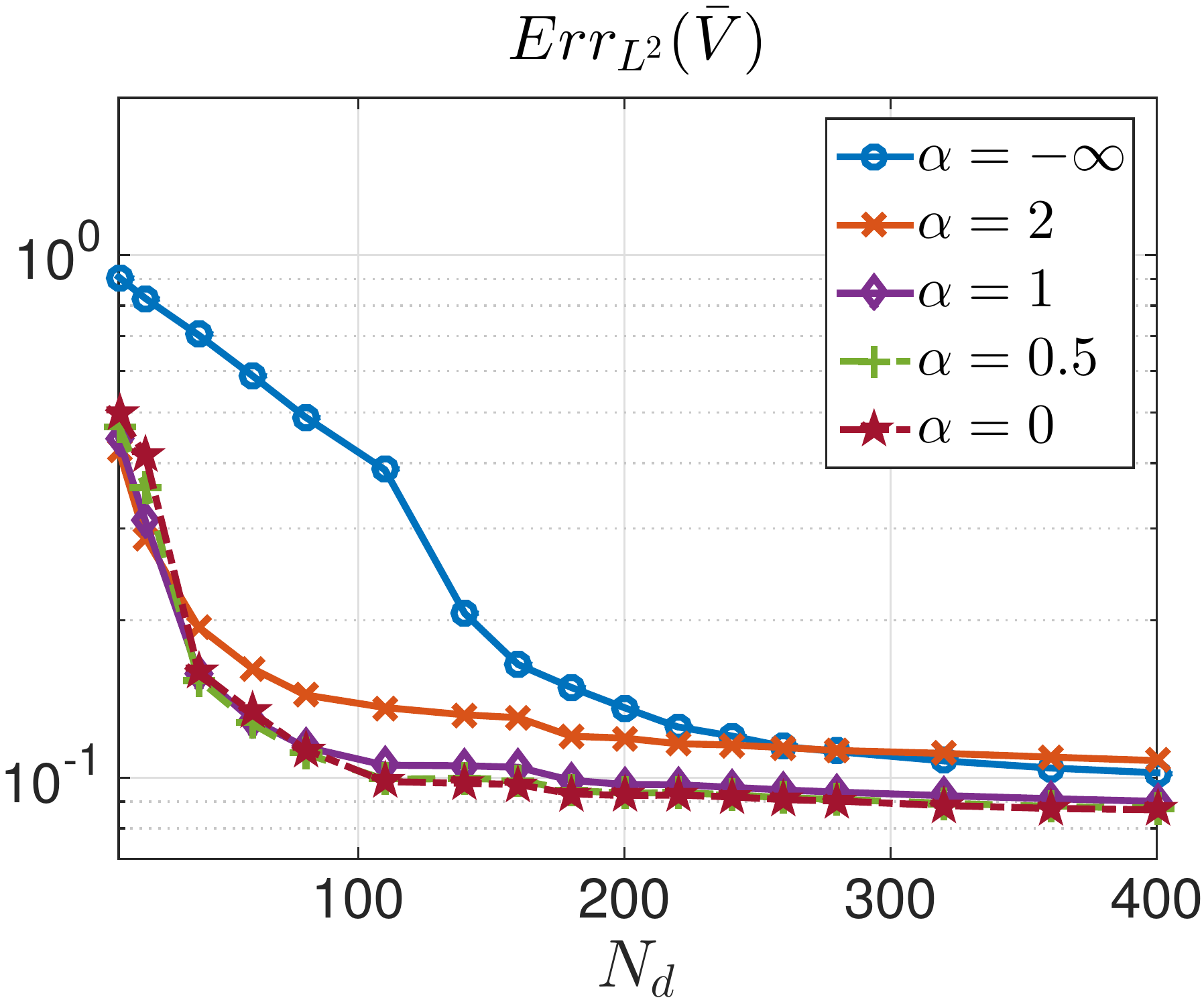}
	}
	\subfigure[$s =8$, $\lambda = 0.008$]
	{
		\includegraphics[height=4cm,width=4cm]{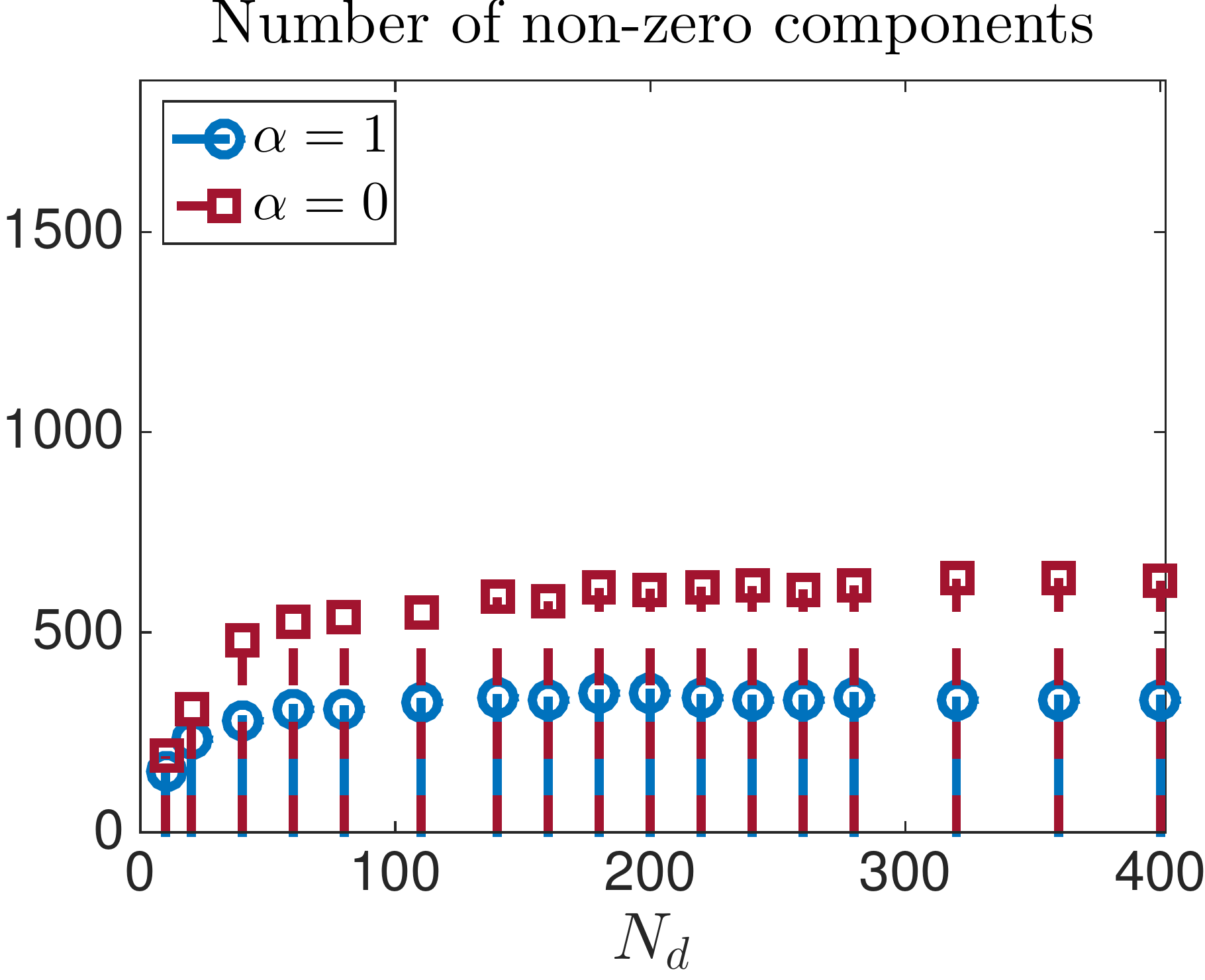}
	}
	\caption{Test 2. Numerical  results for the gradient-augmented polynomial approximation using the Chebyshev  polynomial basis. {\color{black} Due to the inclusion of gradient-augmented information in the training, the validation error in the $H^1$ and $L^2$ norms decrease as the number of samples is moderately increased. The sparsity pattern of the resulting feedback law can be controlled through the parameter $\alpha$ which determines the weight in the $\ell_1$-norm penalty term.}}
	\label{Fig17}
\end{figure}

\begin{figure}[ht!]
	\centering
	
	\subfigure[$s =4$, $\lambda = 0.01$ ]
	{
		\includegraphics[height=4cm,width=4cm]{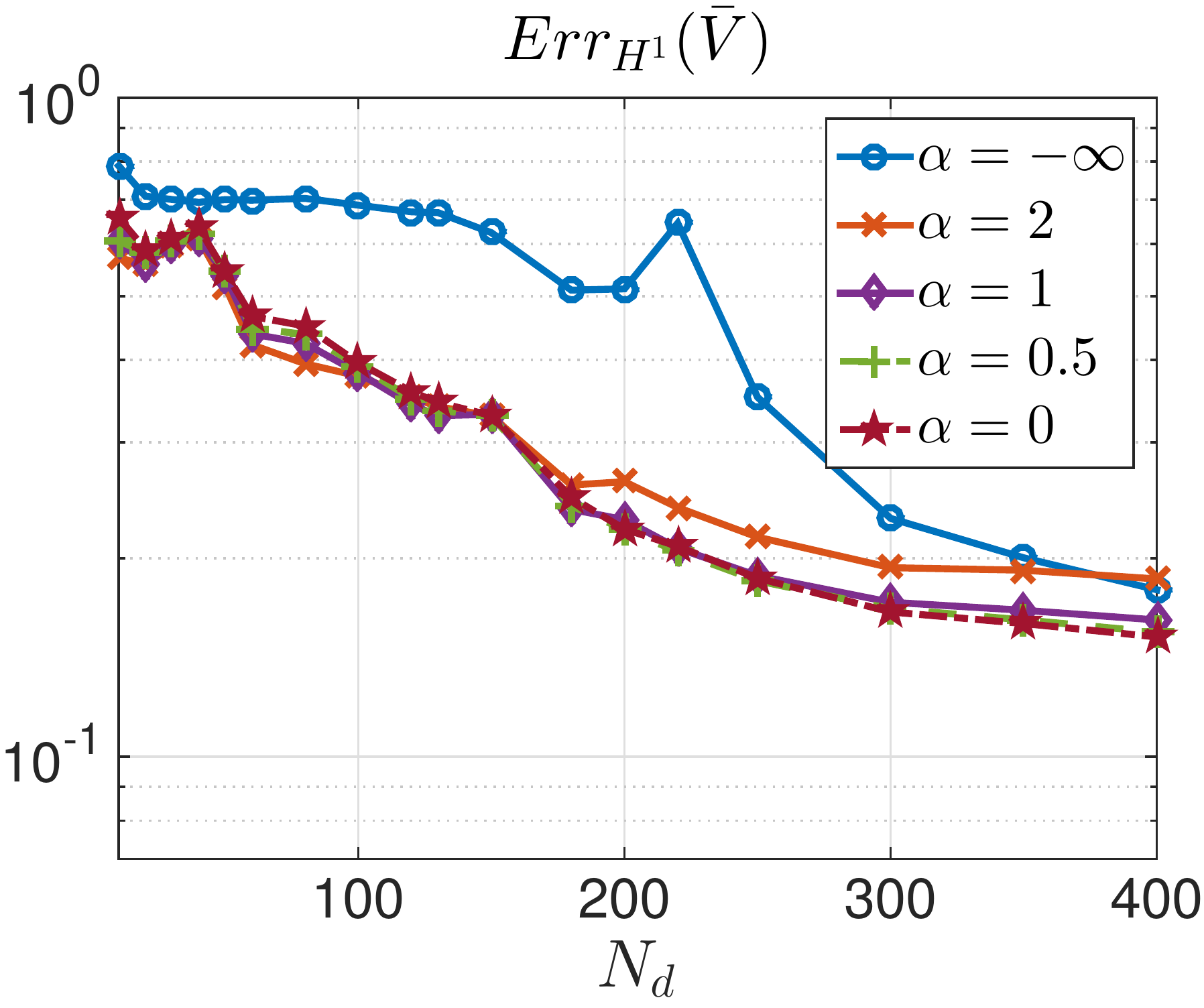}
	}
	\subfigure[$s =4$, $\lambda = 0.01$ ]
	{
		\includegraphics[height=4cm,width=4cm]{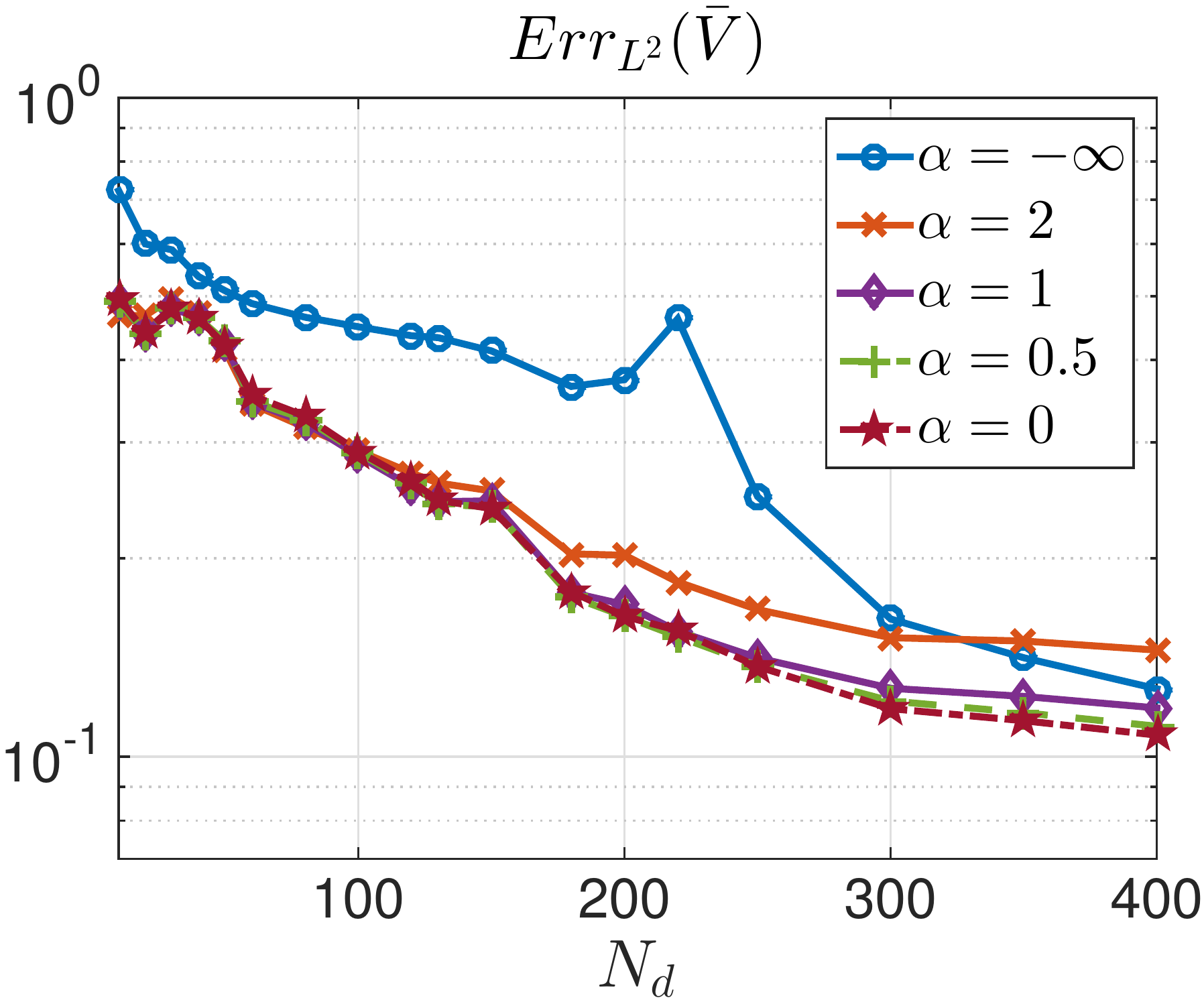}
	}
	\subfigure[$s =4$, $\lambda = 0.01$]
	{
		\includegraphics[height=4cm,width=4cm]{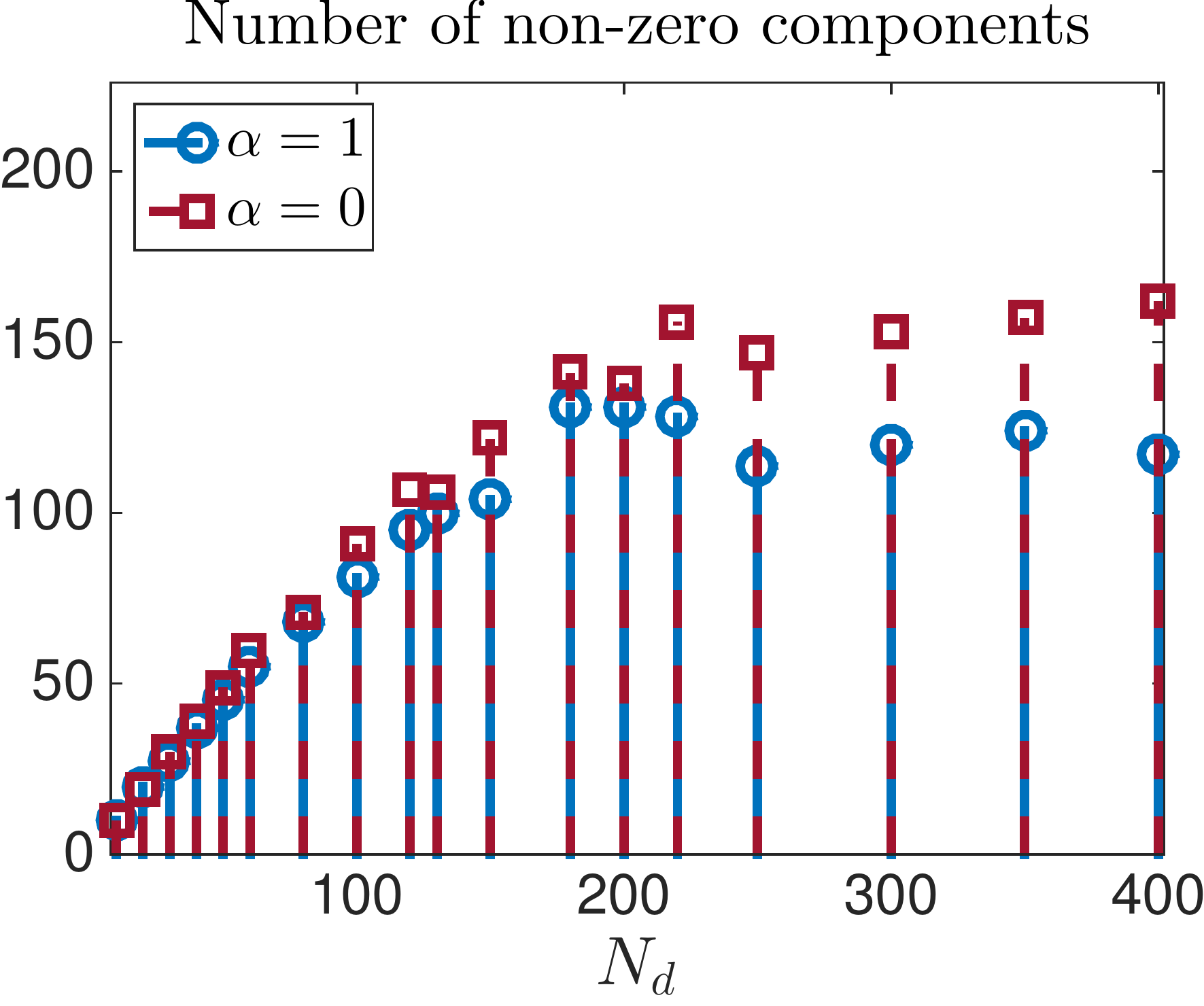}
	}
	\subfigure[$s =8$, $\lambda = 0.01$ ]
	{
		\includegraphics[height=4cm,width=4cm]{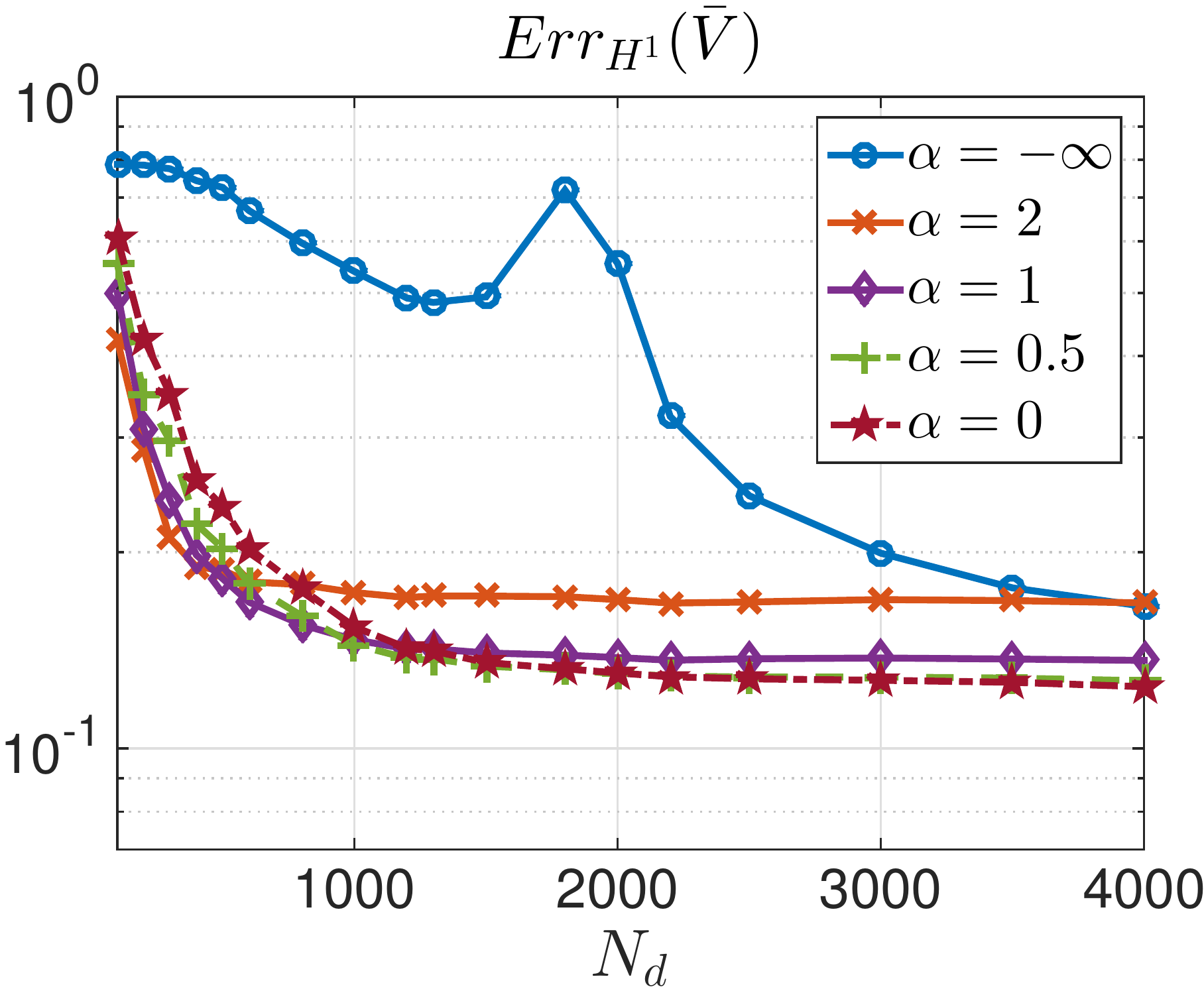}
	}
	\subfigure[$s =8$, $\lambda =0.01$ ]
	{
		\includegraphics[height=4cm,width=4cm]{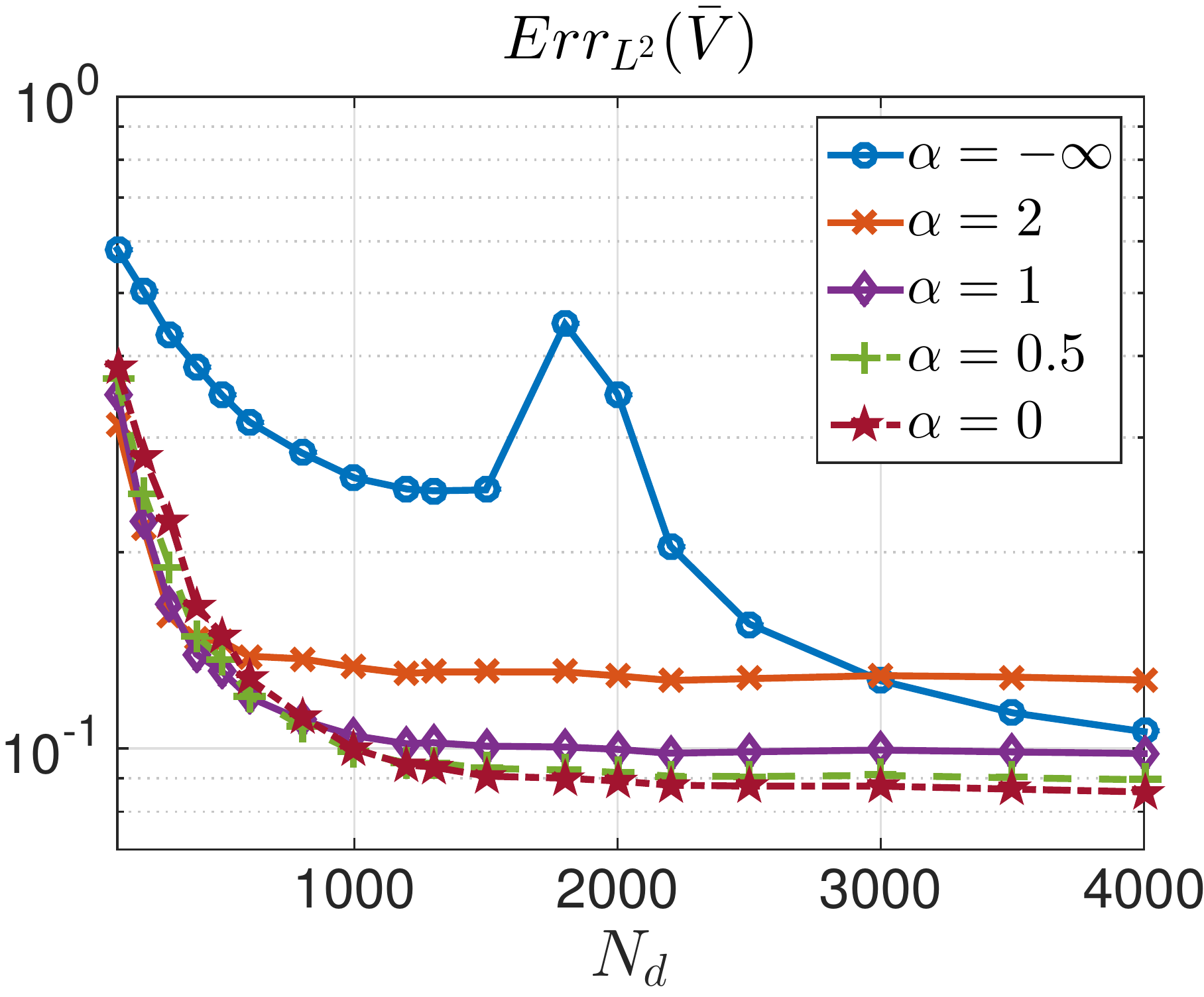}
	}
	\subfigure[$s =8$, $\lambda =0.01$]
	{
		\includegraphics[height=4cm,width=4cm]{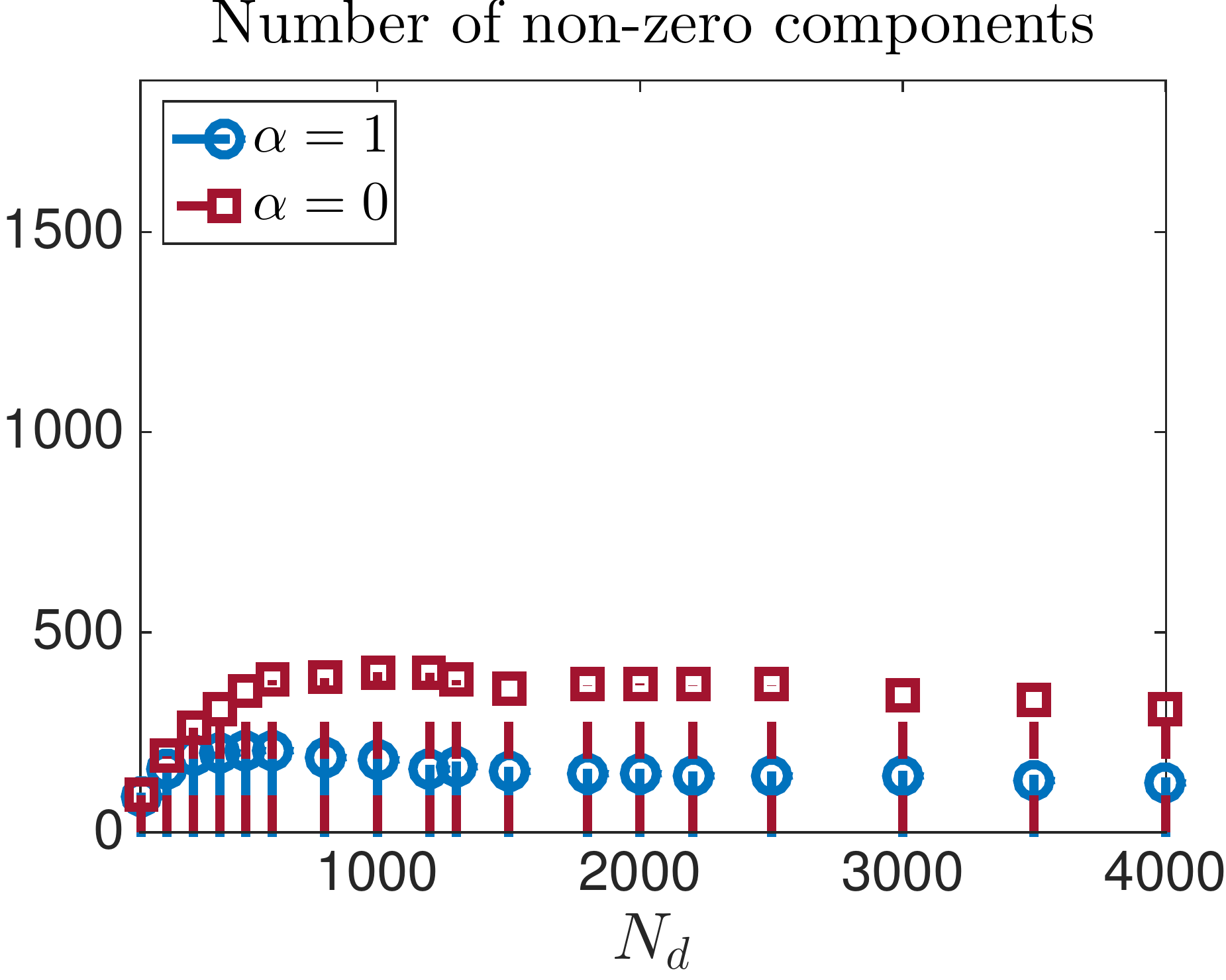}
	}
	\caption{Test 2. Numerical  results for the polynomial approximation without gradient informations using the Chebyshev  polynomial basis. {\color{black} Results are qualitatively similar to Figure \ref{Fig17} with the Legendre polynomial basis.}}
	\label{Fig18}
\end{figure}

For this example we used the two different values  $s = 4$ and $s =8$  in the hyperbolic cross index set  $\mathfrak{I}(s)$. For these choices we have  $|\mathfrak{I}(4)|  =226$ and $|\mathfrak{I}(8)|  =1879$, resulting in  $226$ (resp. $1897$) polynomial  basis functions to approximate the value functions $V: = V(0,\cdot)$.  We computed the solutions $\theta_{\ell_2}$, $\bar{\theta}_{\ell_2}$, $\theta_{\ell_1}$ and $\bar{\theta}_{\ell_1}$ to  problems \eqref{e4} , \eqref{e6}, \eqref{e5}, and  \eqref{e7} for the different choices of $N_d$,  $\alpha$, and polynomial basis.  For problems \eqref{e5} and  \eqref{e7} we show results with $\lambda = 0.01$ and $\lambda = 0.008, 0.04,$ respectively, using Chebyshev polynomials. Similarly as in Test 1, we report   the level of non-sparsity of  $\theta_{\ell_1}$, $\bar{\theta}_{\ell_1}$. The errors \eqref{eq:kk1} are shown for a validation set with $|\mathcal{I}_{val}| = 5000$ samples.
These results are depicted in  Figures \ref{Fig17} and \ref{Fig18}.

Overall, these results allow to draw the same conclusions as for the previous test. In particular, as
${N}_d$ increases, the validation errors are getting smaller. Again, comparing Figures \ref{Fig17} and \ref{Fig18} the  gradient-augmented results reach the lowest errors with significantly smaller datasets $N_d$ than the gradient-free ones.
Figures \ref{Fig17}-\ref{Fig18} also confirm that the errors for  $s = 8$ are smaller  compared to  $s = 4$.  Naturally, for polynomials basis associated to $s = 8$, we need to increase the training data in comparison to $s = 4$.  Comparing rows 1 and 2 in Figure  \ref{Fig17}, we observe that decreasing $\lambda$ results in an increase on the number of non-zero components and in a decrease of the errors.

\begin{figure}[ht!]
	\centering
	\subfigure
	{
		\label{Fig19}
		\includegraphics[height=5cm,width=5cm]{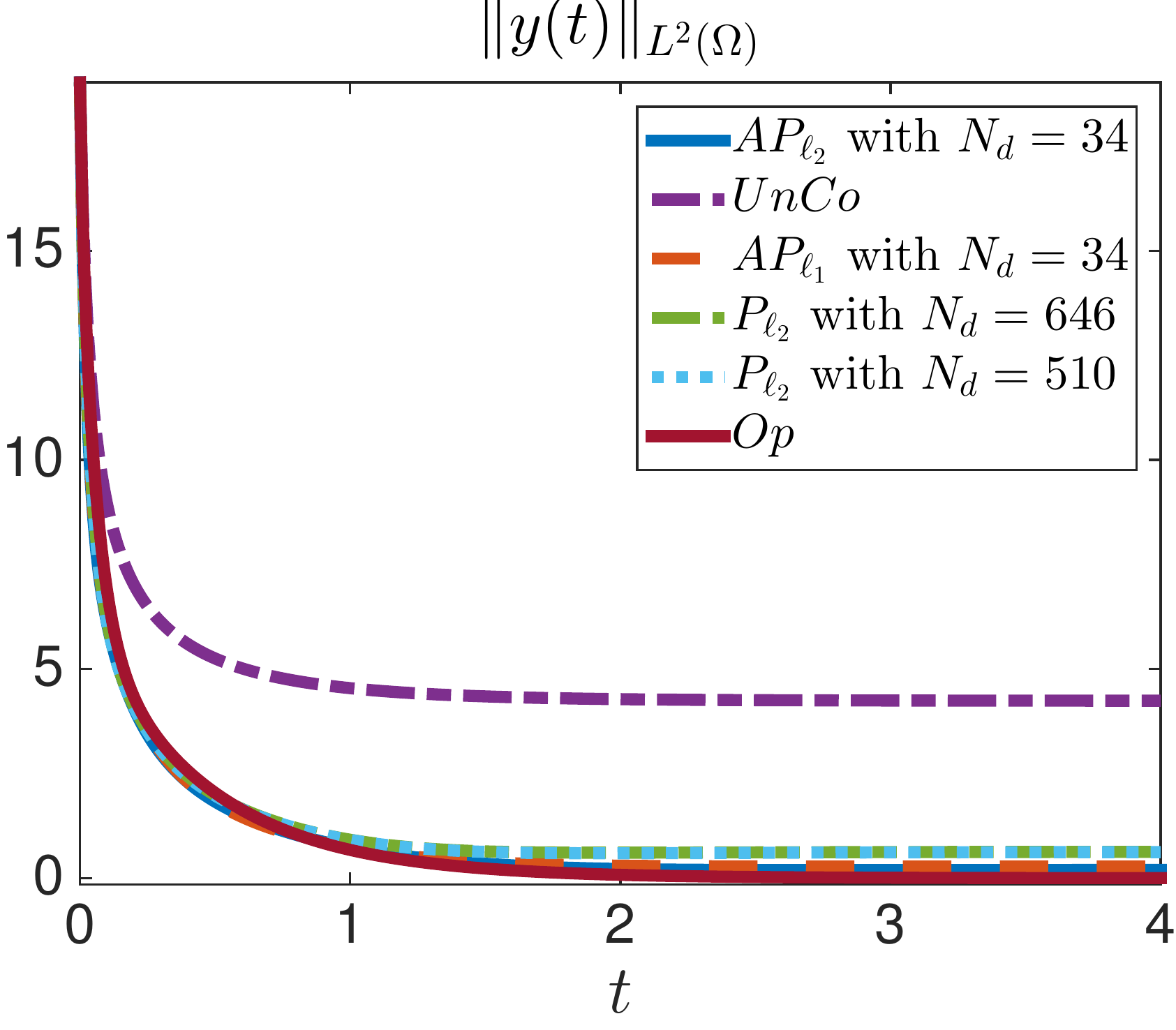}
	}	
	\subfigure
	{
		\label{Fig20}
		\includegraphics[height=5cm,width=5cm]{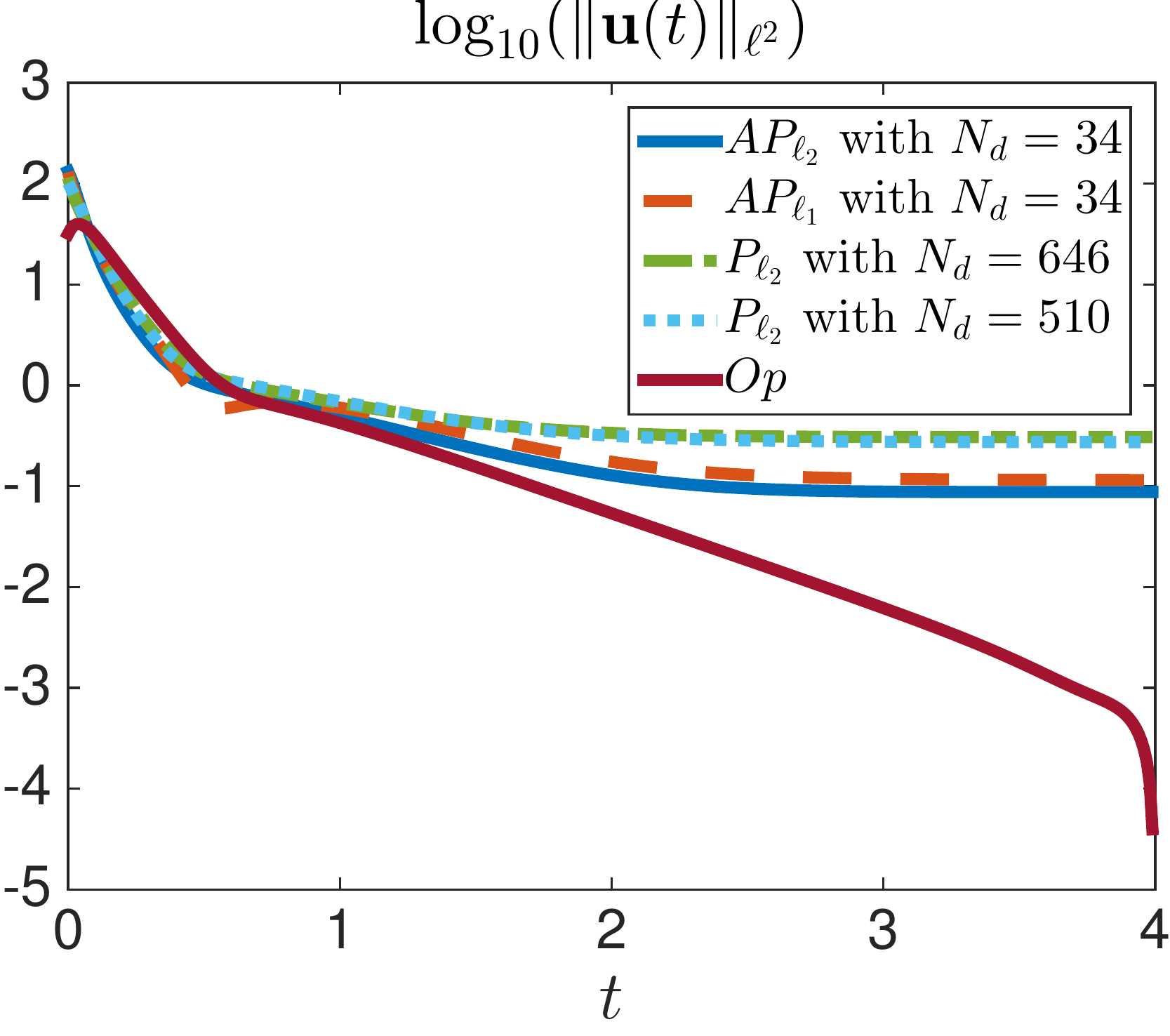}
	}
	\caption{Test 2. Evolution of  $\| y(t) \|_{L^2(\Omega)}$ and  $\log_{10}(\|\mathbf{u}(t)\|)$  for the choices $\lambda = 0.008$, $\alpha = 1$,  and Chebyshev polynomial basis.{\color{black} The gradient-augmented feedback laws outperform, both in stabilization and control energy, the control laws recovered without gradient information.}}
\end{figure}

We approximated the optimal control using \eqref{eq:feedback} and \eqref{e11}, resulting in the feedback law
\begin{equation}
	\label{eq:feedback2}
	\mathbf{u}_{\theta}(\bx)= - \frac{1}{2\beta}\bg^{\top}\sum_{ \mathbf{i} \in  \mathfrak{I}} \theta_{ \mathbf{i}}  \nabla_x \Phi_{\mathbf{i}}(\bx)\,,
\end{equation}
where $ \bg: = [\mathbf{1}_{\omega_1}| \mathbf{1}_{\omega_2}|\mathbf{1}_{\omega_3}]\in\R^{n\times 3}$. We applied this feedback law for the choices  $\theta_{\ell_2}$, $\bar{\theta}_{\ell_2}$, and $\bar{\theta}_{\ell_1}$ on the initial condition
\begin{equation}
	y_0(x) = (x-1)(x+1)+5.
\end{equation}
We report the results for the Chebyshev polynomial basis with $s = 4$ and thus $q = 226$. For the case $\bar{\theta}_{\ell_1}$ we set $\lambda = 0.008$ and  $\alpha = 1$.  The evolution of the norm of the resulting controlled, optimal, and uncontrolled  states are depicted  in Figures \ref{Fig19}, and the associated controls in \ref{Fig20}.

\begin{figure}[ht!]
	\centering
	\subfigure[Uncontrolled ]
	{
		\label{Fig23}
		\includegraphics[height=5cm,width=5cm]{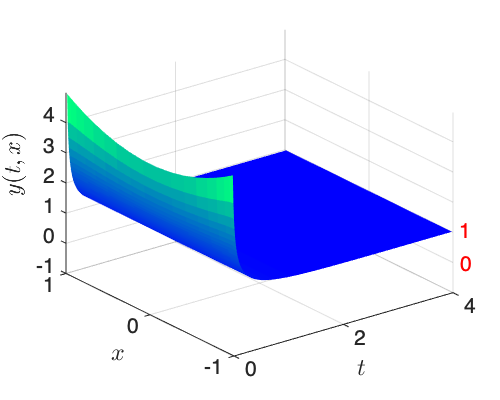}
	}
	\subfigure[Controlled]
	{
		\label{Fig24}
		\includegraphics[height=5cm,width=5cm]{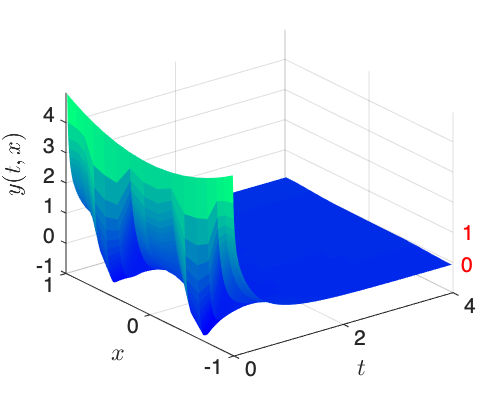}
	}
	\caption{Test 2. The uncontrolled state and  the controlled state by $\mathbf{u}_{\bar{\theta}_{\ell1}}$ for $N_d = 34$, $\lambda = 0.008$, $\alpha = 1$, and $\alpha = 1$. {\color{black}The recovered feedback law effectively stabilizes the state of the system around the equilibrium $y\equiv 0$.}}
\end{figure}

Figure \ref{Fig23} depicts the uncontrolled state. It converges to the stable equilibrium given by the constant function with value 1.  The state controlled  by $\mathbf{u}_{\bar{\theta}_{\ell1}}$ is illustrated  in Figure \ref{Fig24}. Here the state tends to 0 as expected.

\subsection{Test 3: Optimal consensus control in the Cucker-Smale model}\label{exp3}
We conclude with a thorough discussion of a high-dimensional, non-linear, non-local optimal control problem related to consensus control of agent-based dynamics ~\citep{bailo,bongini,caponigro}. We study the  Cucker-Smale model ~\citep{MR2324245} for consensus control with $N_a$ agents with states $(y_i,v_i)\in\R^d\times\R^d$ for  $i = 1, \dots, N_a$, where  $y_i$ and $v_i$ stand for the position and velocity of the $i$-th agent,  respectively, and $d \in \mathbb{N}$ is the dimension of the physical space. Then dynamics of the agents are governed by
\begin{equation}
	\label{e17}
	\begin{split}
		\frac{d y_i}{dt}&=v_i\,,   \\
		\frac{d v_i}{dt}&= \frac{1}{N_a}\sum_{j=1}^{N_a} \frac{v_j-v_i}{1+\|y_i-y_j\|^2}+u_i\,,  \qquad   i = 1, \dots, N_a\,,\\
		y_i(0)&=x_i\,,\quad v_i(0) = w_i\,.
	\end{split}
\end{equation}
The consensus control problem consists of finding a control $\mathbf{u}(t):=(u_1(t),\dots,u_{N_{a}}(t))\in \R^{d\times N_a}$  which steers the system towards the consensus manifold
\begin{align}
	v_i=\bar v=\frac{1}{N_a}\sum\limits_{j=1}^{N_a}v_j\,,\quad\forall i=1,\ldots,N_a\,.
\end{align}
Asymptotic consensus emergence is conditional to the cohesiveness of the initial state $\bx_0=(x_1,\ldots,x_{N_a})$ and $\bw_0=(w_1,\ldots,w_{N_a})$ ~\citep{MR2324245}. To remove this dependence on the initial state, we cast this problem as an optimal control problem by defining the following cost functional
\begin{align}
	\label{eq:objCS}
	J(\mathbf{u};\mathbf{x}_0,\mathbf{w}_0):=\int_0^{\top}\sum_{i=1}^{N_a}\frac{1}{N_a}\| v_i(t)-\bar v \|^2+\beta\| u_i(t)\|^2  \,dt,
\end{align}
and formulating the optimal control problem
\begin{equation}
	\label{e18}
	\tag{$OC(\mathbf{\hat x}_0)$}
	\min_{\mathbf{u} \in L^2(0,T;\mathbb{R}^{d\times N_a})} \{J(\mathbf{u};\mathbf{x}_0,\mathbf{v}_0)\text{ subject to \eqref{e17}} \}.
\end{equation}

For the sake of completeness, in this case the adjoint system is given by ~\citep{bailo}
\begin{align}
	-\frac{dp_{y_i}}{dt} &=\frac{1}{N_a}\sum_{j\neq i}\frac{-2(p_{v_j}-p_{v_i})}{(1+\|y_j-y_i\|^2)^2}[(y_i-y_j)\otimes(v_i-v_j)]  \,,    \\
	-\frac{dp_{v_i}}{dt} &=p_{y_i}+\frac{1}{N_a}\sum_{j\neq i}\frac{p_{v_j}-p_{v_i}}{1+\|y_i-y_j\|^2}+\frac{2}{N}(v_i- \bar v)^{\top}  \qquad  i = 1, \dots, N_a\,, \\
	p_{y_i}(T)&=0 ,\quad p_{v_i}(T)=0\,,
\end{align}
and  the optimality condition reads
\begin{equation}
	p_{v_i}(t)+2\beta u_i^*(t)=0\,  \qquad \forall\, t\in (0,T)\,,\,  i=1,\ldots,N\,.
\end{equation}
We denote the augmented initial state $\hat{\bx}_0=(\bx_0,\bw_0)$, and  we approximate the value function
$V(\hat\bx)=V(0, \hat\bx)$. We set $N_a=20$, $d=2$, $T=10$, and $\beta=0.01$,  and  we compute a dataset  $\{\hat\bx^j,V^j, V^j_{\hat x} \}^{N}_{j =1}$  with $N=10^4$. For every $j$, the initial vectors $\hat\bx^j \in \mathbb{R}^{80}$ were chosen quasi-randomly from the hypercube  $[-3,3 ]^{80}$. The dataset was computed by  solving open-loop problems with a time discretization using the fourth order Runge-Kutta method with step-size $\Delta t = 0.01$.

\begin{figure}[ht!]
	\centering
	
	\subfigure[$\lambda = 0.0005$]
	{
		\includegraphics[height=4cm,width=4cm]{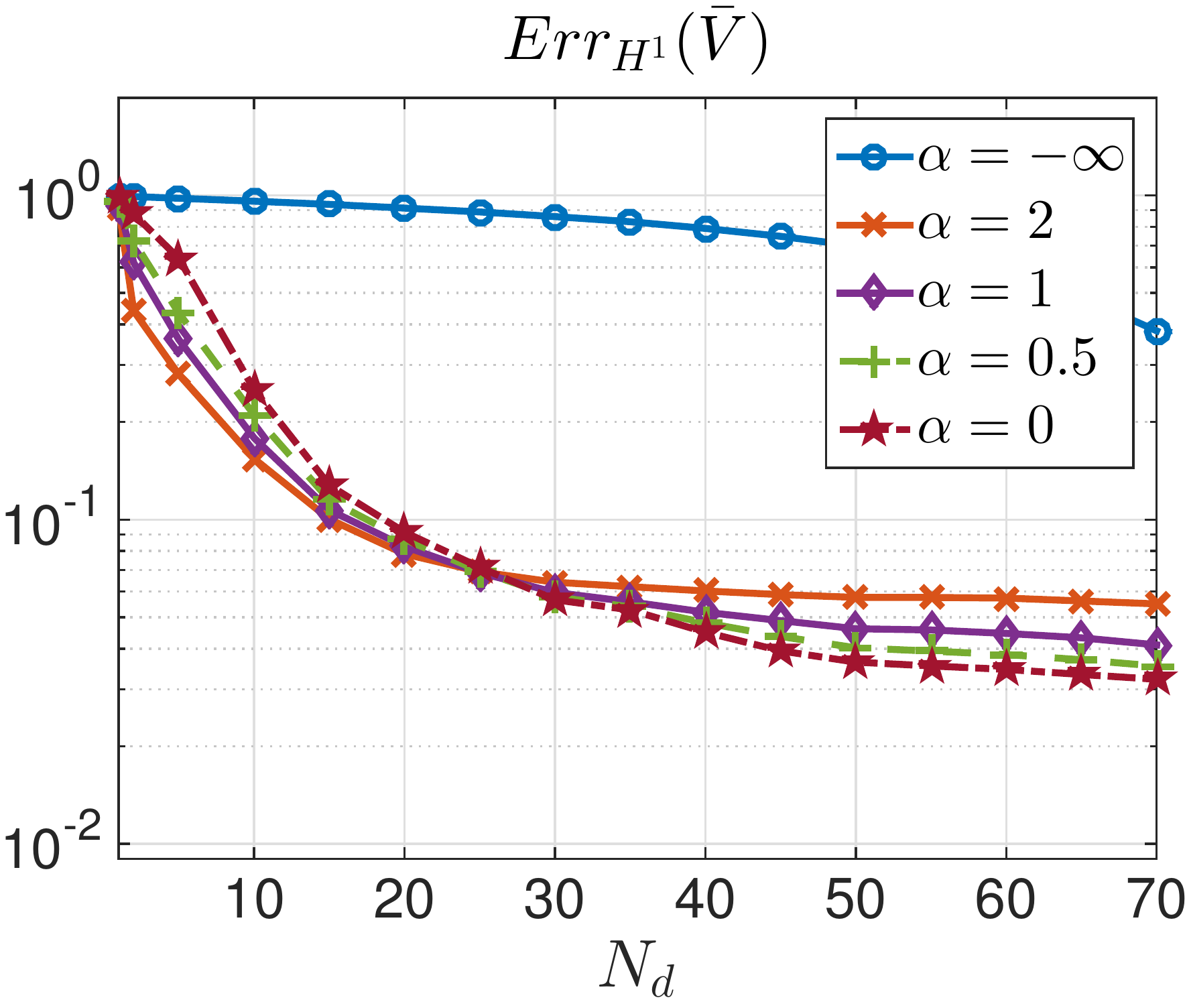}
	}
	\subfigure[$\lambda =0.0005$]
	{
		\includegraphics[height=4cm,width=4cm]{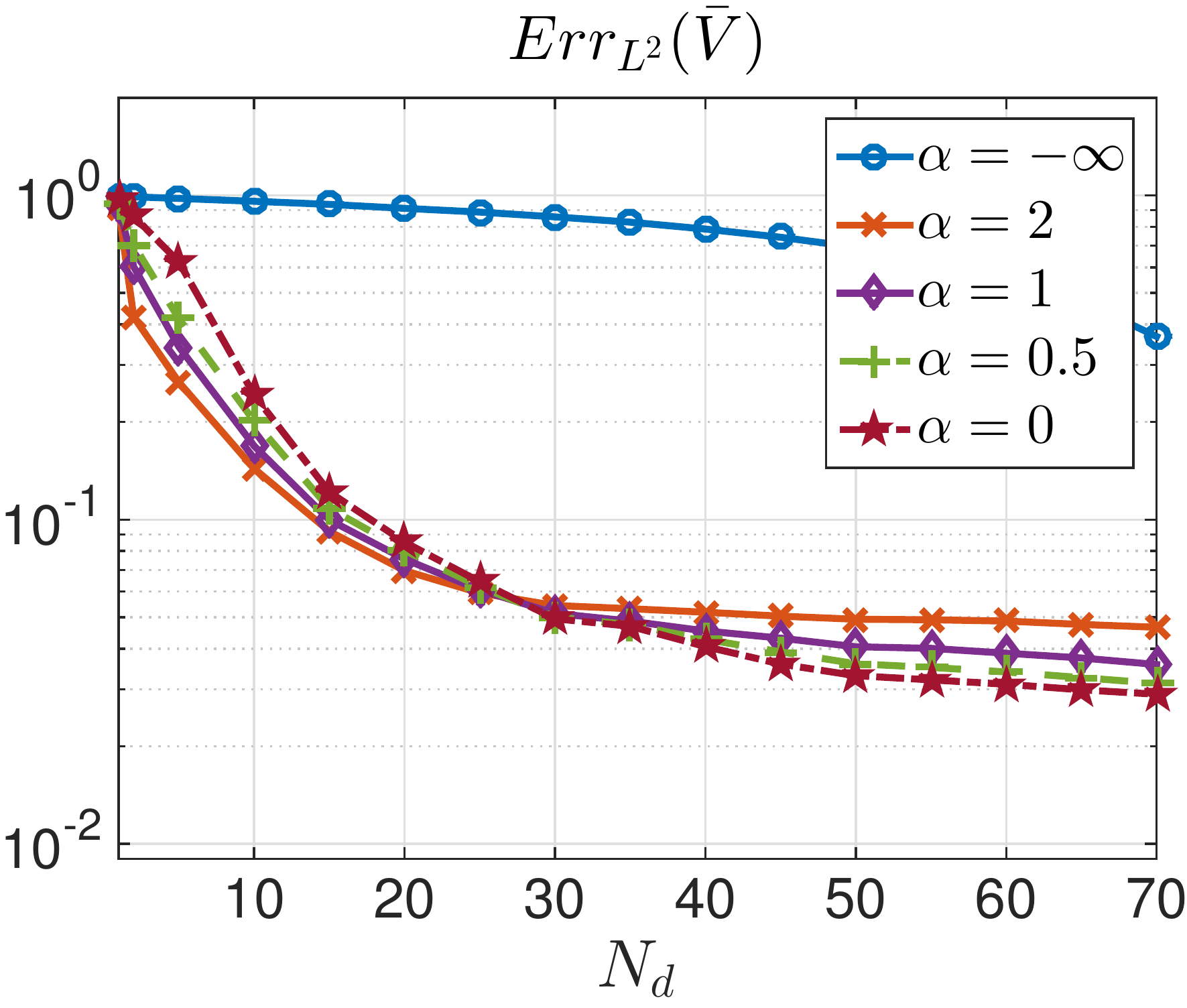}
	}
	\subfigure[$\lambda = 0.0005$]
	{
		\includegraphics[height=4cm,width=4cm]{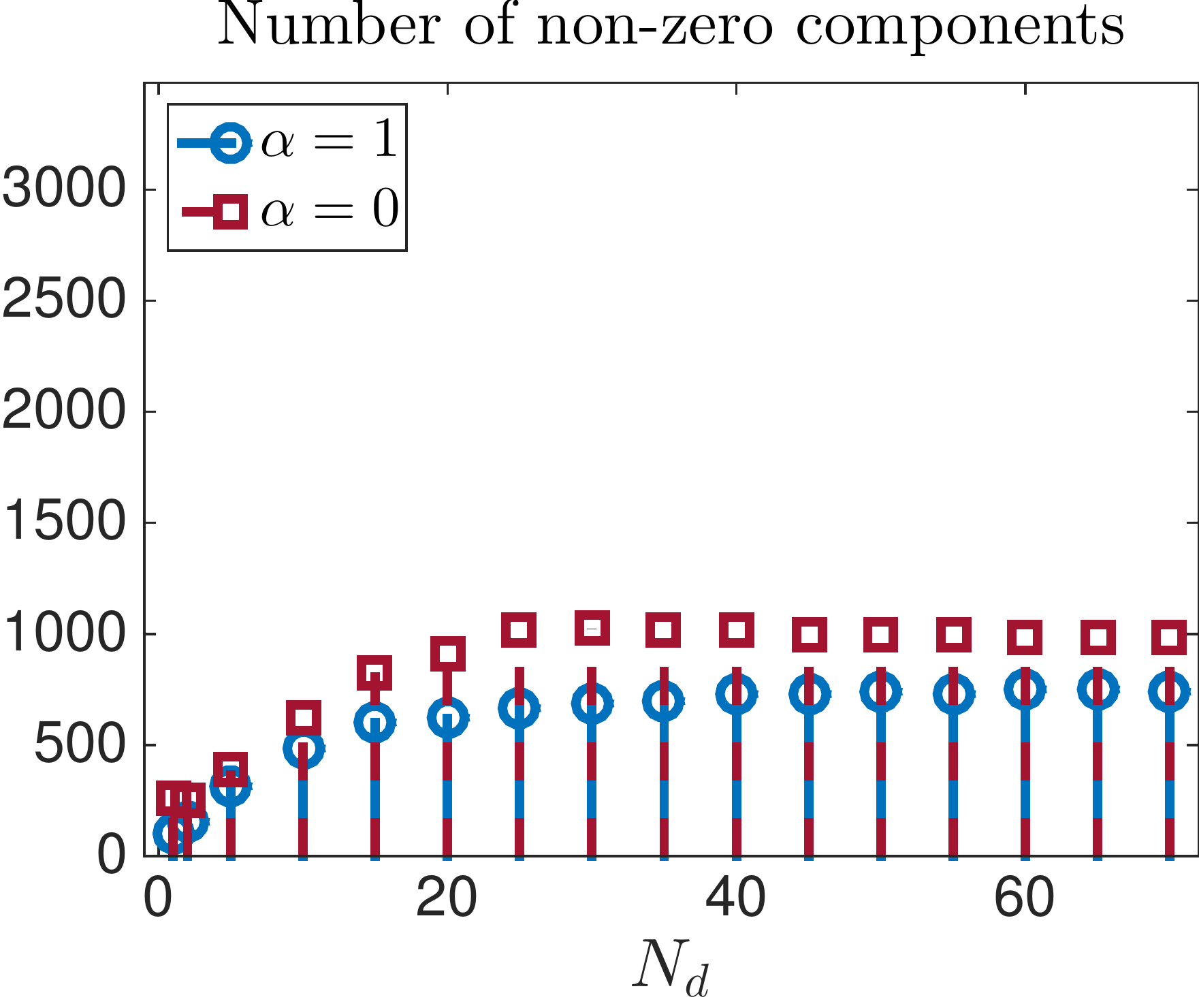}
	}
	\subfigure[$\lambda = 0.0001$]
	{
		\includegraphics[height=4cm,width=4cm]{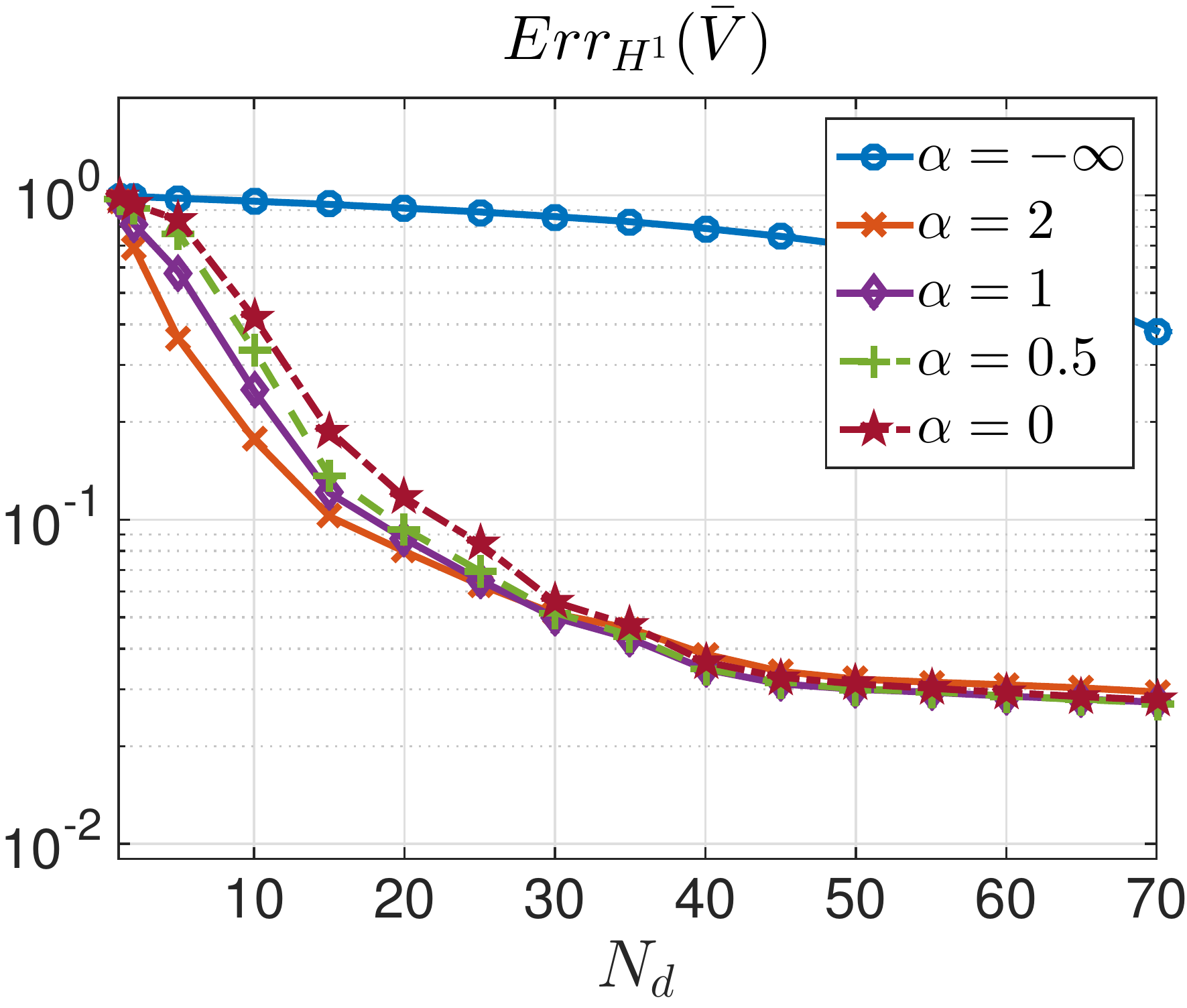}
	}
	\subfigure[$\lambda =0.0001$]
	{
		\includegraphics[height=4cm,width=4cm]{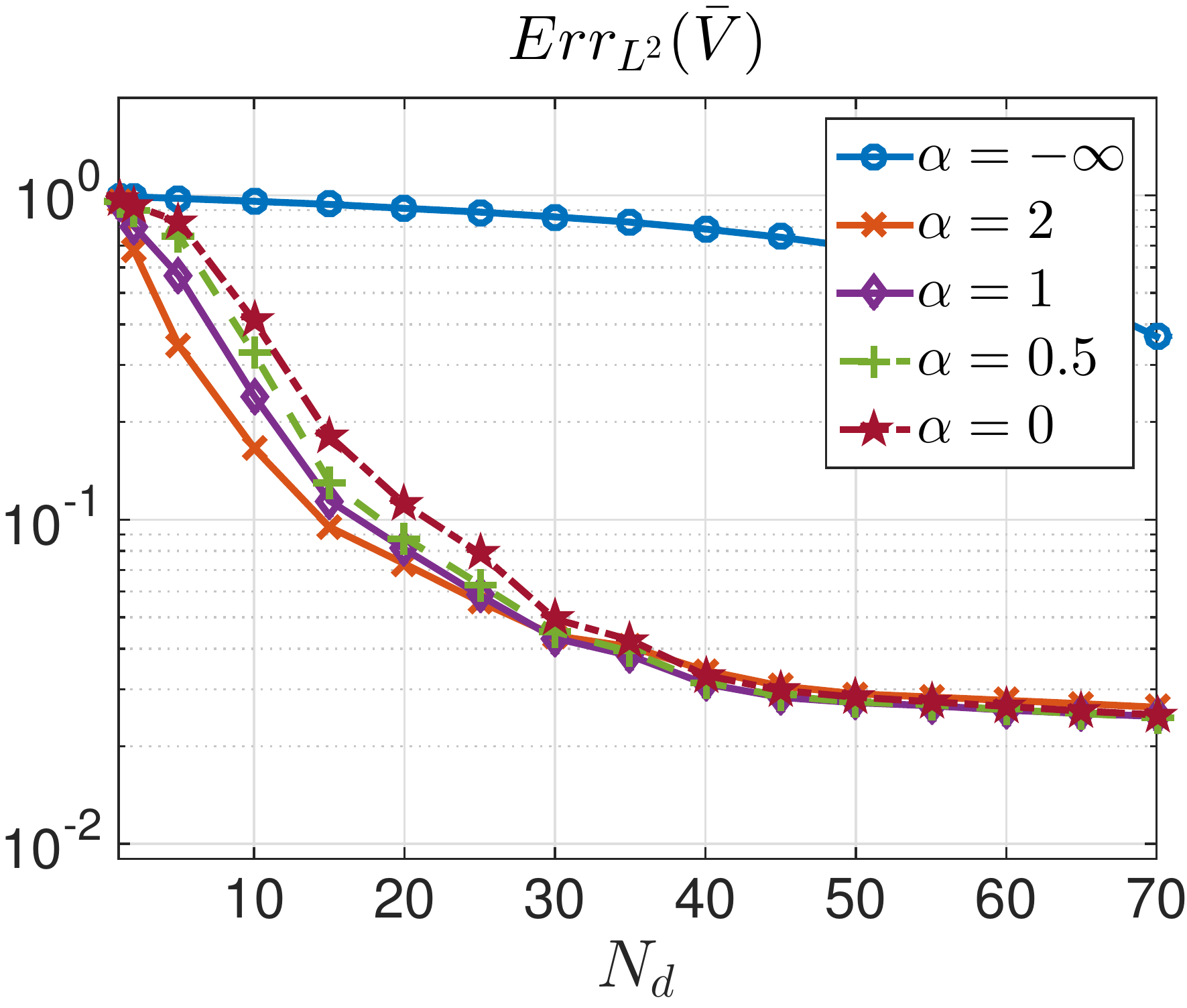}
	}
	\subfigure[$\lambda = 0.0001$]
	{
		\includegraphics[height=4cm,width=4cm]{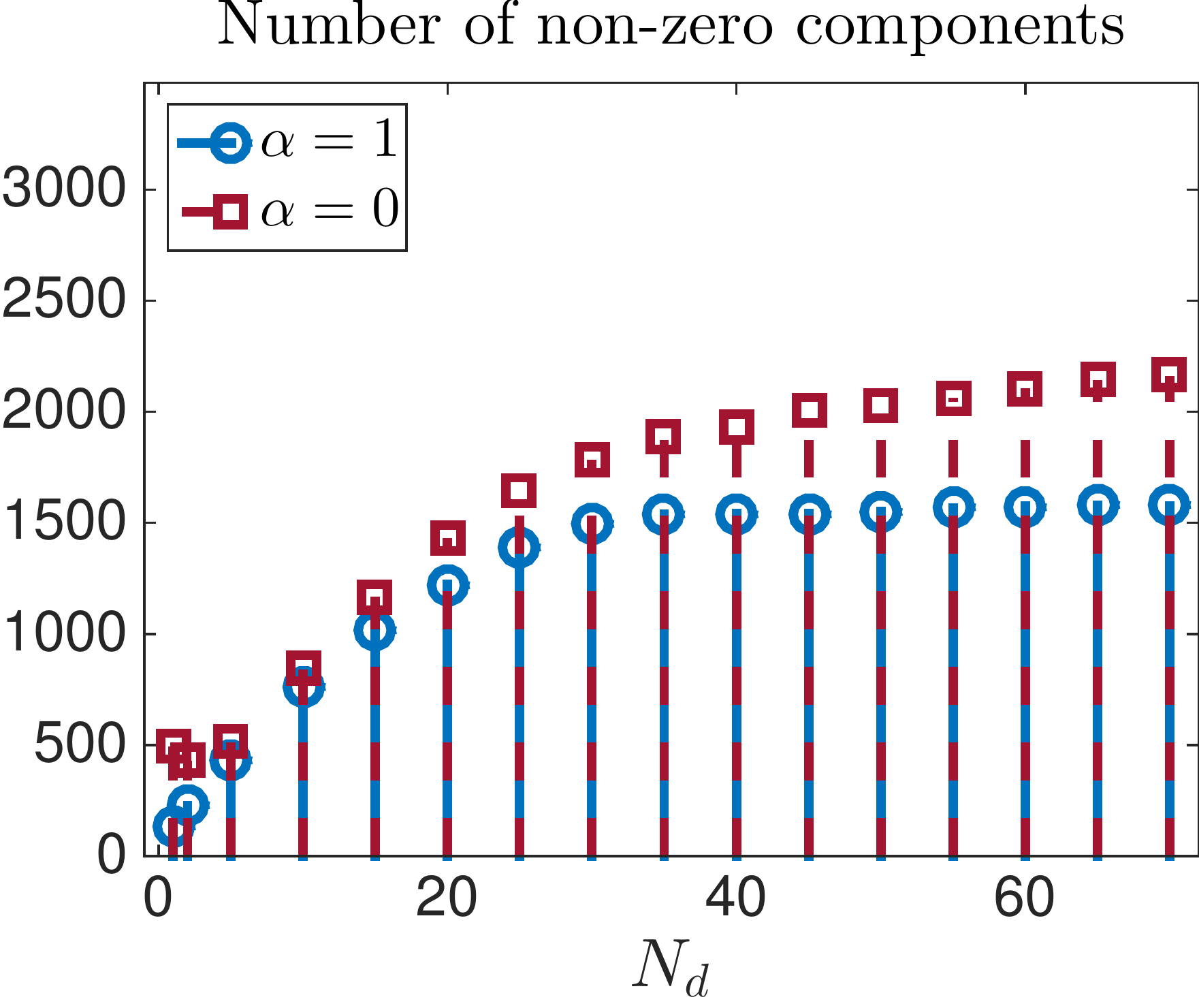}
	}
	\caption{Test 3. Numerical  results for the gradient-augmented polynomial approximation using the Legendre  polynomial basis. {\color{black} In this high-dimensional test with a feedback law in 80 dimensions, we observe that the use of gradient-augmented information is crucial to obtain small training errors with few samples.}}
	\label{Fig25}
\end{figure}

\begin{figure}[ht!]
	\centering
	\subfigure[$\lambda = 0.0005$]
	{
		\includegraphics[height=4cm,width=4cm]{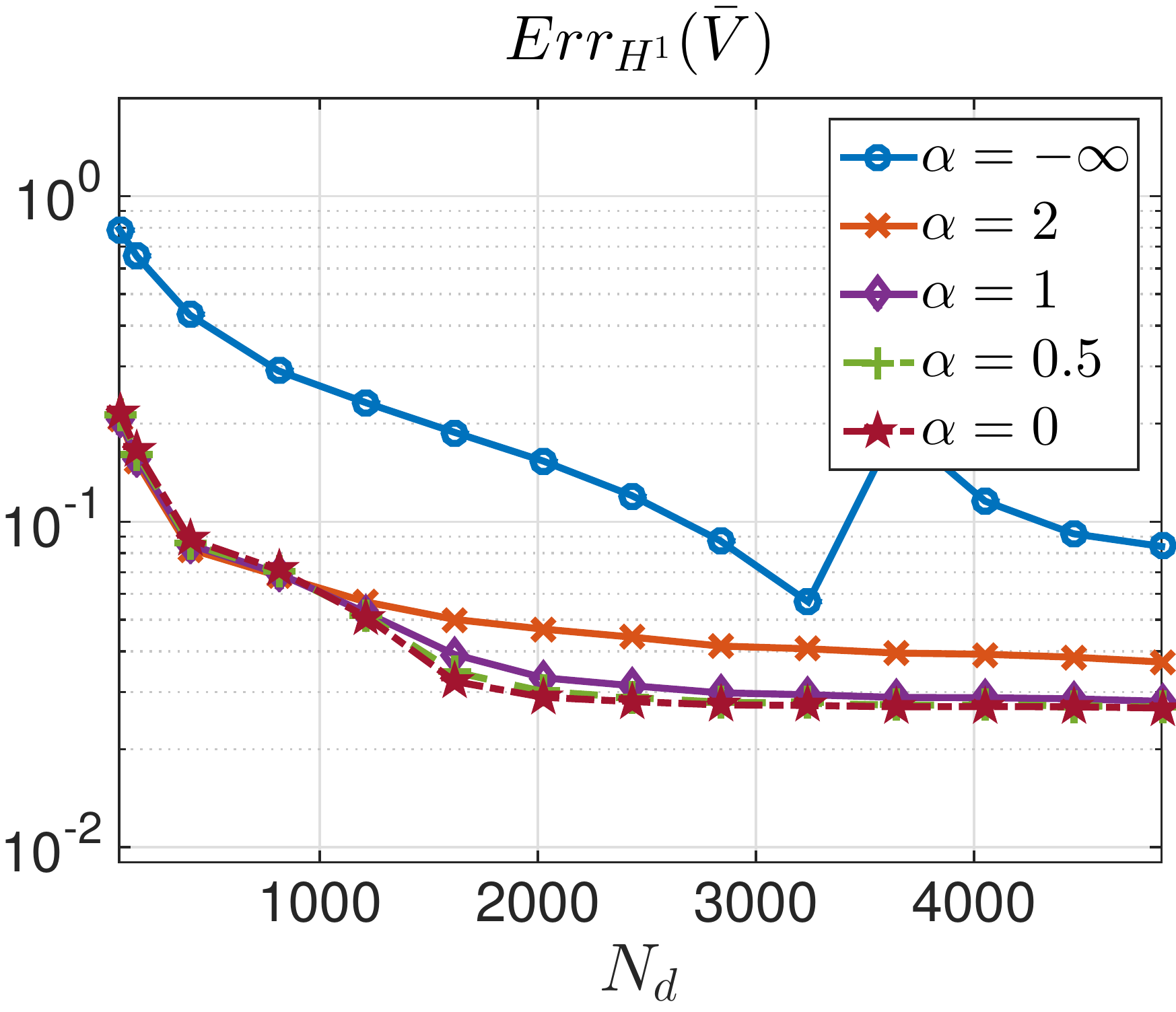}
	}
	\subfigure[$\lambda =0.0005$]
	{
		\includegraphics[height=4cm,width=4cm]{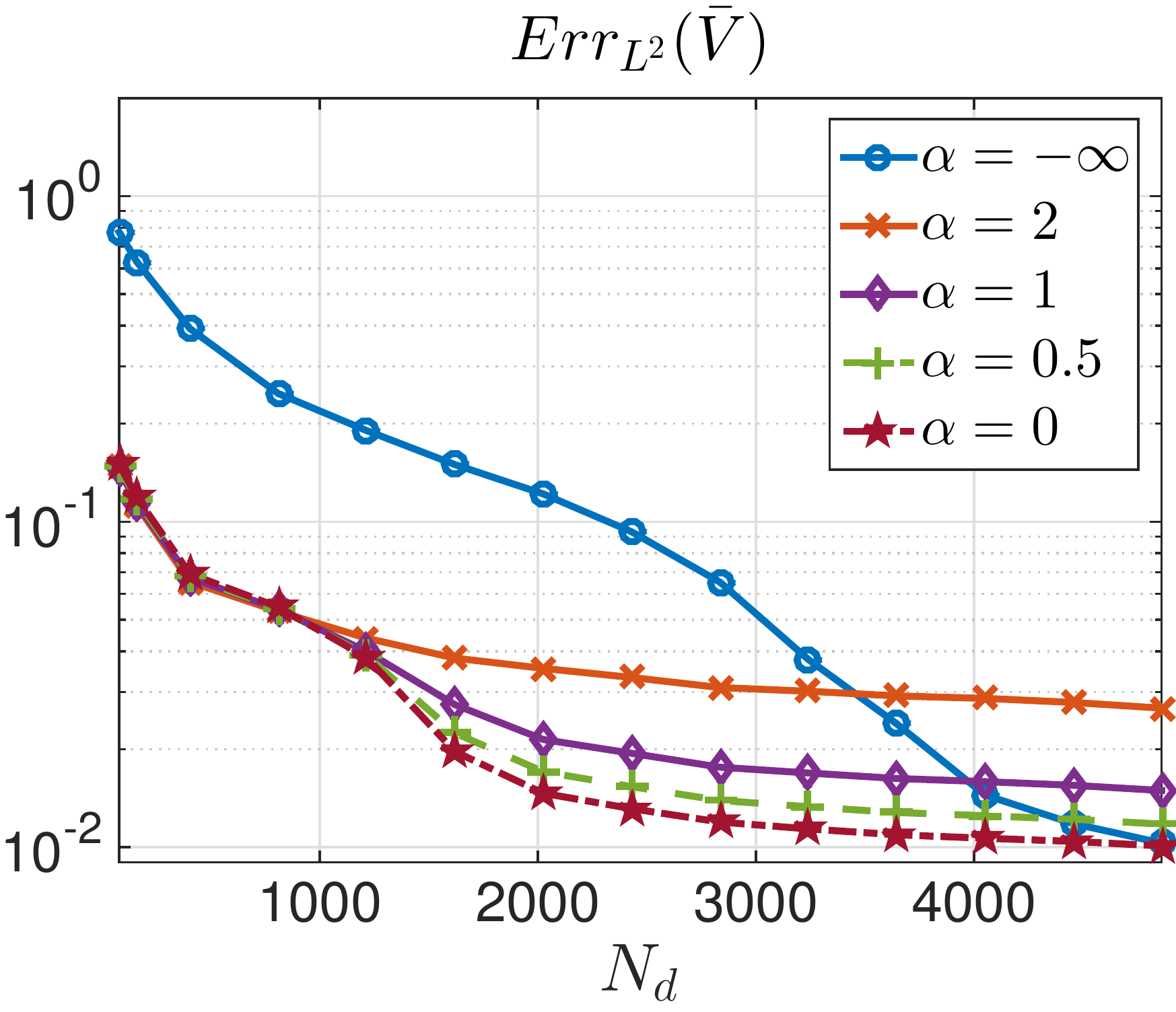}
	}
	\subfigure[$\lambda = 0.0005$]
	{
		\includegraphics[height=4cm,width=4cm]{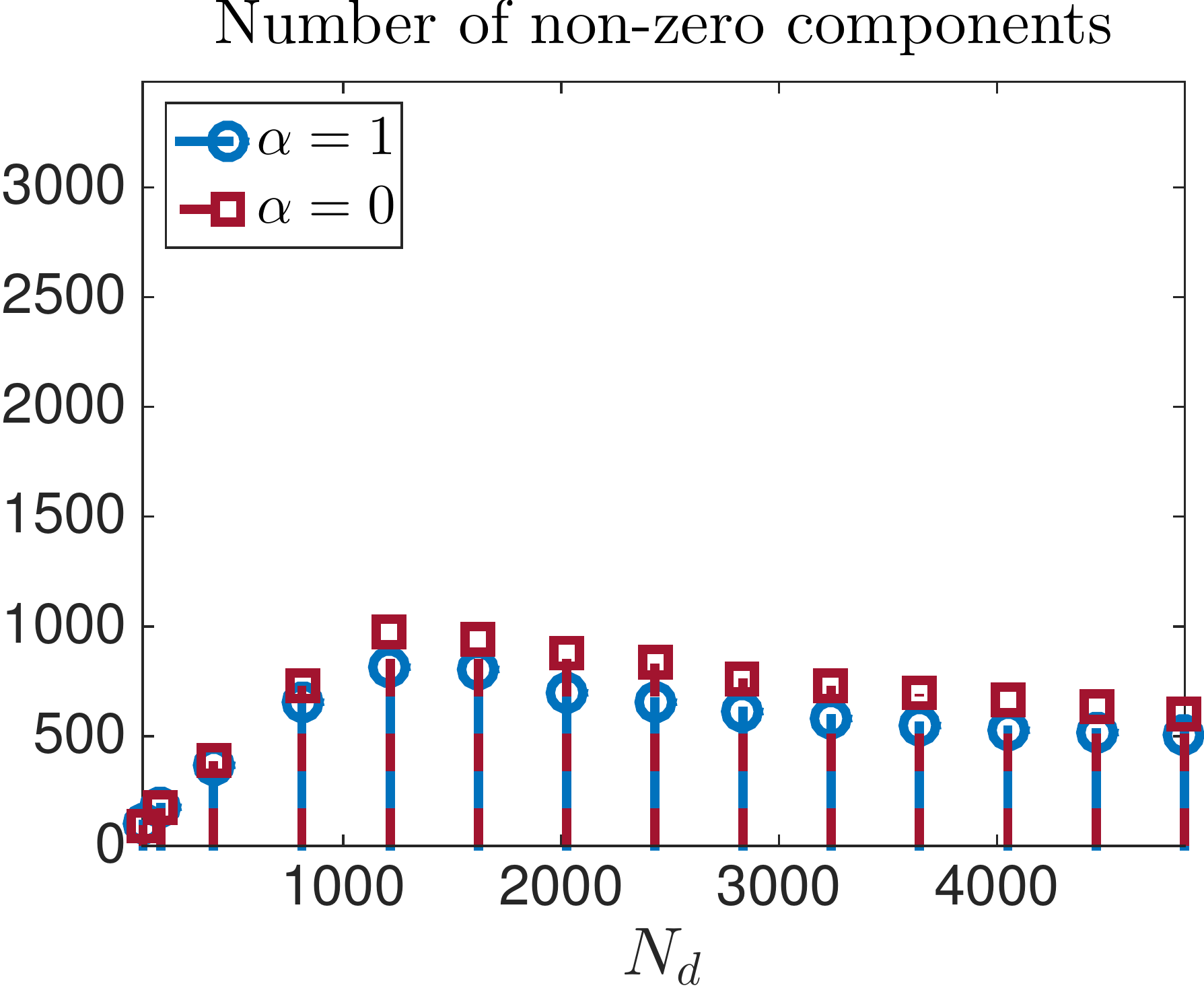}
	}
	
	\caption{Test 3. Numerical  results of the polynomial approximation using the Legendre basis  without gradient information. {\color{black} Compared to Figure \ref{Fig25}, the number of samples required to reach similar errors is almost 2 orders of magnitude larger, as this regression does not include gradient measurements.}}
	\label{Fig26}
\end{figure}

The choice $s = 4$ in the hyperbolic cross index set  $\mathfrak{I}(s)$, results in  $|\mathfrak{I}(4)|  =3481$ polynomial basis functions for $V$. As in the previous examples, we computed the solutions $\theta_{\ell_2}$, $\bar{\theta}_{\ell_2}$, $\theta_{\ell_1}$ and $\bar{\theta}_{\ell_1}$  for different values of $N_d$ and $\lambda$.  We report results for the errors and levels of non-sparsity for the case of Legendre polynomials as  basis functions, $\lambda =  5\times 10^{-4}, 10^{-4}$ and  $|\mathcal{I}_{val}| = 5000$, in Figures \ref{Fig25} and \ref{Fig26}, with and without gradient information, respectively.
For this 80-dimensional problem,   the  observations from the previous examples are confirmed, and the gradient-augmented sparse regression requires two orders of magnitude less of training samples to achieve the same error levels of the gradient-free counterpart.

\begin{figure}[ht!]
	\centering
	\subfigure[State]
	{
		
		\includegraphics[height=5cm,width=5cm]{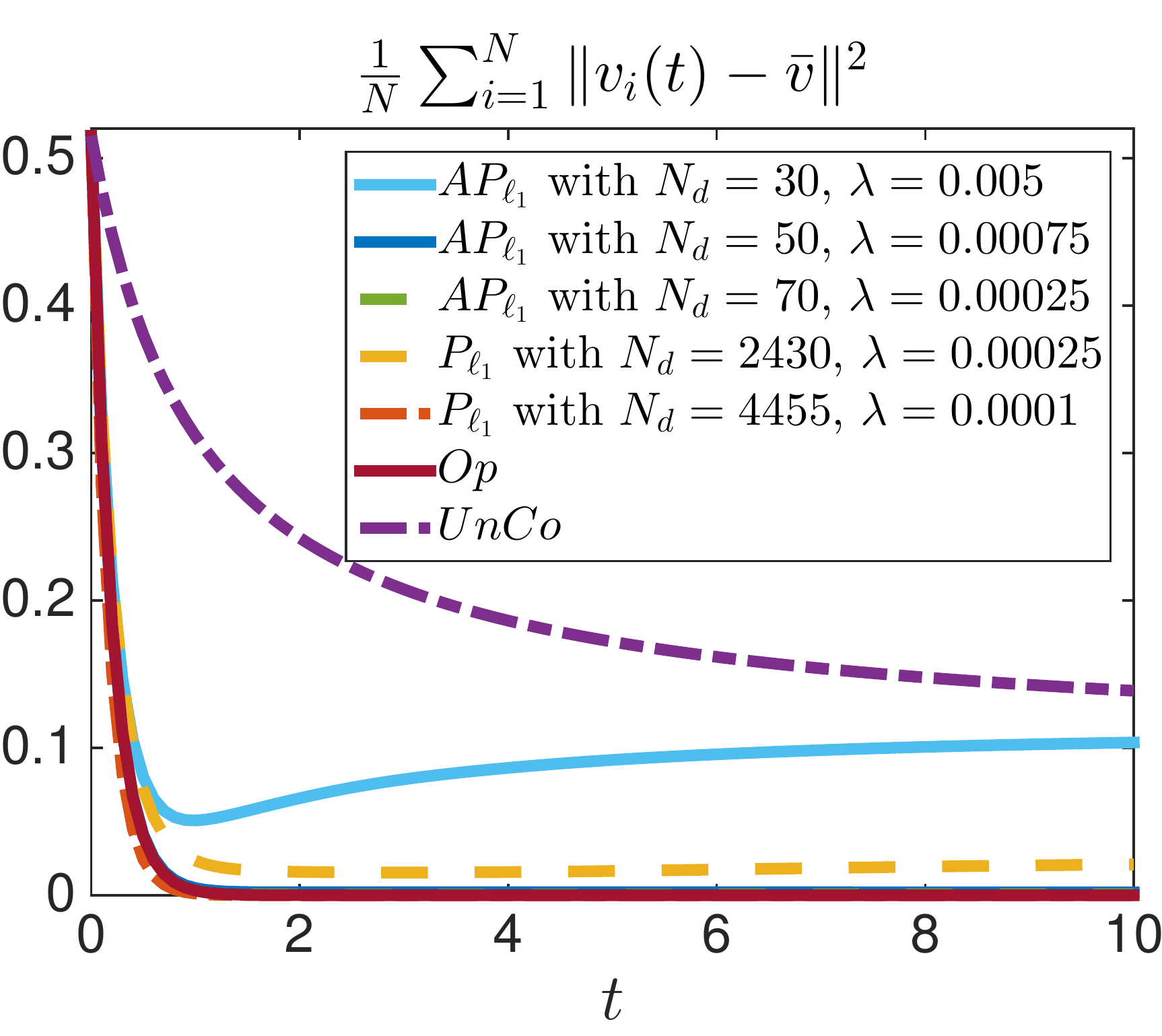}
	}
	\subfigure[Control]
	{
		\includegraphics[height=5cm,width=5cm]{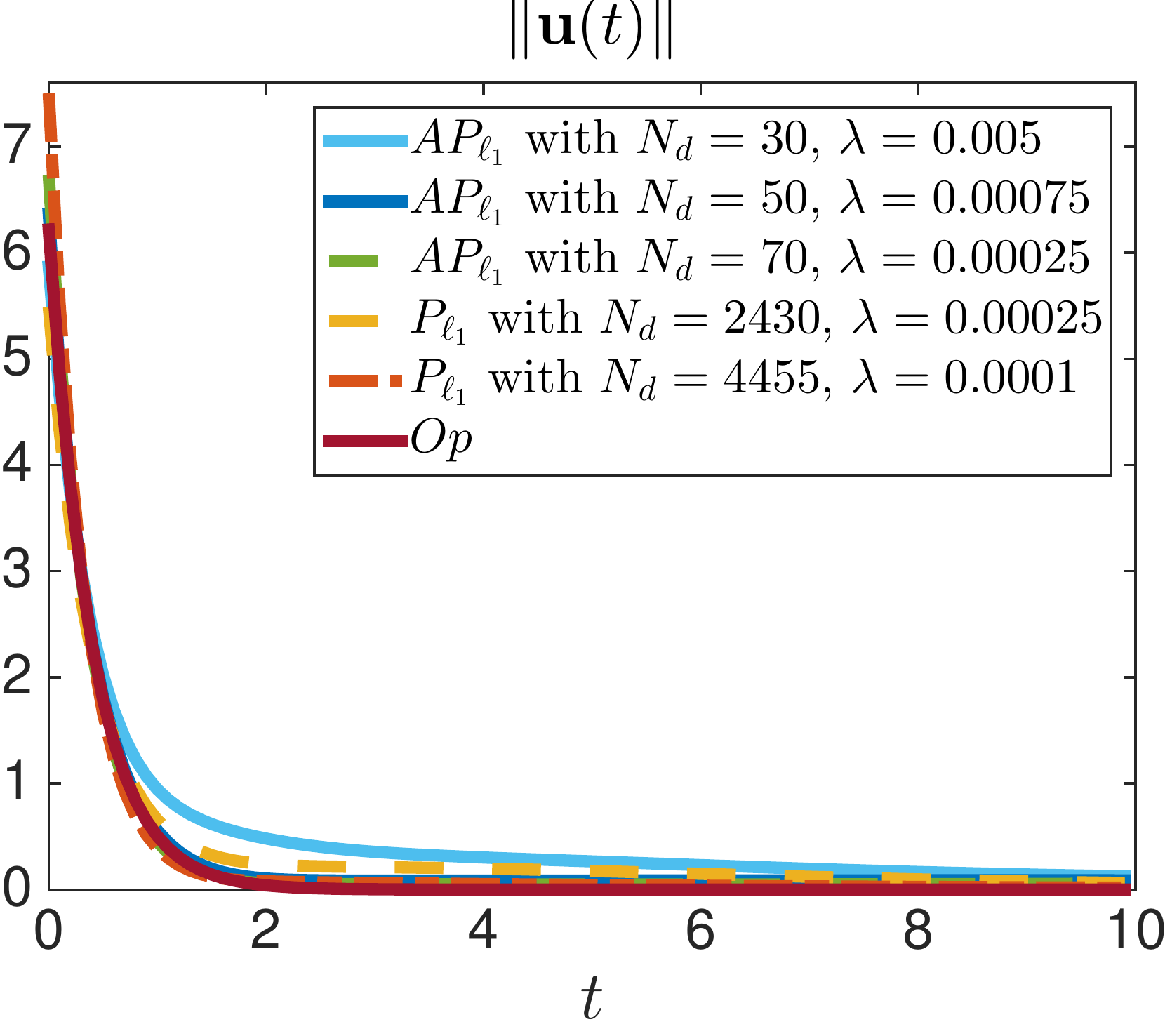}
	}
	\caption{Test 3. (a) Evolution of  the consensus variance $\frac{1}{N_a}\sum_{i=1}^{N_a}\|v_i(t)-\bar{v}\|^2$ for different control laws and (b) the norm   $\|\mathbf{u}(t)\|$ of the associated control signals. Here \textsl{Op} stands for the exact optimal control, and \textsl{UnCo} for the uncontrolled solution. {\color{black}Recovered feedback laws including gradient information and a small number of samples (50 and 70 in this case) effectively stabilize the dynamics around consensus, while a similar control law without gradient information requires over 4000 training samples.}}
	\label{Fig27}
\end{figure}

Furthermore, we computed the approximation  of  the optimal control according to \eqref{eq:feedback2} with $\bg =[\mathbf{0};I]^{\top}$ for different setting of  $\theta=\theta_{\ell_1}$ and $\theta=\bar{\theta}_{\ell_1}$, and present stabilization results in Figure \ref{Fig27}.

\begin{figure}[ht!]
	\centering
	\subfigure[Uncontrolled State]
	{
		\label{Fig30a}
		\includegraphics[height=3.5cm,width=3.5cm]{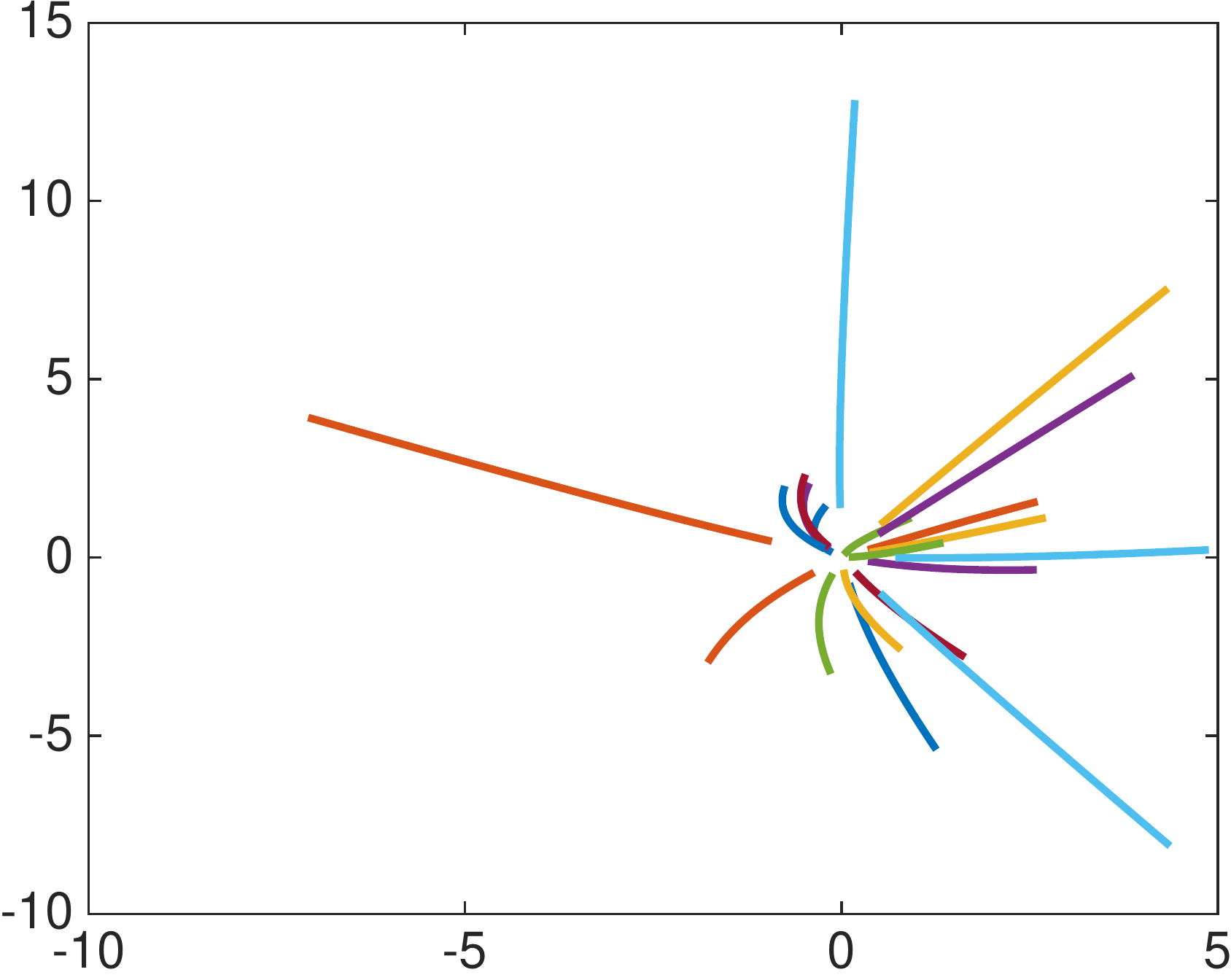}
	}
	\subfigure[Optimal state]
	{
		\label{Fig30}
		\includegraphics[height=3.5cm,width=3.5cm]{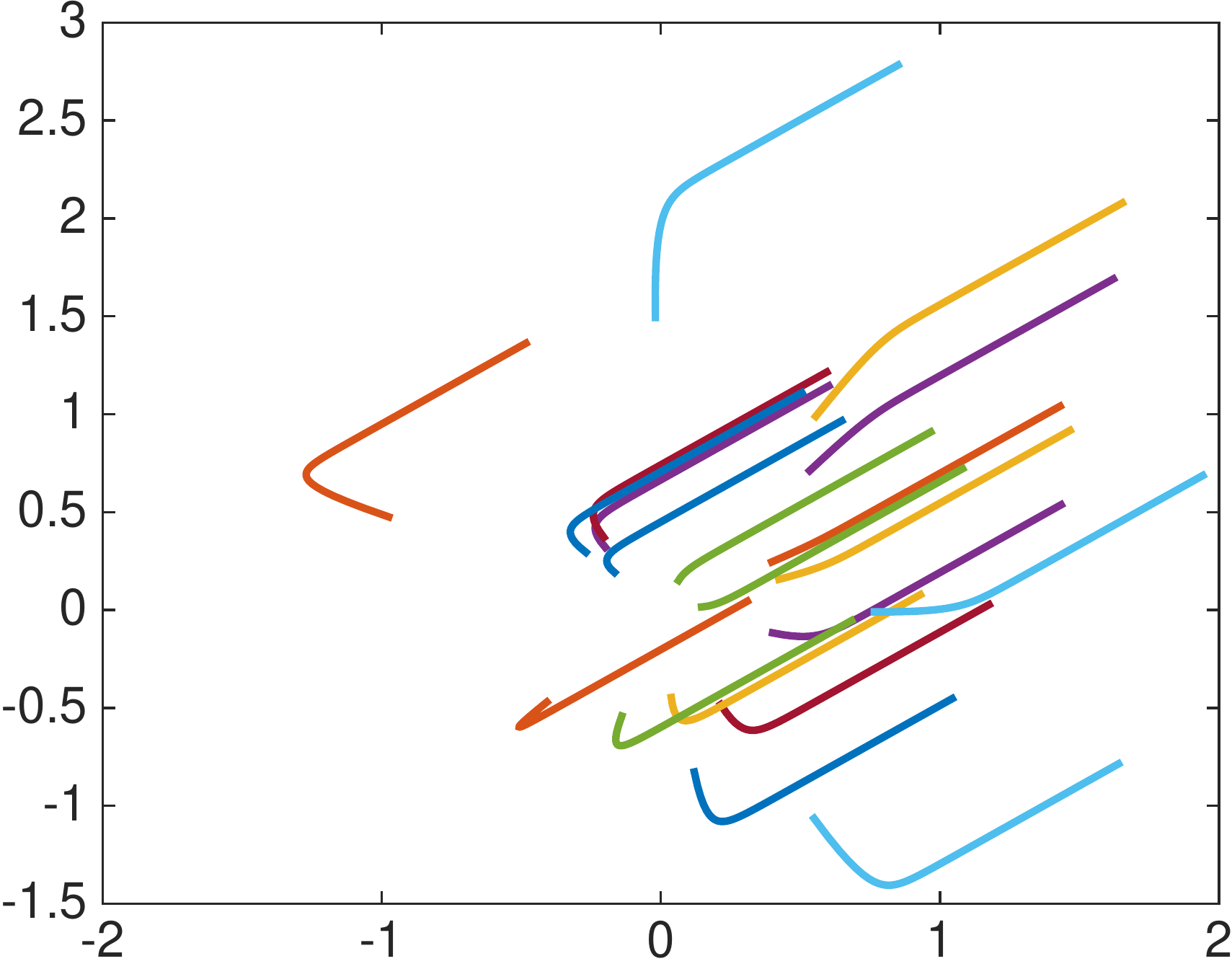}
	}
	\subfigure[$AP_{\ell_1}$ with $N_d =70$, $\lambda =2.5\times 10^{-4}$]
	{
		\label{Fig32}
		\includegraphics[height=3.5cm,width=3.5cm]{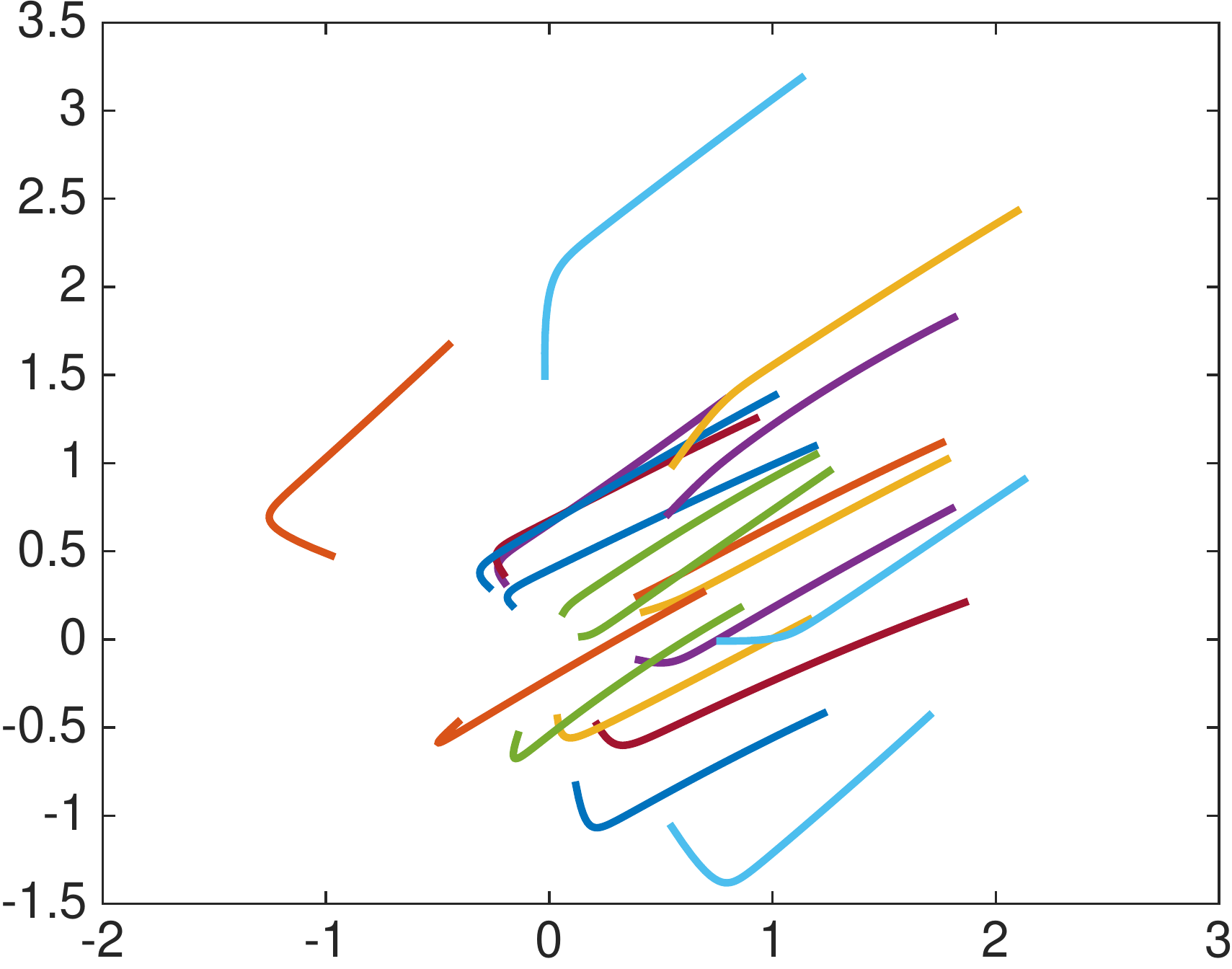}
	}
	\subfigure[$P_{\ell_1}$ with $N_d =2430$, $\lambda = 2.5\times 10^{-4}$]
	{
		\label{Fig33}
		\includegraphics[height=3.5cm,width=3.5cm]{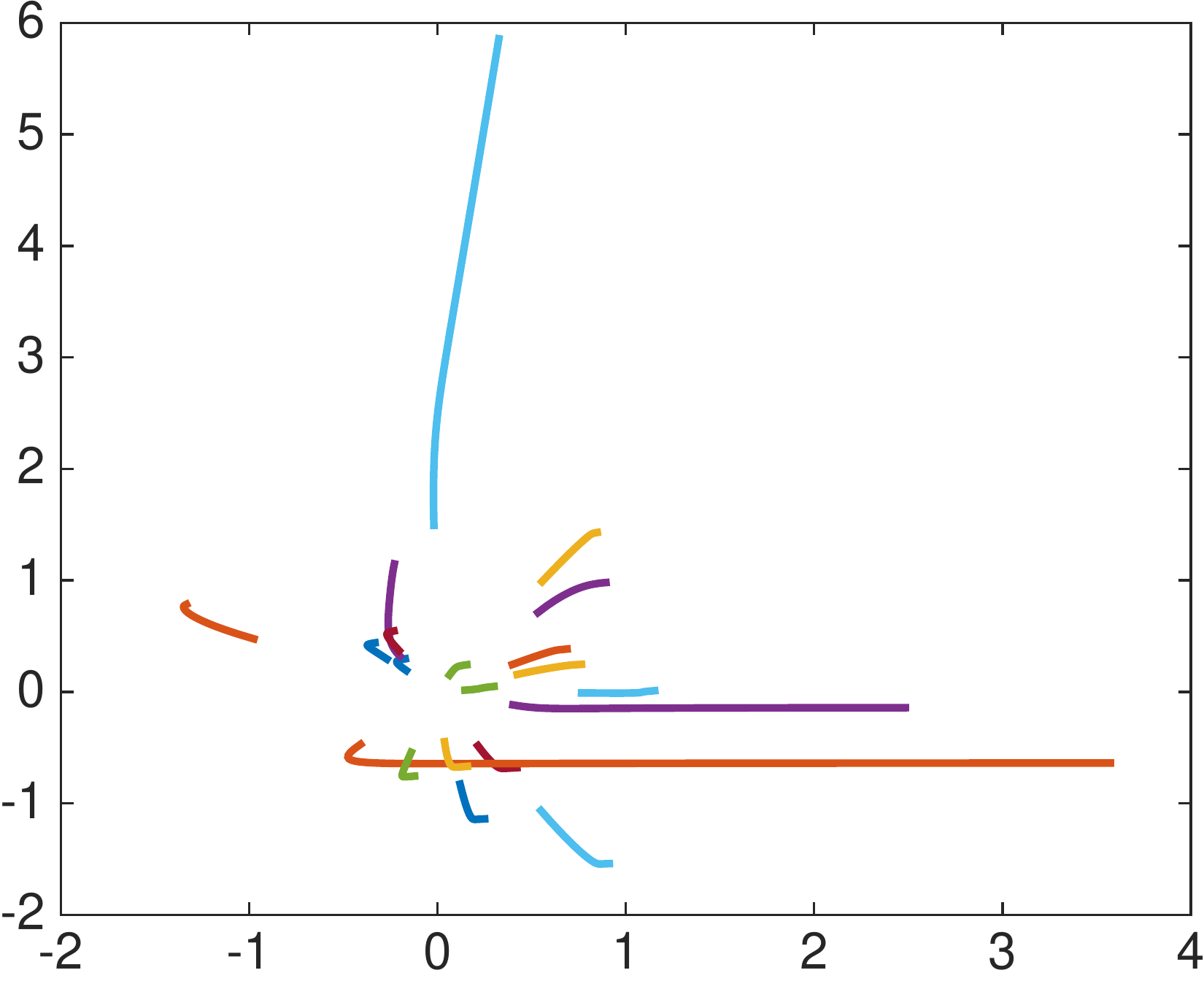}
	}
	
	\caption{Test 3. Trajectories generated by different  control laws. (a) Uncontrolled trajectories diverging in space. (b) Optimal trajectory for reference. In this high-dimensional problem ($n=80$), the sparse, gradient-augmented regression for $V$ (c) yields an feedback law which approaches the optimal trajectory (b) with few training samples $N_d=70$. {\color{black} In (d) we observe a controlled trajectory with a feedback law without gradient information and an insufficient number of samples, which fails to stabilize the dynamics.}}
\end{figure}

Figure \ref{Fig30a} shows the uncontrolled  dynamics of the agents for a specific  initial state $\hat\bx_0$.
The dynamics of the  optimal state and approximations, corresponding to $OC(\hat\bx_0)$, are plotted
in Figures \ref{Fig30}, \ref{Fig32}, and \ref{Fig33}. Colored trajectories are uncontrolled, and red trajectories represent the controlled evolution. By comparing these Figures, it can be seen that the dynamics of the optimal state  and its approximation obtained by $\mathbf{u}_{ \bar{\theta}_{\ell_1}}$ for $N_d =70$, $\lambda =2.5\times 10^{-4}$ are almost identical. Note that the controlled trajectories in these two subplots \ref{Fig30} and \ref{Fig32} achieve consensus, unlike \ref{Fig33}, where the feedback obtained with a sparse regression without gradient information does not stabilize the dynamics despite the large training dataset. This observation  is also supported by Table  \ref{table8} and Figure \eqref{Fig27}. Table \ref{table8} shows that the smallest validation error is achieved for the gradient-augmented sparse regression of the value function with only 70 training samples, with approximately 20\% of nonzero components. Figure \eqref{Fig27} depicts the evolution of the tracking term $\frac{1}{N_a}\sum_{i=1}^{N_a}\| v_i(t)-\bar v \|^2$ and the norm of the control of every feedback law. From this figure, we can see that  the control $\mathbf{u}_{ \bar{\theta}_{\ell_1}}$ associated to $N_d =70$, $\lambda =2.5\times 10^{-4}$ delivers the best approximation for the optimal control of $OC(\hat\bx_0)$ among the different control laws.

\begin{table}[htbp!]
	\begin{center}
		\begin{tabular}{cccc}
			\hline\hline
			&  $Err_{L^2} $& $Err_{H^1} $&  Nonzero components  \\
			\cmidrule(lr){2-2}\cmidrule(lr){3-3}\cmidrule(lr){4-4} \morecmidrules\cmidrule(lr){2-2}\cmidrule(lr){3-3}\cmidrule(lr){4-4}
			$\bar{V}_{\ell_1}$ for  $N_d =50$, $\lambda = 7.5\times10^{-4}$ & $5.40 \times 10^{-2}$ &$6.37 \times 10^{-2}$&$393/3481$  \\
			$\bar{V}_{\ell_1}$ for $N_d =70$, $\lambda =2.5\times10^{-4}$ &$3.56\times 10^{-2}$ &$4.11\times10^{-2}$&$738/3481$ \\
			$V_{\ell_1}$ for $N_d =2430$, $\lambda = 2.5\times10^{-4}$  &$7.46 \times 10^{-2}$ &$9.38\times10^{-2}$&$656 / 3481$ \\
			\hline\hline
		\end{tabular}
	\end{center}
	\caption{Test 3. Validation errors for different regressions. {\color{black} The use of gradient-augmented information for the regression leads to the recovery of high-dimensional feedback laws with a reduced number of training samples. The inclusion of a sparsity promoting term in the regression reduces to the number of non-zero components in the control law to less than 15\%.}}
	\label{table8}
\end{table}

\subsection*{Concluding remarks}

We have presented a sparse polynomial regression framework for the approximation of feedback laws arising in nonlinear optimal control. The main ingredients of our approach include: the generation of a gradient-augmented dataset for the value function associated to the control problem by means of PMP solves, a hyperbolic cross polynomial ansatz for recovering the value function and its feedback law, and a sparse optimization method to fit the model. Through a series of numerical tests, we have shown that the proposed approach can approximate high-dimensional control problems at moderate computational cost. {\color{black} It is worth to note that our numerical tests correspond to smooth value functions, where gradients are regular and their inclusion in the regression is well-posed}. The gradient-augmented dataset reduces the number of open-loop solves required to recover the optimal control, and the sparse regression provides a feedback law of reduced complexity, which is an appealing feature for real-time implementations. The effectiveness of the proposed methodology suggests different research directions:
\begin{itemize}
	\item The deep neural network ansatz proposed in ~\citep{nakazim} can be combined with a sparsity-promoting loss function along the lines of our work. In the light of recent results discussed in ~\citep{adcockNN}, it is a pertinent question to find whether deep neural networks or polynomial approximants are more effective ansatz for the value function. As we have previously mentioned, there are different control-theoretical arguments which support the case for having a polynomial approximation of the value function. 
	\item {\color{black}The study of the regression framework in control problems with a fully nonlinear control structure. For the sake of simplicity, we have restricted our problem formulation to the control-affine case, however, the link between PMP and HJB holds in the fully nonlinear case under standard regularity assumptions ~\citep[Section 3.2]{Subbotina}. Nevertheless, the nonlinear control structure can deteriorate the approximation power of the hyperbolic cross ansatz for the value function. Moreover, dataset generation in the fully nonlinear case can be costly as the structure of \eqref{eq:tpbvp} becomes more complex.}
	\item The extension of the presented results in the context to time-dependent, and second-order stochastic control problems where the representation formula given by the PMP is replaced by a backward stochastic differential equation ~\citep{gprbsde}. Finally, the ideas proposed in this work regarding sparse polynomial regression can be implemented in the context of approximate dynamic programming in reinforcement learning ~\citep{bertsekasrl}.
\end{itemize}

% Acknowledgements should go at the end, before appendices and references

%\acks{Dante Kalise  was supported by a public  grant  as  part  of the Investissement d'avenir project, reference ANR-11-LABX-0056-LMH, LabEx LMH. Karl Kunisch was supported in part by the ERC advanced grant 668998 (OCLOC) under the EU’s H2020 research program.}


\vskip 0.2in
\bibliography{polynomshort}

\end{document}